\documentclass[12pt, a4paper]{amsart}
\usepackage{url}
\usepackage{t1enc}
\usepackage{ifthen,epsfig,array}
\usepackage{amsmath,amsfonts,amsthm,amssymb}
\usepackage[mathscr]{eucal}
\usepackage{latexsym}
\usepackage[T1]{fontenc}
\usepackage[latin1]{inputenc}
\usepackage[english]{babel}
\usepackage{fancyhdr}
\usepackage{verbatim}
\usepackage{graphicx}
\usepackage{epstopdf}
\usepackage{longtable}
\usepackage{fancybox}
\usepackage{float}
\usepackage{geometry}
\usepackage{enumerate}
\usepackage{bm}
\usepackage{color}
\usepackage{dcolumn}
\usepackage{slashbox}
\usepackage[noend]{algorithmic} 
\usepackage{algorithm, caption}
\algsetup{indent=1.5em} 
\usepackage{subfig}

\theoremstyle{myplain}
\newtheorem{theorem}{Theorem}[section]

\theoremstyle{myremark}

\newcommand{\N}{\mathbb{N}}

\newcommand{\R}{\mathbb{R}}

\newcommand{\Vjk}{\frac{X_j}{t}}
\newcommand{\Ij}{1_{(X_j>t)}}

\newcommand{\E}{\mathbb{E}}
\newcommand{\Prob}{\mathbb{P}}

\newcommand{\talpha}{\tilde{\alpha}}
\newcommand{\tlambda}{\tilde{\lambda}}
\newcommand{\ttau}{\tilde{\tau}}
\newcommand{\halpha}{\widehat{\alpha}}
\newcommand{\betain}{\beta_{\infty}}

\newcommand{\overF}{\overline{F}}

\DeclareMathOperator*{\argmin}{arg\,min}
\DeclareMathOperator*{\argmax}{arg\,max}

\usepackage{tabularx,ragged2e,booktabs,caption}

\makeatletter
\@namedef{subjclassname@1991}{Mathematics Subject Classification}
\makeatother

\makeatletter
\newcommand{\addresseshere}{%
  \enddoc@text\let\enddoc@text\relax\vfill{0}
}
\makeatother

\begin{document}
\title[Tempered Pareto-type modelling using Weibull distributions]{Tempered Pareto-type modelling using Weibull distributions}

\author[H. Albrecher]{Hansj\"{o}rg Albrecher}
\address{(H. Albrecher) Department of Actuarial Science, Faculty of Business and Economics, University of Lausanne, Switzerland and Swiss Finance Institute}
\email{hansjoerg.albrecher@unil.ch}

\author[J. Araujo-Acuna]{Jos\'{e} Carlos Araujo-Acuna}
\address{(JC. Araujo-Acuna) Department of Actuarial Science, Faculty of Business and Economics, University of Lausanne, Switzerland}
\email{josecarlos.araujoacuna@unil.ch}

\author[J. Beirlant]{Jan Beirlant}
\address{(J. Beirlant) Department of Mathematics, LStat and LRisk, KU Leuven, Belgium and Department of Mathematical Statistics and Actuarial Science, University of the Free State, South Africa}
\email{jan.beirlant@kuleuven.be}

\begin{abstract}
In various applications of heavy-tail modelling, the assumed Pareto behavior is tempered ultimately in the range of the largest data. In insurance applications, claim payments are influenced by claim management and  claims may for instance be subject to a higher level of inspection at highest damage levels leading to weaker tails than apparent from modal claims. 
Generalizing earlier results of Meerschaert et al.\ \cite{ Meerschaert2012parameter} and Raschke \cite{raschke2020alternative}, in this paper we consider tempering of a Pareto-type distribution with a general Weibull distribution in a peaks-over-threshold approach. This requires to modulate the tempering parameters as a function of the chosen threshold. Modelling such a tempering effect is important in order to avoid overestimation of risk measures such as the Value-at-Risk ($VaR$) at high quantiles. We use a pseudo maximum likelihood approach to estimate the model parameters, and consider the estimation of extreme quantiles. We derive basic asymptotic results for the estimators, give illustrations with simulation experiments and apply the developed techniques to fire and liability insurance data, providing insight into the relevance of the tempering component in heavy-tail modelling.  
\end{abstract}

\maketitle


\section{Introduction}
Probability distributions with power-law tails are extensively used in various fields of applications including insurance, finance, information technology, mining of precious stones and language studies (see e.g.\ \cite{nair} for a recent overview). In extreme value methodology such applications are appropriately modelled  using the concept of Pareto-type models such that a variable $X$ of interest satisfies 
\begin{eqnarray}
\Prob(X > x) = x^{-\alpha}\ell(x),
\label{Patype}
\end{eqnarray}
with $\alpha >0$ and some slowly varying function $\ell$ satisfying 
\begin{eqnarray}
\frac{\ell(tx)}{\ell(t)} \rightarrow 1 \text{ as }  t\rightarrow\infty  \text{ for every } x>0. 
\label{ell}
\end{eqnarray}
In addition to the (pure) Pareto distribution, further examples from this model are the Burr, Fr\'echet, $t$ and log-gamma distribution (see  Beirlant et al.\  \cite[Ch. 2]{beirlant1996practical} for an overview). Often the power-law behaviour does not extend indefinitely due to some truncation or tapering effects. In Beirlant et al.\ \cite{beirlant2016tail}, estimation of truncated tails was developed in a peaks-over-threshold (POT) approach for Pareto-type tails, and other max-domains of attraction were dealt with in Beirlant et al.\ \cite{beirlant2017fitting}. Inspired by applications in geophysics and finance, Meerschaert et al.\ \cite{ Meerschaert2012parameter} discussed parameter estimation under exponential tempering of a simple Pareto law  with survival function
\begin{eqnarray}
\Prob(X > x) = c x^{-\alpha}e^{-\beta x},
\label{expPa}
\end{eqnarray} 
where  $\alpha,\beta >0$ and $c>0$ is a scale parameter.
 In the context of insurance data, Raschke \cite{raschke2020alternative} recently discussed the use of the more general Weibull tempering of a simple power law  with survival function
\begin{eqnarray}
\Prob(X > x) = c x^{-\alpha}e^{-(\beta x)^{\tau}},
\label{weiPa}
\end{eqnarray}
with $c, \alpha,\beta, \tau >0$.\\

\noindent
However, typically the power-law behaviour only sets in from some threshold $t$ on, rather than from the lowest measurements as assumed when using the simple Pareto model. The Pareto-type model \eqref{Patype} allows for flexible modelling of this behaviour. In this paper we therefore want to study {\it Weibull tempered Pareto-type distributions} with survival function
\begin{eqnarray}
\Prob(X > x) =  x^{-\alpha} \ell (x) e^{-(\beta x)^{\tau}},
\label{WtempPa}
\end{eqnarray}
with $\ell$ a slowly varying function, $\alpha=1/\gamma >0$ controlling the power-law tail with extreme value index $ \gamma$, and $\beta,\tau$ governing the Weibull tempering.\\

\begin{figure}[ht]
\centering
\centering
	\subfloat{\includegraphics[width=0.5\textwidth]{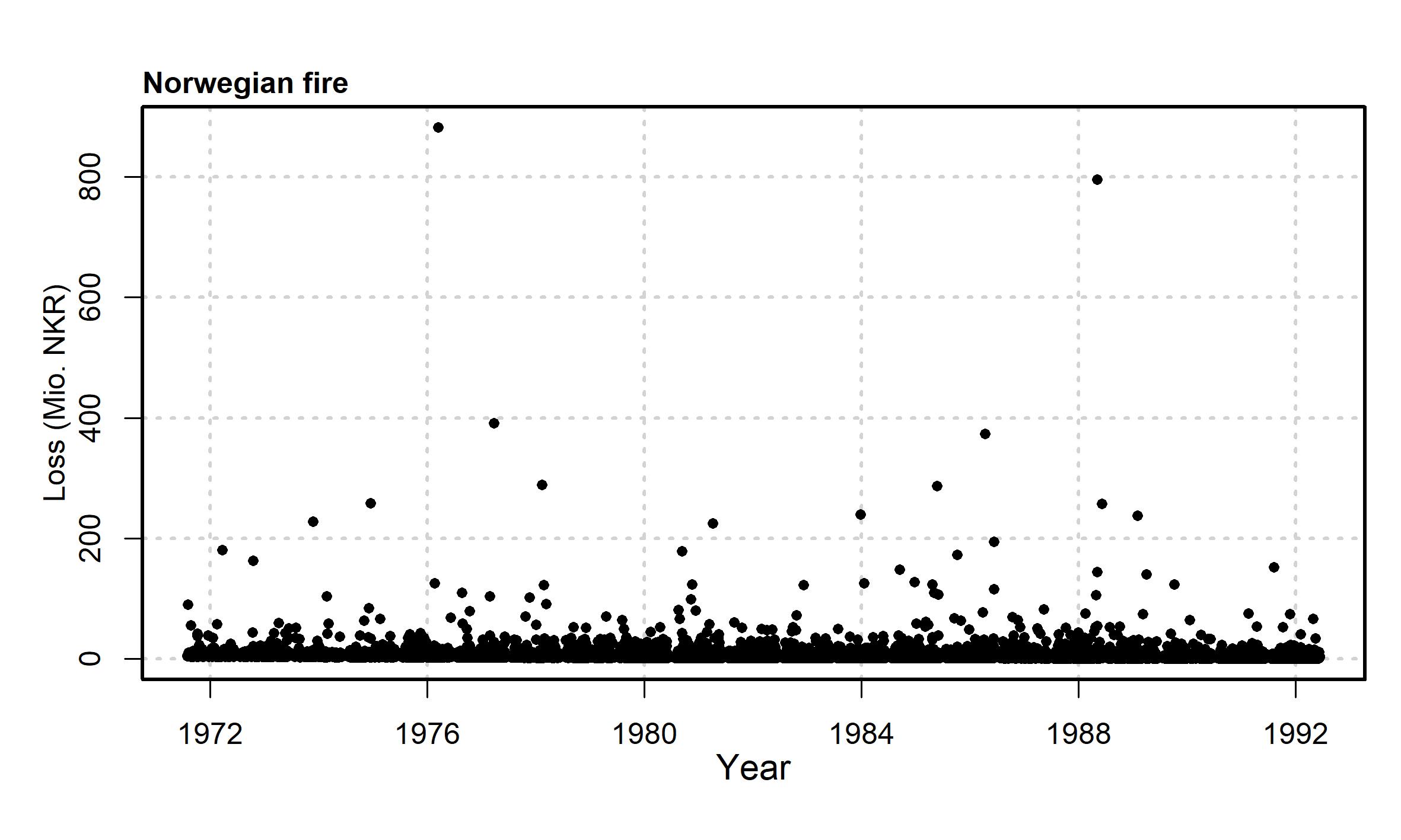}}
	\subfloat{\includegraphics[width=0.5\textwidth]{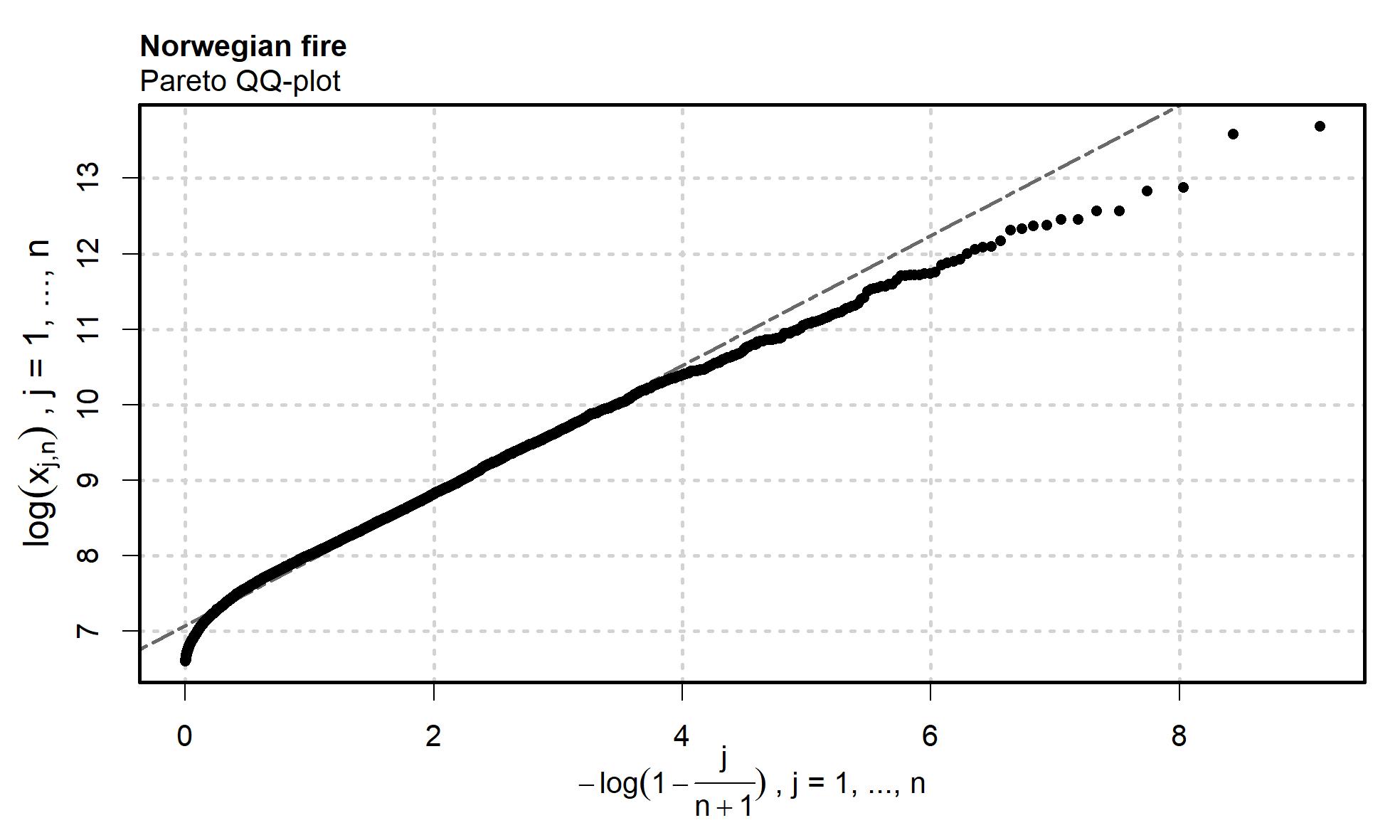}} \\
	\centering\subfloat{\includegraphics[width=0.5\textwidth]{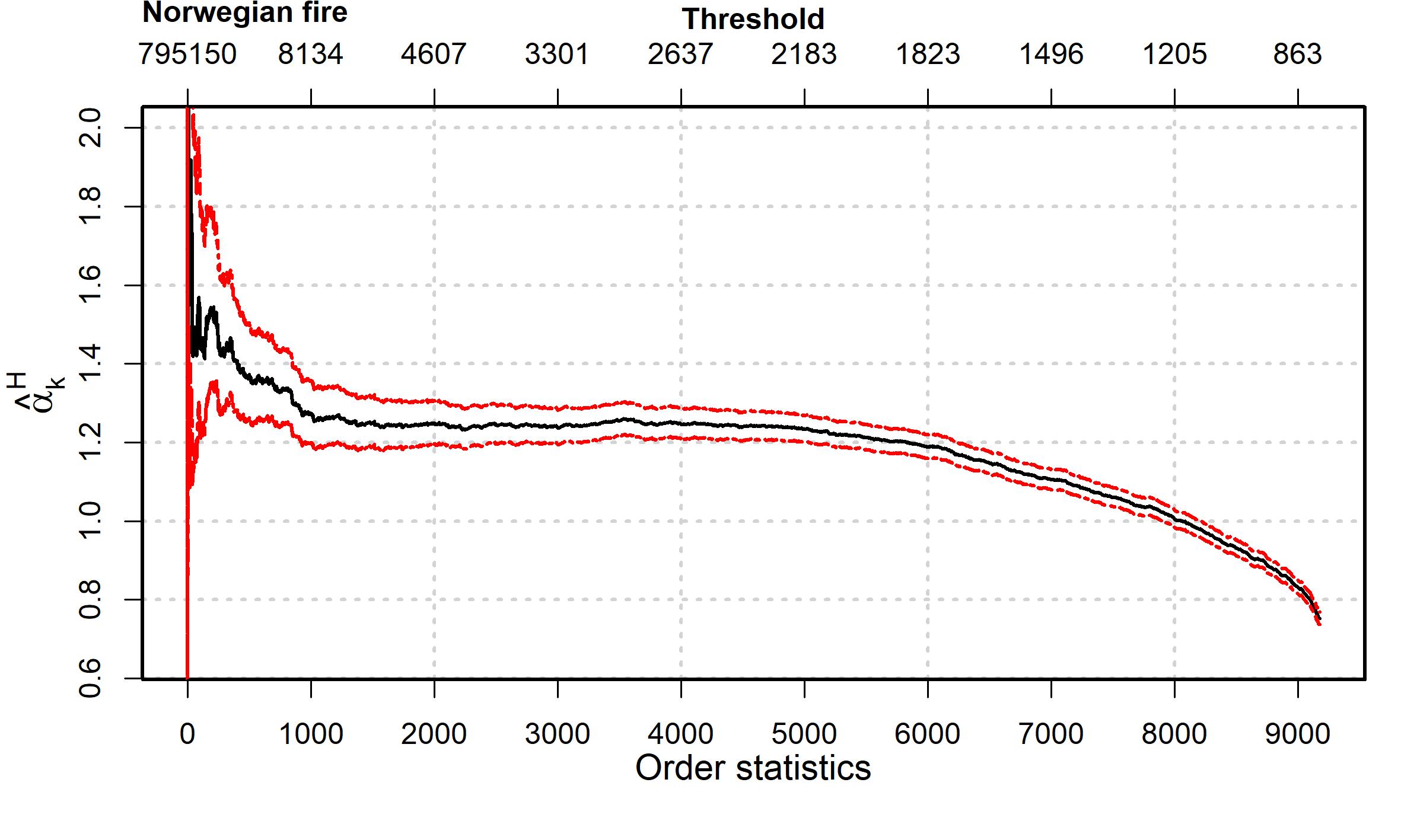}}\\
	\caption{Norwegian Fire claim data: claim sizes as a function of occurrence time (top left), log-log plot (top right) and $\hat{\alpha}^H_k$ estimates with 95\% confidence interval  (bottom).}
	\label{fig_Norwegian1}
\end{figure}

\begin{figure}[ht]
\centering
\centering
	\subfloat{\includegraphics[width=0.5\textwidth]{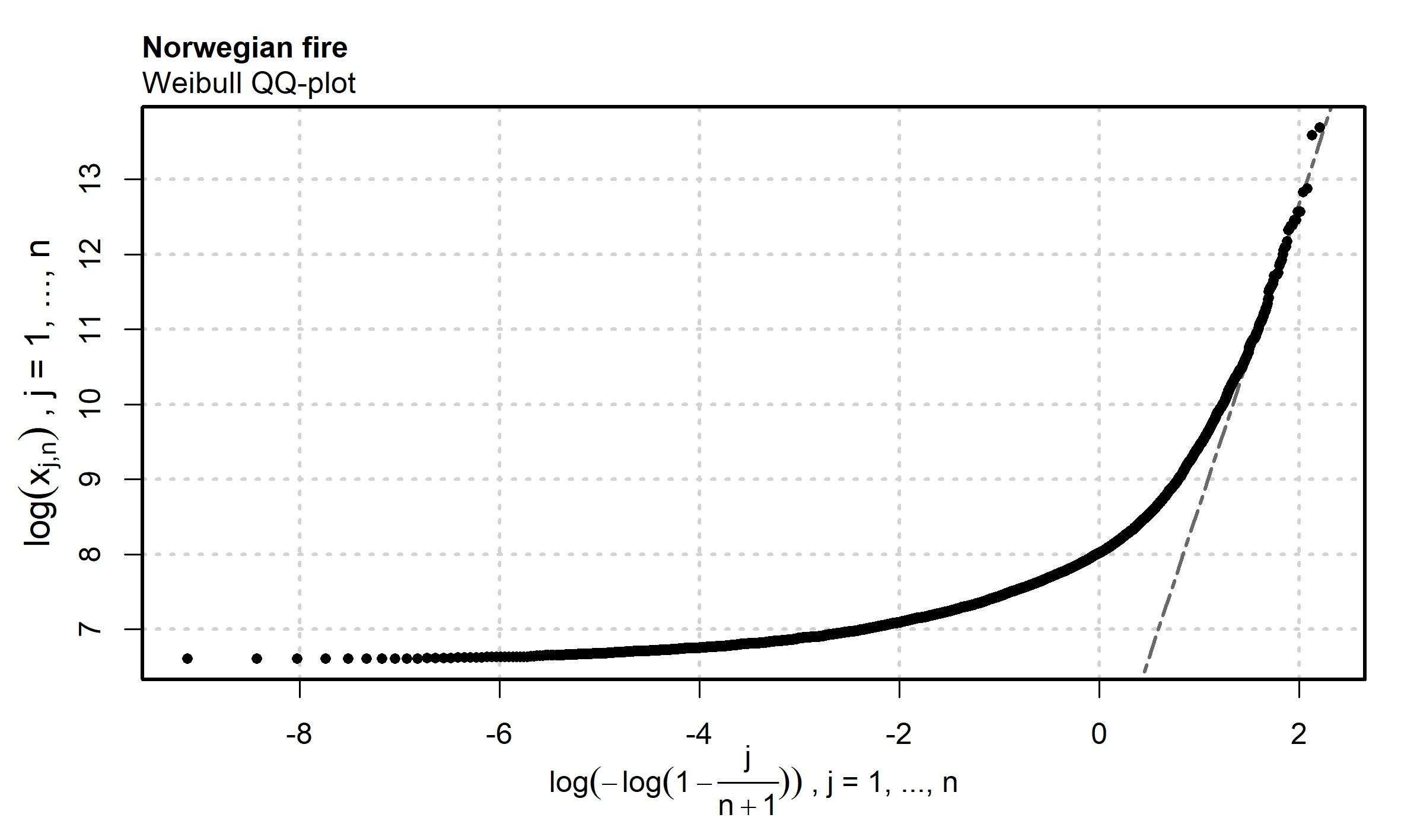}}
	\subfloat{\includegraphics[width=0.5\textwidth]{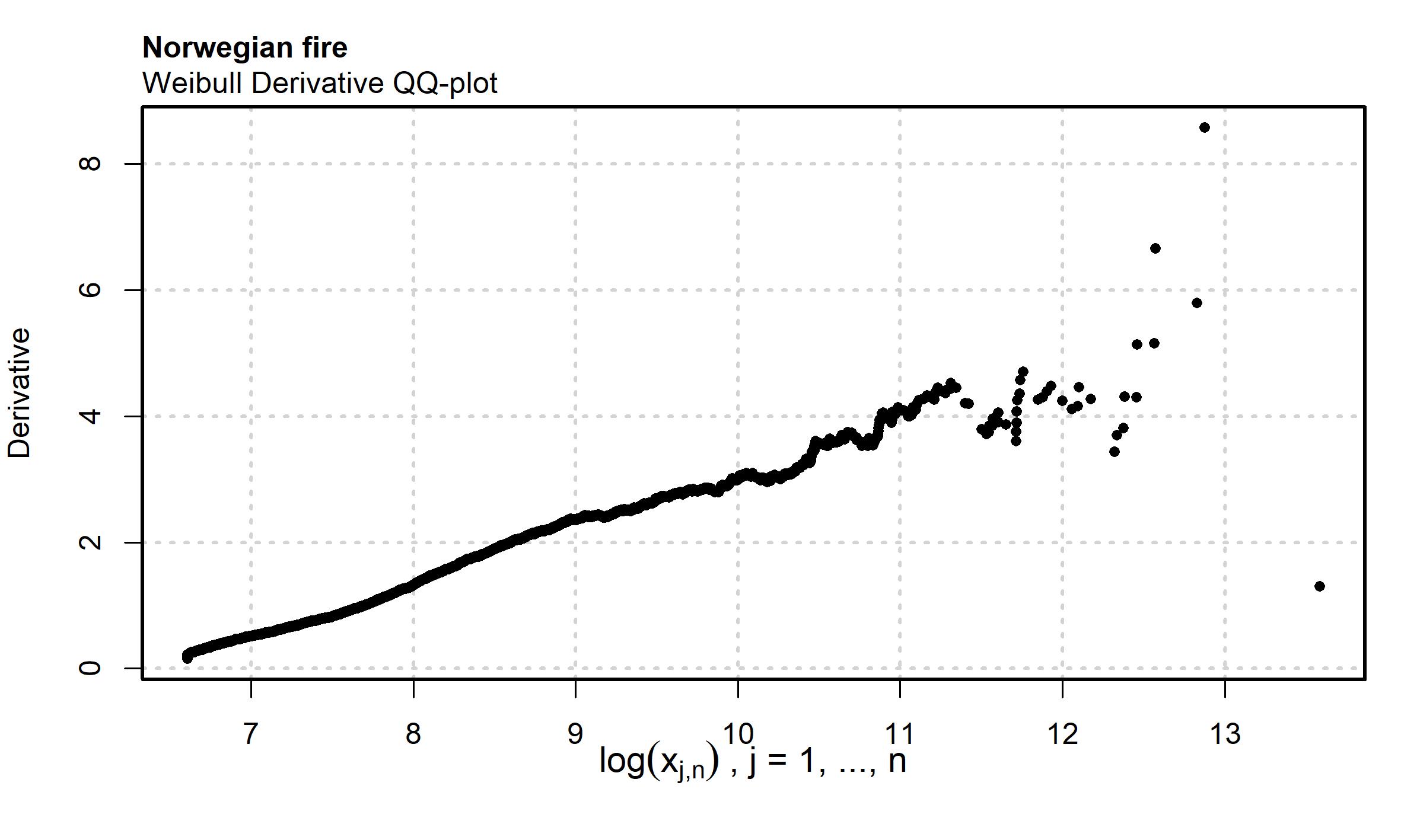}} \\
	\caption{Norwegian Fire claim data: Weibull QQ-plot (left) and Weibull derivative plot (right).}
	\label{fig_Norwegian2}
\end{figure}

\noindent
We  illustrate the need for such Weibull tempering of a Pareto-type tail  with the Norwegian fire insurance  data set discussed in Beirlant et al.\  \cite{beirlant1996practical}, which contains the year of occurrence of the claim  and the claim value (in thousand Krones) from 1972 until 1992, see also Brazauskas and Kleefeld \cite[Sec. 2]{brazauskas2016norwegian} for a detailed description of the data. In Figure \ref{fig_Norwegian1} these data are plotted by year of occurrence, next to a log-log plot (Pareto QQ-plot)
\[
\left(-\log \left(1- \frac{j}{n+1}\right),  \log X_{j,n} \right), \; j=1,\ldots,n,
\]
where $X_{1,n} \leq X_{2,n} \leq \ldots \leq X_{n,n}$ denote the ordered data from a sample of size $n$. Strict Pareto behaviour corresponds to an overall linear log-log plot, but linearity only arises approximately at the top 5000 observations. Note also the bending at the largest observations in the upper right corner in the log-log plot. This tapering near the highest observations often occurs  with insurance claim data and typically is due to a stricter claim management policy for the larger claims.  This tapering is also visible when plotting the pseudo maximum likelihood estimator $\hat{\alpha}^H_k=1/H_{k,n}$ of $\alpha$ under \eqref{Patype} (cf.\ bottom plot in Figure \ref{fig_Norwegian1}), where $H_{k,n}$ denotes the Hill estimator \cite{Hill1975hill} 
\begin{equation}
H_{k,n} =  {1 \over k}\sum_{j=1}^k\log {X_{n-j+1,n} \over X_{n-k,n}}.
\label{Hill}
\end{equation} 
The latter can be considered as an estimator of the slope in the log-log plot when restricting to the top $k+1$ observations. In that sense, the statistics $H_{k,n}$ can be considered as derivatives of the Pareto QQ-plot at the top $k$ observations.  Here, the values $\hat{\alpha}_k$ exhibit a stable area for $1000 \leq k \leq 5000$ which expresses power-law behaviour beyond $X_{n-100,n}$, and make a sharp increase at the smallest $k$ values due to tapering.
Following the QQ- and derivative plot methodology from Chapter 4 in Albrecher et al.\ \cite{albrecher2017reinsurance}, one can construct a Weibull QQ-plot $(\log (-\log (1-{j \over n+1}) ),\log X_{j,n}),\; j=1,\ldots,n$, and its derivative plot in order to verify the Weibull nature of the tempering as proposed in  \eqref{WtempPa}. A Weibull tail is observed when a linear behaviour is apparent in that QQ-plot at some top portion of the data, which can then be confirmed by a constant derivative plot in that region. For the present case, Figure \ref{fig_Norwegian2} shows that the derivative plot becomes constant on average when $\log X > 11$, corresponding to a linear Weibull pattern in the QQ-plot at the top observations with vertical coordinate larger than 11.  

\begin{figure}[ht]
\centering
\centering
	\subfloat{\includegraphics[width=0.5\textwidth]{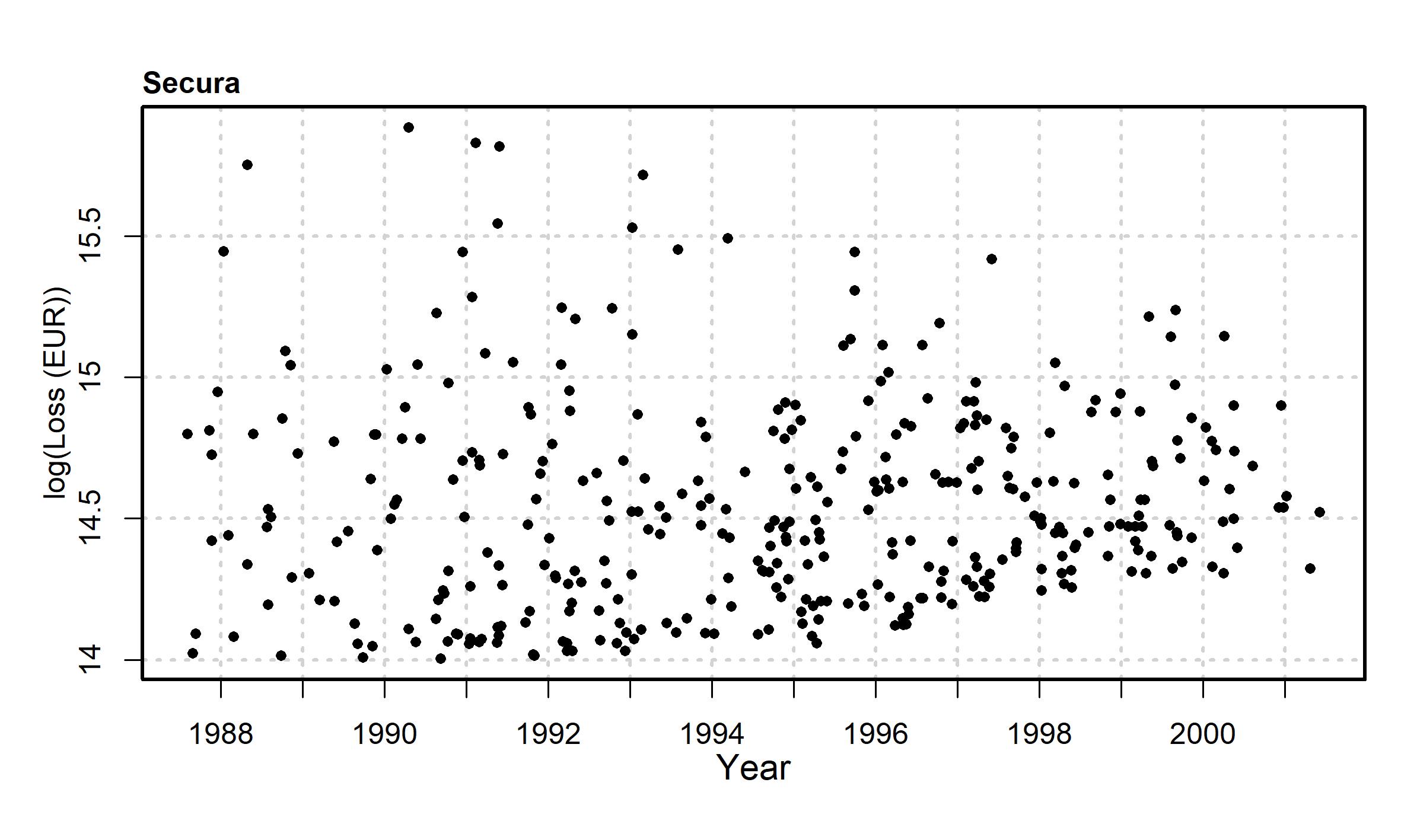}}
	\subfloat{\includegraphics[width=0.5\textwidth]{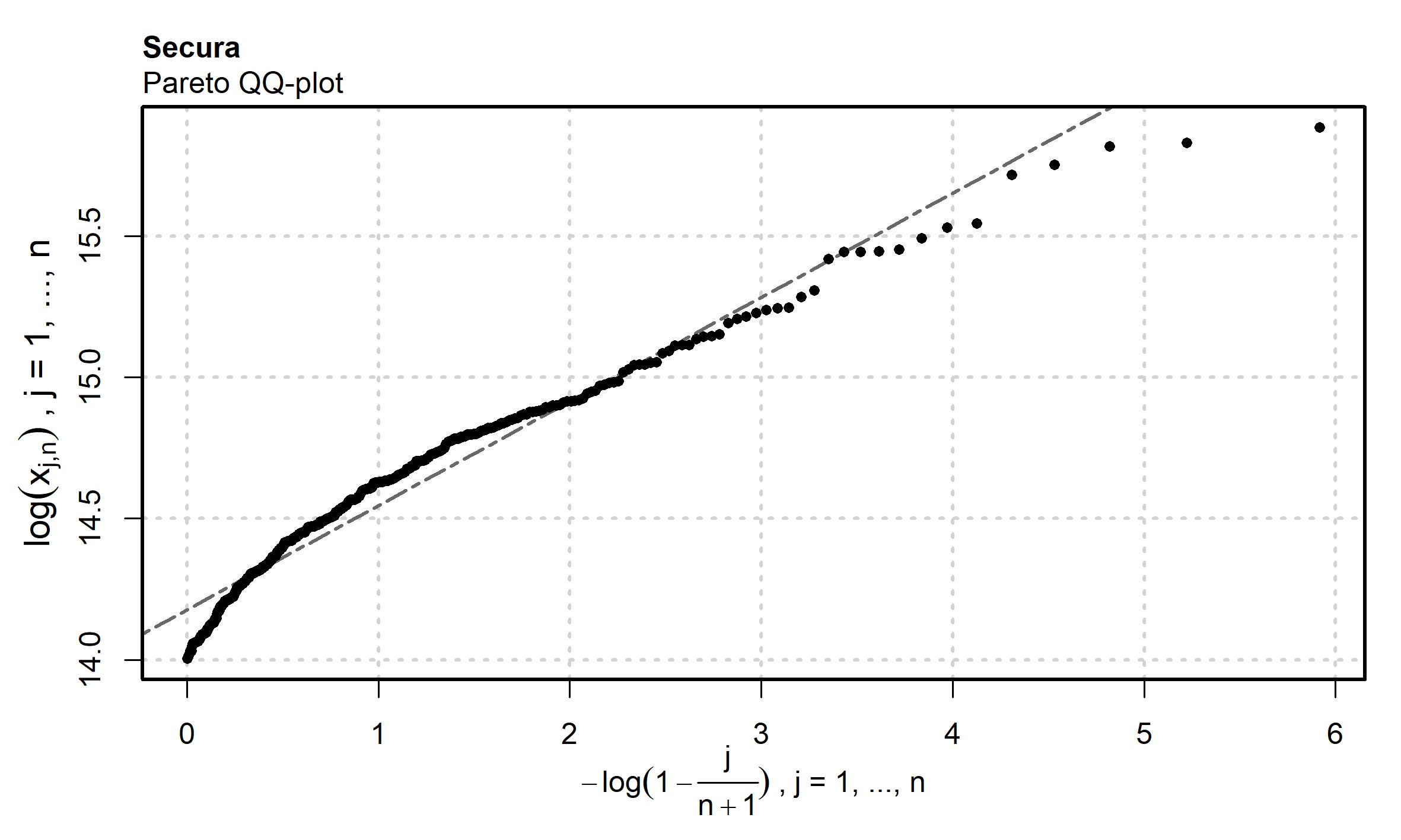}} \\
	\centering\subfloat{\includegraphics[width=0.5\textwidth]{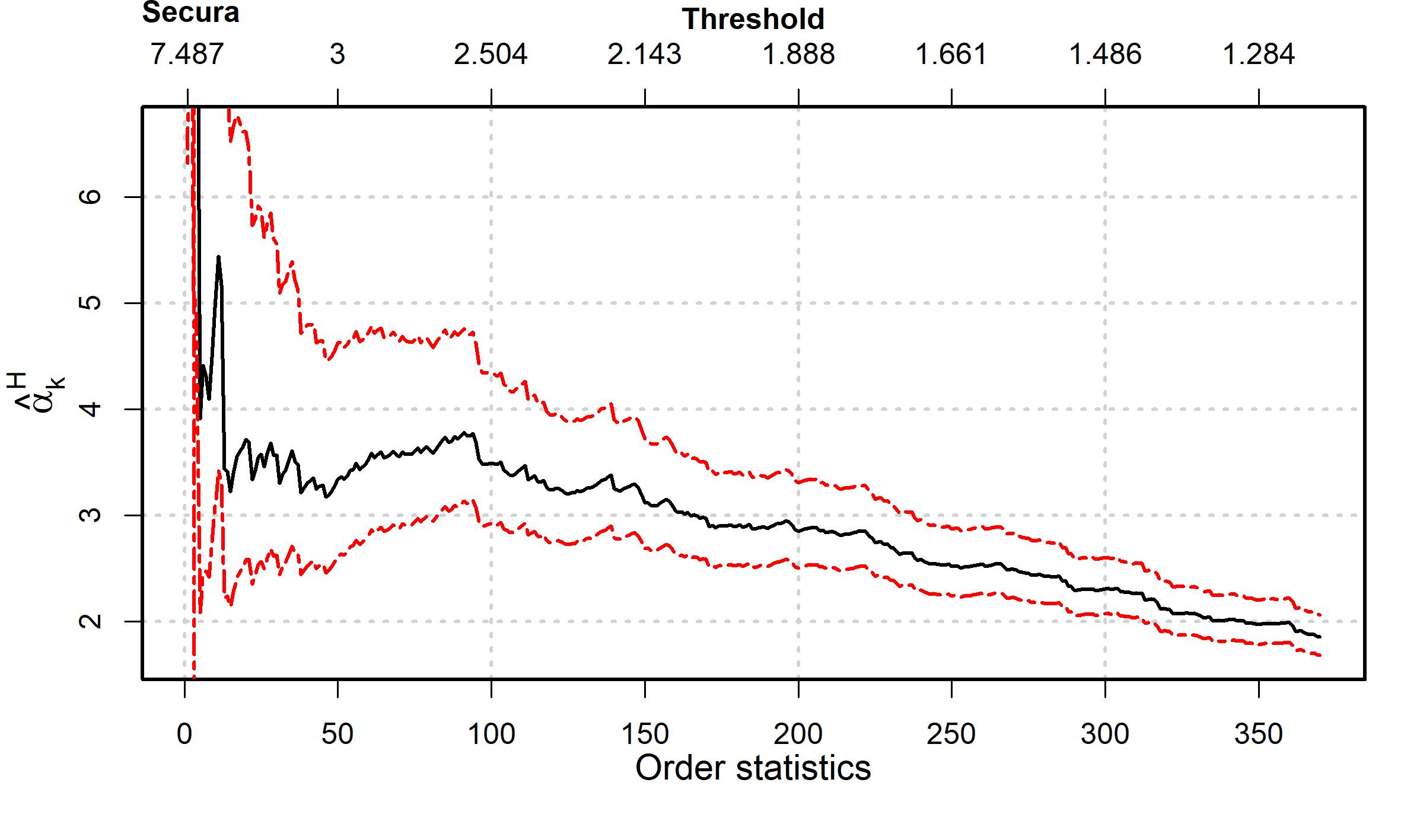}}\\
	\subfloat{\includegraphics[width=0.5\textwidth]{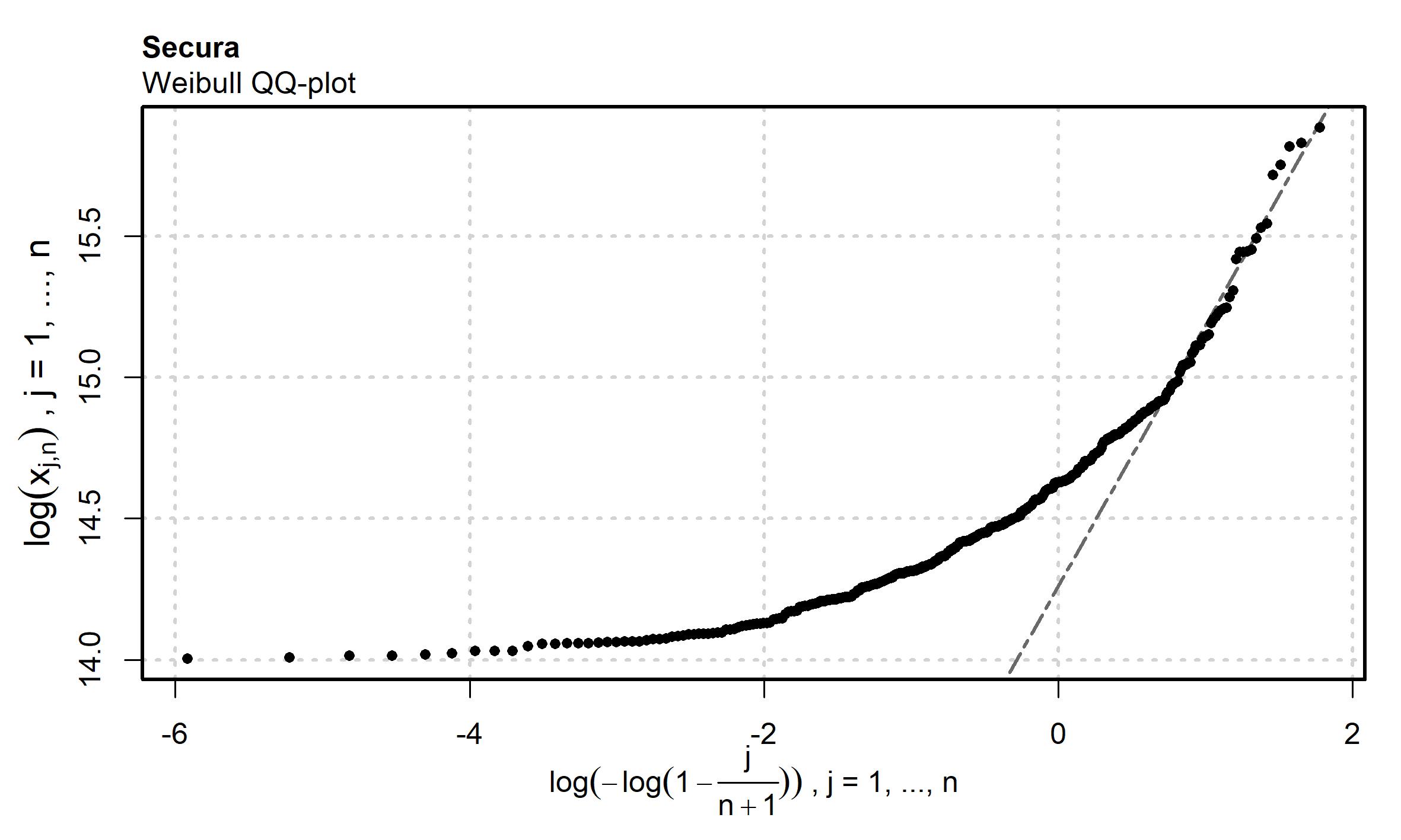}}
	\subfloat{\includegraphics[width=0.5\textwidth]{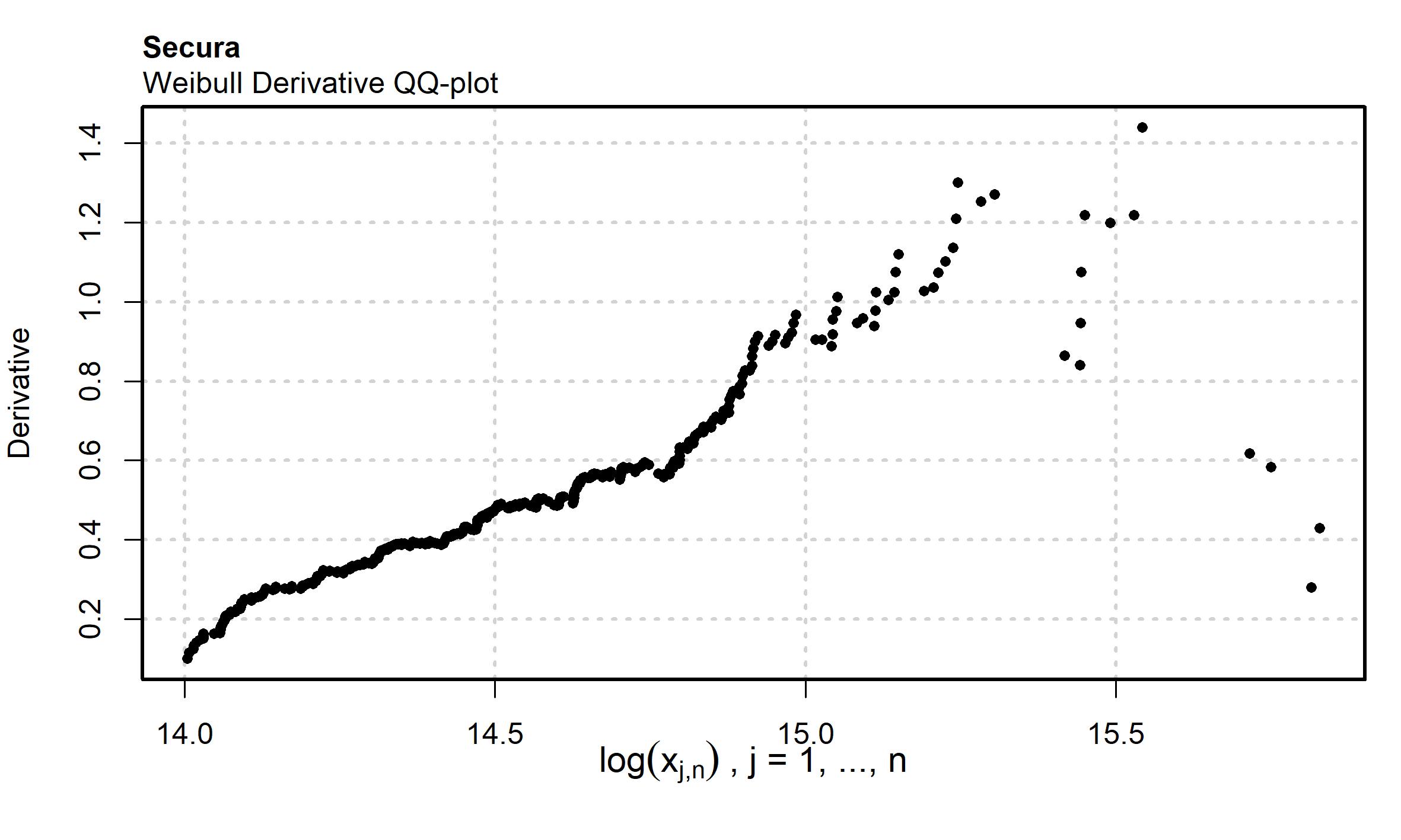}} \\
	\caption{Secura Belgian Re claim data: $\log$(Claim sizes) as a function of the year of occurrence (top left), log-log plot (top right), $\hat{\alpha}^H_k$ estimates with 95\% confidence interval  (middle), Weibull QQ-plot (bottom left) and Weibull derivative plot (bottom right).}
	\label{fig_secura1}
\end{figure}

As a second example, a tapering effect is  also observed in the Secura Belgian Re data set from  Beirlant et al.\  \cite{beirlant1996practical}. We refer the reader to Beirlant et al. \cite[Sec. 1.3.3 \& Sec. 6.2]{beirlant2006statistics} for further details about the data set. The Pareto QQ-plot in Figure \ref{fig_secura1} shows a linear pattern on from $\log X >15$, but bending is visible near the top 10 observations, leading to higher values of $\hat{\alpha}^H_k$ at $k\leq 10$. The Weibull derivative plot  shows an ultimately decreasing behaviour at the largest 10 observations. This then  could lead to  truncated Pareto  modelling rather than Weibull tempering of a Pareto-type tail, as discussed in detail in Beirlant et al.\ \cite{beirlant2016tail}.   

\noindent
In this paper, we complement the graphical and exploratory analysis of Weibull tempering of Pareto-type tails as illustrated above with a mathematical analysis of model \eqref{WtempPa}.  
 This can be considered as an alternative to the truncated Pareto-type distributions $X$ discussed in \cite{beirlant2016tail} which were defined by $X=_d Y|Y<T$ for some high value of $T$ and $Y$ satisfying Pareto-type behaviour \eqref{Patype}. Truncation also leads to tapering and appears for instance in modelling of earthquake energy levels on the basis of the Gutenberg-Richter law. 
From the viewpoint of truncation, model \eqref{WtempPa} corresponds to $X=\min(Y,W)$ with $Y$ and $W$ independent, $Y$ being Pareto-type distributed and $W$ Weibull distributed with $\Prob (W>x)= e^{-(\beta x)^\tau}$. 
Such a model is intended to describe situations  where  a gradual transit from a power-law decay to an exponentially fast decay is observed as one goes further into the tail. In view of the general nature of the Pareto-type models \eqref{Patype}, this approach will not be able to capture the characteristics over the whole range of the distribution but focuses rather on the largest observations above some threshold $X_{n-k,n}$. However, if appropriate such tempered tail fits could be spliced with different methods to describe the data below the chosen $X_{n-k,n}$, as it was done before to obtain composed models with a Pareto or generalized Pareto tail fit; see for instance  Reynkens et al. \cite{reynkens2017modelling} for mixed Erlang compositions with Pareto tails, Brazauskas and Kleefeld \cite{brazauskas2016norwegian} for  log-normal and Weibull models spliced with Pareto tail fits, and  Raschke \cite{raschke2020alternative} for Pareto-Pareto or cascade Pareto modelling.  Albrecher et al. \cite{albrecher2020matrix} considered a parsimonious and versatile family of distributions for the modelling of heavy-tailed risks using the class of matrix Mittag-Leffler distributions.\\

\noindent
In Section \ref{sec2}, we position the tempered Pareto-Weibull model in a POT approach allowing $\beta \to 0$ as the threshold $t \to \infty$, and study pseudo maximum likelihood estimation providing basic asymptotic theory. We also discuss estimation of extreme return levels and return periods. Proofs of mathematical results are deferred to the Appendix. In Section \ref{sec3} we provide simulation results, and in Section \ref{sec4} we complete the analysis of the Norwegian fire and the Belgian liability insurance data sets based on the obtained results. Section \ref{secconc} concludes. 


\section{Tempered Pareto-type modelling and estimation}\label{sec2}

Let $X = \min(Y, W)$ with $Y$ and $W$ independent, where $Y$ is Pareto-type distributed following \eqref{Patype} and
\begin{eqnarray*}
\Prob(W > x) = e^{-(\beta x)^\tau} \text{ for } x>0.
\end{eqnarray*}
The survival function of $X$ is then given by
\begin{eqnarray*}
\Prob(X > x) := \overline{F}(x) = x^{-\alpha}\ell(x) e^{-(\beta x)^\tau}.
\end{eqnarray*}
For the POT distribution $\frac{X}{t} \big| X>t$
for some threshold $t>0$, we obtain for $x>1$
\begin{align*}
\overline{F}_t(x) &:= \Prob ({X \over t}>x |X>t)\\
&=  \frac{\Prob(X> tx)}{\Prob(X>t)}\\ &= \frac{(tx)^{-\alpha}}{t^{-\alpha}}\frac{\ell(xt)}{\ell(t)} \frac{e^{-(\beta xt)^\tau}}{e^{-(\beta t)^\tau}} \\
&= x^{-\alpha} \frac{\ell(xt)}{\ell(t)} e^{-(\beta t)^\tau (x^\tau-1)}.
\end{align*}
By definition ${\ell(xt)}/{\ell(t)}\approx 1$ for large enough thresholds $t$. We then assume that at some large values of $t$, the parameter $\beta$ is inversely proportional to $t$, so that a simple Pareto-Weibull model \eqref{weiPa} provides an appropriate fit to the POTs $X/t$ ($X>t$), at least better than the simple Pareto fit with distribution function $1-x^{-\alpha}$ as used in classical extreme value methodology for Pareto-type tails. In order to formalize the above, one takes the limit for $t \to \infty$ which necessarily requires $\beta=\beta_t \downarrow 0$ as $t \uparrow \infty$. The model considered in this paper is then formally given by 
\\

\noindent
$(\mathcal{M})$
 The POT distribution $\overline{F}_t$ satisfies
\begin{eqnarray*} 
\overline{F}_t(x) \to \overline{F}_{\alpha, \beta_{\infty},\tau}(x) := x^{-\alpha}e^{-\beta^\tau_{\infty}(x^\tau-1)} ,  \mbox{ as } t \to \infty \mbox{ for every } x>1,
\end{eqnarray*}
where
\begin{itemize}
\item[a)] \textbf{(rough tempering)} $\beta = \beta_{t}$ satisfies $\beta_t t\rightarrow \beta_{\infty} >0$, corresponding to the situation where the deviation from the Pareto behavior due to Weibull tempering will be visible in the data from $t$ on and the approximation of the POT distribution using the limit distribution $\overline{F}_{\alpha, \beta_{\infty},\tau}$ appears more appropriate than using  $\overline{F}_{\alpha, 0,\tau} = x^{-\alpha},$ the simple Pareto distribution;
\item[b)] \textbf{(light tempering)} $\beta = \beta_{t}$ satisfies $\beta_t t\rightarrow 0$, corresponding to
\begin{eqnarray*}
\overline{F}_t(x) \to x^{-\alpha}, \; x >1,
\end{eqnarray*}
in which case the tempering is hardly or not visible in the data above $t$. It will then be practically impossible to discriminate light tempering from no tempering.
\end{itemize}

\vspace{0.3cm}\noindent

Given a particular threshold $t$, the quasi-likelihood procedure consists of fitting the limit distribution in $(\mathcal{M})$   to the POT data
\begin{eqnarray*}
 \frac{X_j}{t} \mbox{ when } X_j>t, \; j=1,\ldots,n.
\end{eqnarray*}
We also use the notation $\lambda =  \beta^\tau_\infty$, so that the limit distribution in  $(\mathcal{M})$ is given by 
\[
\overline{F}_{\alpha, \lambda,\tau}(x) = x^{-\alpha}
e^{-\lambda \,  (x^\tau -1)},\; x>1.
\]
The log-likelihood is then given by
\begin{eqnarray}
\begin{split}
\log L(\alpha, \lambda,\tau) &= -(1+\alpha)\sum_{j=1}^{n} \log\left(\Vjk\right)\Ij -\lambda \sum_{j=1}^{n} \left( \left(\Vjk\right)^\tau -1\right)\Ij \\&\hspace{0.5cm}+ \sum_{j=1}^{n} \log\left(\alpha + \lambda\tau \left(\Vjk\right)^\tau \right)\Ij .  
\end{split}
\label{WtempPa_ll}
\end{eqnarray}
In extreme value methodology the choice of a threshold $t$ is an important matter. A common practice is to select the $(k+1)$-largest observation $x_{n-k,n}$ for some $k \in \{4,\ldots,n-1\}$ as the threshold $t$, and to plot the resulting estimates as a function of the inverse rank $k$.  Many authors then suggest to find $k$ in a stable portion of these plots, if available. Data driven choices of $k$ are sometimes available minimizing the asymptotic mean squared error based on asymptotic results that describe the bias and variance for intermediate $k$ sequences. While an asymptotic result is presented below in Theorem 2.1, we here present an approach focusing on the goodness-of-fit of the tempering model to the POT data above the different thresholds $x_{n-k,n}$, using a QQ-plot approach. Then, for a given value of $\tau$, one finds the least-squares line that minimizes
\begin{equation}\label{eqline}
\left(-\log \left(1- \hat{F}_k\left(V_{j,k}\right)\right),  \alpha \log V_{j,k} + \tau \beta_{\infty}^{\tau} h_{\tau}\left( V_{j,k}\right) \right), \; j=1,\ldots,k,
\end{equation}
with $h_\tau (x) = (x^\tau -1)/\tau$, the POT data $V_{j,k} = X_{n-j+1, n} / X_{n-k,n}, ~ j = 1, \ldots, k,$, and $\hat{F}_k$ denoting the empirical distribution function based on those POTs. Therefore, since $ \hat{F}_k\left(V_{j,k}\right)= \frac{j}{k+1}$, one is led to minimize
\begin{align}
WLS(V_{j,k}; \alpha_k, \delta_k, \tau_k) := \sum_{j=1}^k w_{j,k} \left( \frac{1}{\alpha} \log \frac{k+1}{k-j+1} - \log V_{j,k} - \delta h_{\tau}\left( V_{j,k}\right) \right)^2,
\label{WtempPar_QQ}
\end{align}
with respect to $\alpha$ and $\delta=\tau \beta_{\infty}^{\tau}$, where $\{w_{j,k}, j = 1, \ldots, k\}$ are appropriate weights. In particular, if $w_{j, k} = 1/\log \left( \frac{k+1}{k-j+1} \right)$ when $\delta \downarrow 0$, i.e.\ without tempering, we recover the classical Hill estimator $H_{k,n}$.\\

Optimization using \eqref{WtempPar_QQ} also leads to an adaptive selection method for choosing $k$ which gives  appropriate estimates for $(\alpha, \tau, \beta_{\infty})$, choosing the $k$ for which the WLS value is minimal:
\begin{equation}
\hat{k} = \argmin_k SS_k \label{kopt}
\end{equation}
with
\begin{equation}
SS_k= \sum_{j=1}^k  {1 \over \log \left( \frac{k+1}{k-j+1} \right) }
\left( \frac{1}{\hat\alpha_k^W} \log \left( \frac{k+1}{k-j+1}\right) - \log V_{j,k} - \hat\delta _k^W h_{\hat\tau _k^W}\left( V_{j,k}\right) 
\right)^2.\label{SSk}
\end{equation}

Since for $\tau \to 0$ the parameters $\alpha$ and $\tau$ become non-identifiable, numerical issues will arise during the statistical estimation procedure when directly optimizing the likelihood, or when minimizing \eqref{WtempPar_QQ}. However, fixing a value of $\tau$ during the calibration procedure  reduces numerical instabilities. The  optimization procedure Algorithm \ref{alg:WLS} which is used in the simulations and cases,  leads to weighted least-squares estimates  $(\widehat{\alpha}^{W}_k, \widehat{\beta}^{W}_{\infty,k}, \widehat{\tau}^{W}_k)$ and maximum likelihood estimates $(\widehat{\alpha}^{M}_k, \widehat{\beta}^{M}_{\infty,k}, \widehat{\tau}^{M}_k)$, starting from a grid of $m$ initial $\tau$ values $\ttau_1 < \ttau_2 < \cdots < \ttau_m$, $m\in \N$.
\begin{algorithm}[ht]
\caption{\textbf{Estimation of  $(\widehat{\alpha}^{W}_k, \widehat{\beta}^{W}_{\infty,k}, \widehat{\tau}^{W}_k)$ and $(\widehat{\alpha}^{M}_k, \widehat{\beta}^{M}_{\infty,k}, \widehat{\tau}^{M}_k)$}  }
\label{alg:WLS}
\algsetup{indent=3em}
\begin{algorithmic}[1]
\STATE{ \textbf{set} $\ttau_1 < \ttau_2 < \cdots < \ttau_m$, $m\in \N$}
\FOR{$k = 1, 2, $ \TO $n-1$} 
{\FOR{$i = 1, 2, $ \TO $m$} 
\STATE {\textbf{Optimization step.} Set $$\left(\widehat{\alpha}_{k,\ttau_i}, \widehat{\delta}_{ k,\ttau_i}\right) := \argmin_{(\alpha > 0, \delta > 0)} WLS(V_{j,k}; \alpha, \delta, \ttau_i)$$} 
\STATE {$\widehat{WLS}_{k, \ttau_i} \leftarrow WLS(V_{j,k}; \widehat{\alpha}_{k,\ttau_i}, \widehat{\delta}_{ k,\ttau_i}, \ttau_i)$}
\STATE {$\widehat{\lambda}_{k, \ttau_i} \leftarrow \delta_{k, \ttau_i} / \ttau_i$}
\STATE {$\widehat{\log L}_{k, \ttau_i} \leftarrow \log L(V_{j,k}; \widehat{\alpha}_{k,\ttau_i}, \widehat{\lambda}_{k, \ttau_i},\ttau_i)$}
\ENDFOR}
\STATE {Set $$(\widehat{\alpha}^{W}_k, \widehat{\delta}^{W}_{k}, \widehat{\tau}^{W}_k) := \argmin_{(\widehat{\alpha}_{k,\ttau_i}, \widehat{\delta}_{ k,\ttau_i}, \ttau_i)} \left\{ \widehat{WLS}_{k, \ttau_i}; i = 1, \ldots, m \right\}$$ }
\STATE { $\widehat{\beta}^{W}_{\infty,k} \leftarrow  \left( \widehat{\delta}^{W}_k / \ttau^W_k\right)^{1/\ttau^W_k}$ }
\STATE {Set $$(\widehat{\alpha}^{M}_k, \widehat{\lambda}^{M}_{k}, \widehat{\tau}^{M}_k) := \argmax_{(\widehat{\alpha}_{k,\ttau_i}, \widehat{\lambda}_{ k,\ttau_i}, \ttau_i)} \left\{ \widehat{\log L}_{k, \ttau_i}; i = 1, \ldots, m \right\}$$ }
\STATE { $\widehat{\beta}^{M}_{\infty,k} \leftarrow  (\widehat{\lambda}^{M}_k)^{1/\ttau^M_k}$ }
\ENDFOR
\RETURN{$(\widehat{\alpha}^{W}_k, \widehat{\beta}^{W}_{\infty,k}, \widehat{\tau}^{W}_k)$ and $(\widehat{\alpha}^{M}_k, \widehat{\beta}^{M}_{\infty,k}, \widehat{\tau}^{M}_k),$ for $k = 1, 2, \ldots, n-1.$}
\end{algorithmic}
\end{algorithm}

In order to estimate return periods of the type $1/\Prob (X>z)$ for some large outcome level $z$, we use the approximation
\[
\frac{\Prob (X>tx)}{\Prob (X>t)}\approx x^{-\alpha}e^{-\lambda \tau h_{\tau}(x)}
\]
with $t$ large, so that setting $tx=z$ and $t=x_{n-k,n}$ for some $k$, we obtain the estimators for $\Prob (X>z)$
\begin{equation}
\hat{P}_{z,k}^{W} = {k+1 \over n+1}\left( {z \over x_{n-k,n}}\right)^{-\hat{\alpha}^{W}_k }\exp \left(-\hat{\lambda}_k \hat{\tau}^{W}_k h_{\hat{\tau}^{W}_k} (z/x_{n-k,n})\right)
\label{prob}
\end{equation}
and similarly $\hat{P}_{z,k}^{M}$,
 where $\Prob (X>t) = \Prob (X>x_{n-k,n})$ is estimated using the empirical proportion $(k+1)/(n+1).$
\\

\noindent
The value $z=\hat{Q}^W_{p,k}$ solving the equation 
\begin{equation}
{k+1 \over n+1}\left( { z \over x_{n-k,n}}\right)^{-\hat{\alpha}^W_k }\exp \left(-\hat{\lambda}^W_k \hat{\tau}^W_k h_{\hat{\tau}^W_k} (z/x_{n-k,n})\right)=p,
\label{quant}
\end{equation}
and similarly $\hat{Q}^M_{p,k}$, for a given value $p \leq {1 \over n}$ then yields an estimator for an extreme quantile or return level $Q(1-p)$, and hence for Value-at-Risk  (${\text{VaR}_p}$) risk measures at extreme quantile levels $1-p$.

We end this section stating the asymptotic distribution of the maximum likelihood estimators $\hat\alpha _t, 
\hat\lambda _t, \hat\tau _t$. 
The likelihood equations in $(\alpha,\lambda,\tau)$ are given by 
\begin{align*}
\sum_{j=1}^n \left\{\alpha + \lambda\tau  \left(\Vjk\right)^\tau \right\}^{-1}\Ij &= \sum_{j=1}^{n} \log \left(\Vjk\right)\Ij,\\
\sum_{j=1}^n \frac{ \left(\Vjk\right)^\tau}{\alpha + \lambda\tau  \left(\Vjk\right)^{\tau} }\Ij &= \sum_{j=1}^{n} h_\tau  \left(\Vjk\right)\Ij, \\
\sum_{j=1}^n \frac{ \left(\Vjk\right)^\tau \log  \left(\Vjk\right)}{\alpha + \lambda\tau  \left(\Vjk\right)^{\tau} }\Ij &=
\sum_{j=1}^n  \left(\Vjk\right)^\tau \log  \left(\Vjk\right) \Ij . 
\end{align*}
We further assume  classical second order slow variation
\begin{eqnarray}
\frac{\ell(ty)}{\ell(t)} = 1+D t^{\rho} h_{\rho}(y), \text{ with } D\in \R, \rho <0,
\label{2ndorder}
\end{eqnarray}
and set $\hat{\boldsymbol{\theta}}_t = 
(\hat{\alpha}_t,\hat{\lambda}_t,\hat{\tau}_t)^t$ and 
$\boldsymbol{\theta} = 
(\alpha,\lambda,\tau)^t$.

\begin{theorem}
Under $\overline{F}(x)=x^{-\alpha}\ell (x)e^{-\beta x^\tau}$ satisfying $(\mathcal{M})$ with $\beta_\infty >0$ and  $\ell$ satisfying \eqref{2ndorder}, we have as $n, t \to \infty$ such that $n\overline{F}(t) \to \infty$ and $\sqrt{n\overline{F}(t)}t^\rho \to \nu >0 $ that
\begin{eqnarray*}
 \sqrt{n\overline{F}(t)}\left( \hat{\boldsymbol{\theta}}_t - \boldsymbol{\theta}\right) 
 \to_d  
 \mathcal{N}_3 \left(D\nu  {\bf I}^{-1}{\bf b} ,{\bf I}^{-1}\right)
\end{eqnarray*}
with ${\bf I} \in \R^{3\times 3}$ symmetric and ${\bf b} \in \R^{3\times 1}$ and
\begin{eqnarray*}
I_{1,1} &=& \int_1^\infty \frac{u^{-\alpha -1}e^{-\lambda \tau h_{\tau}(u)}}{\alpha+\lambda\tau u^{\tau}} du, \\
 I_{2,2} &=& \tau^2 \int_1^\infty \frac{u^{2\tau-\alpha -1}e^{-\lambda\tau h_{\tau}(u)}}{\alpha+\lambda\tau u^{\tau}} du,\\
 I_{3,3} &=& \lambda \int_1^\infty 
 \left\{ \log^2 (u)
 -\frac{ 2\log u}{\alpha+\lambda \tau u^{\tau}} 
+\frac{\lambda u^{\tau}(1+2\tau \log u)-\alpha\tau  (\log u)^2}{(\alpha+\lambda \tau u^{\tau})^2}
 \right\}
 u^{\tau-\alpha -1}e^{-\lambda\tau h_{\tau}(u)}(\alpha+\lambda\tau  u^{\tau}) du,
 \\
 I_{1,2} &=& \tau\int_1^\infty \frac{u^{\tau-\alpha -1}e^{-\lambda\tau h_{\tau}(u)}}{\alpha+\lambda \tau u^{\tau}} du,\\
I_{1,3} &=& \lambda \int_1^\infty (1+\tau\log u)\frac{ \, u^{\tau-\alpha -1}e^{-\lambda\tau h_{\tau}(u)}}{\alpha+\lambda \tau u^{\tau}} du, \\
I_{2,3} &=&  \int_1^\infty 
\left\{ \log u - \frac{\alpha(1+\tau\log u)}{(\alpha+\lambda \tau u^{\tau})^2}  \right\}
u^{\tau-\alpha -1}e^{-\lambda\tau h_{\tau}(u)}(\alpha+\lambda\tau u^{\tau}) du,\\
b_1 &=&\int_1^\infty \left({1 \over \alpha+\lambda \tau u^{\tau}}-\log u \right)u^{-\alpha -1}e^{-\lambda\tau h_{\tau}(u)} [h_\rho (u) (\alpha+\lambda \tau u^{\tau})-u^\rho]du, \\
b_2 &=& \int_1^\infty \left({\tau u^\tau \over \alpha+\lambda\tau u^{\tau}}-\tau h_\tau (u) \right)u^{-\alpha -1}e^{-\lambda\tau h_{\tau}(u)} [h_\rho (u) (\alpha+\lambda\tau u^{\tau})-u^\rho]du, \\
b_3 &=& \lambda \int_1^\infty 
\left( \frac{1+\tau\log u}{\alpha +\lambda\tau u^\tau}
-\log u \right)
u^{\tau -\alpha -1}e^{-\lambda\tau h_{\tau}(u)} [h_\rho (u) (\alpha+\lambda\tau u^{\tau})-u^\rho]du .
\end{eqnarray*}
\end{theorem}
The derivation of this result is postponed to the Appendix.


\section{Simulation results}\label{sec3}
The finite sample behavior of the estimators $(\hat{\alpha}^W_{k},\hat{\tau}^W_{k})$ and $(\hat{\alpha}^M_{k},\hat{\tau}^M_{k})$ and the resulting tail probabilities $\hat{P}^W_{z,k}$, $\hat{P}^M_{z,k}$ and extreme quantiles $\hat{Q}^W_{p,k}$,  $\hat{Q}^M_{p,k}$ resulting from Algorithm \ref{alg:WLS}, \eqref{prob} and \eqref{quant} respectively have been studied through an extensive Monte Carlo simulation procedure. For each setting, 500 runs with sample size $n=500$ were performed. The mean and root mean squared error (RMSE) of the estimators are presented for the following models: 
\begin{itemize}
\item[(a)] Burr-Weibull($\alpha, \xi, \tau, \beta$) model with Burr survival distribution given by
\begin{align*}
F_Y(y) = 1-\left(1+y^{-\xi\alpha}\right)^{1/\xi}, ~ y>0, ~ \alpha>0,~  \xi < 0.
\end{align*}
Here \eqref{2ndorder} is satisfied with  $\rho=\xi\alpha$. We used $(\alpha, \xi, \tau, \beta)= (2, -1, 1.50, 0.50)$ and $(2, -1, 0.50, 0.20)$.
\item[(b)] Fr\'echet-Weibull($\alpha, \tau, \beta$) model with the Fr\'echet distribution function
\begin{align*}
F_Y(y)=\exp(-y^{-\alpha}), ~ y>0, ~\alpha>0.
\end{align*}
Here \eqref{2ndorder} is satisfied with  $\rho=-\alpha$. We used $(\alpha, \tau, \beta)=(2,2,0.50)$ and $(2,0.50,0.20)$.
\item[(c)] Pareto-Weibull($\alpha, \tau, \beta$) model using the Pareto distribution 
\begin{align*}
F_Y(y)=y^{-\alpha}, ~ y>1, ~\alpha>0.
\end{align*}
Here $\ell(x)=1$. We used $(\alpha, \tau, \beta)= (1,2,0.20)$.
\item[(d)] In order to study the behaviour of the estimators  under Weibull tempering of a heavy tailed distribution outside the Pareto-type family  we simulated from a tempered log-normal distribution with parameters $\mu = 0$ and $\sigma = 10$.
\end{itemize}
In the plots concerning the estimation of $\alpha$ we also plot the results for the Hill estimator $H_{k,n}$, while in case of the tail quantile estimates  $\hat{Q}^W_{p,k}$ and $\hat{Q}^M_{p,k}$ we also provide the results for the Weissman \cite{Weissman1978estimation} estimator $\hat{Q}^H_{p,k}= X_{n-k,n} \left( {k \over np}\right)^{1/\hat{\alpha}_k^H}$. Finally, we also present the boxplots of the estimates when using the adaptive choice $\hat k$ given in  \eqref{kopt} for $k$. The characteristics for the tail probability estimators  $\hat{P}^W_{z,k}$, $\hat{P}^M_{z,k}$ are quite comparable to those of the extreme quantiles, and are omitted here. 

\begin{figure}[ht]
\centering
	\subfloat{\includegraphics[width=0.5\textwidth, height = 0.5\textheight]{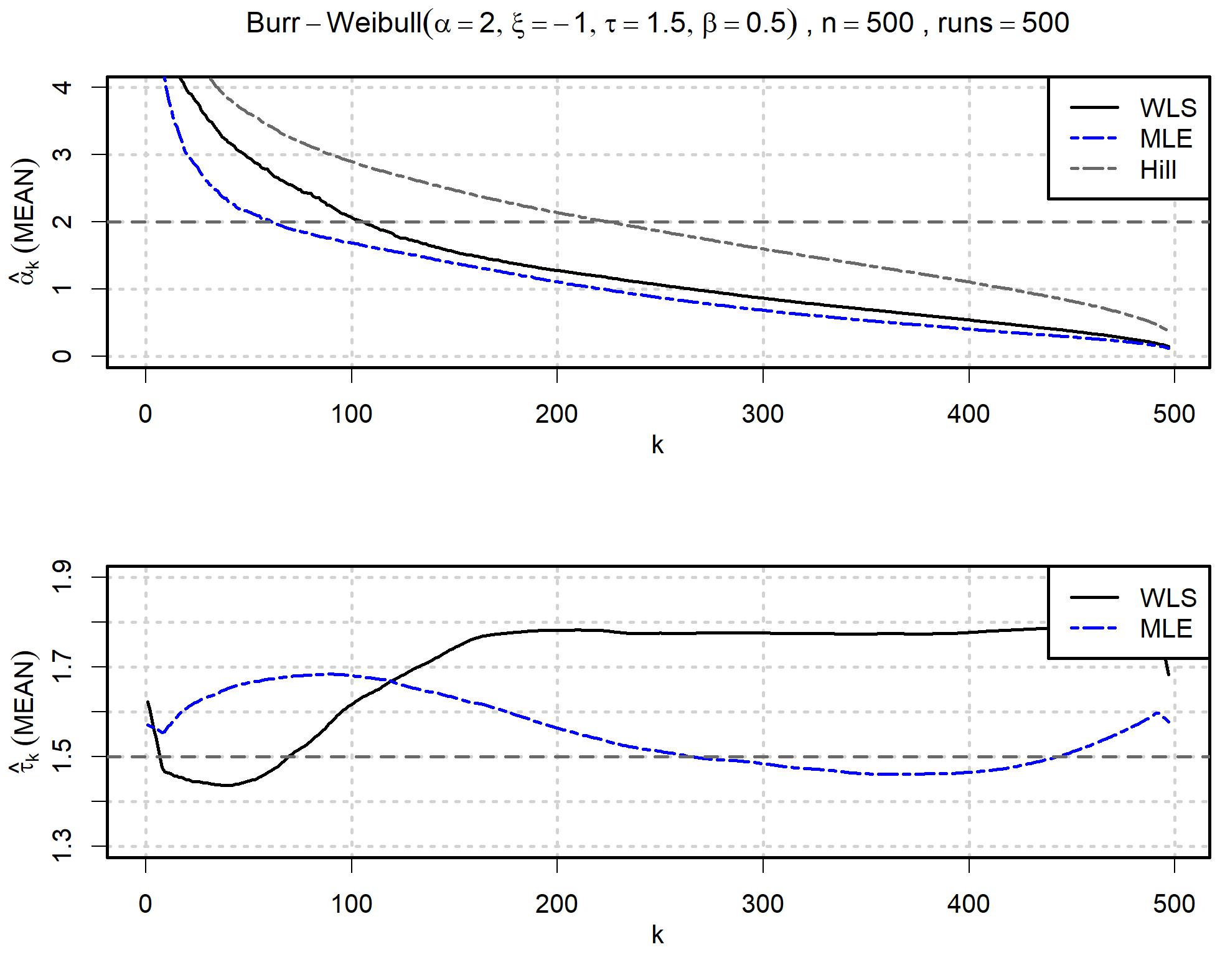}} 
	\subfloat{\includegraphics[width=0.5\textwidth, height = 0.5\textheight]{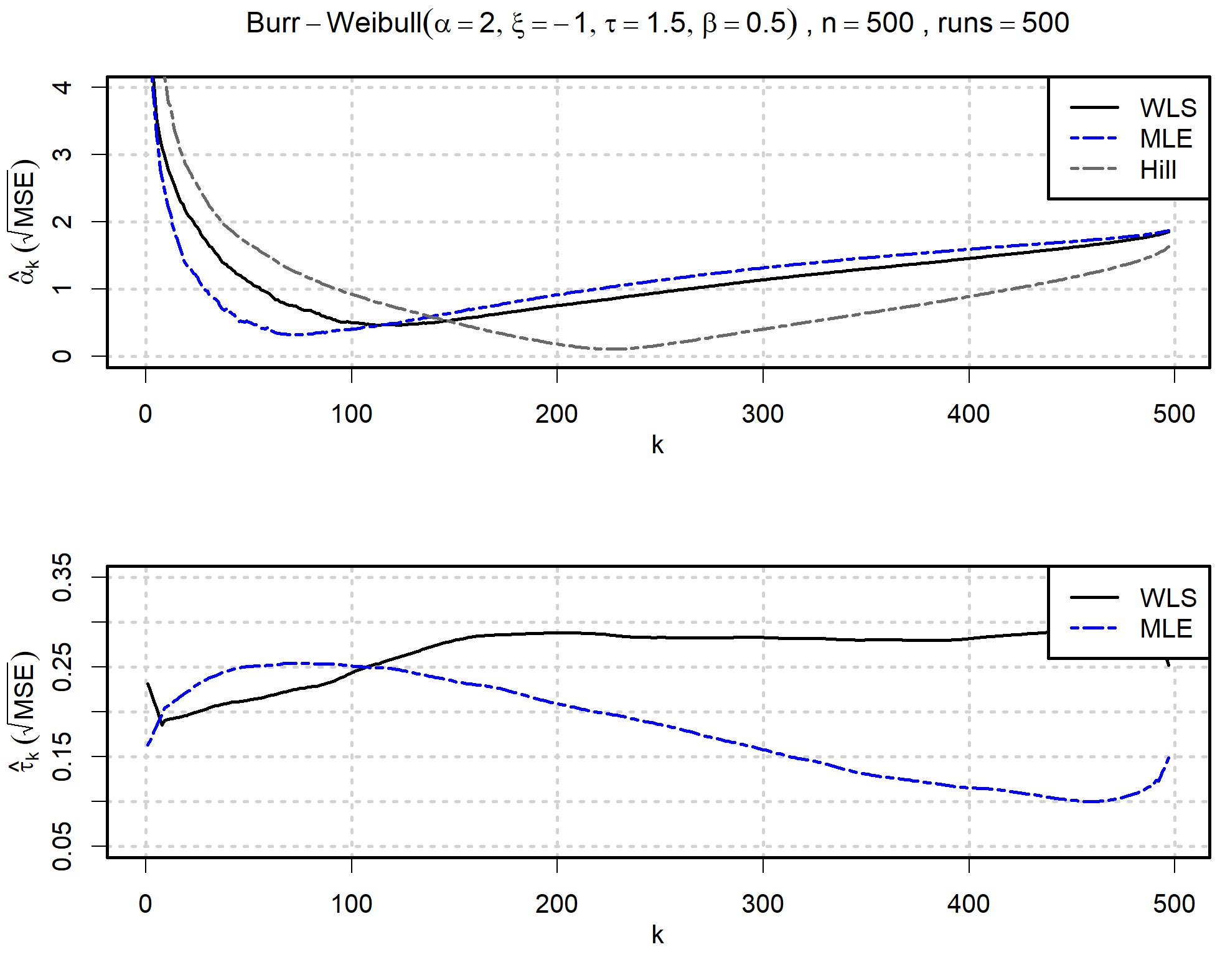}} \\
	\subfloat{\includegraphics[width=0.5\textwidth]{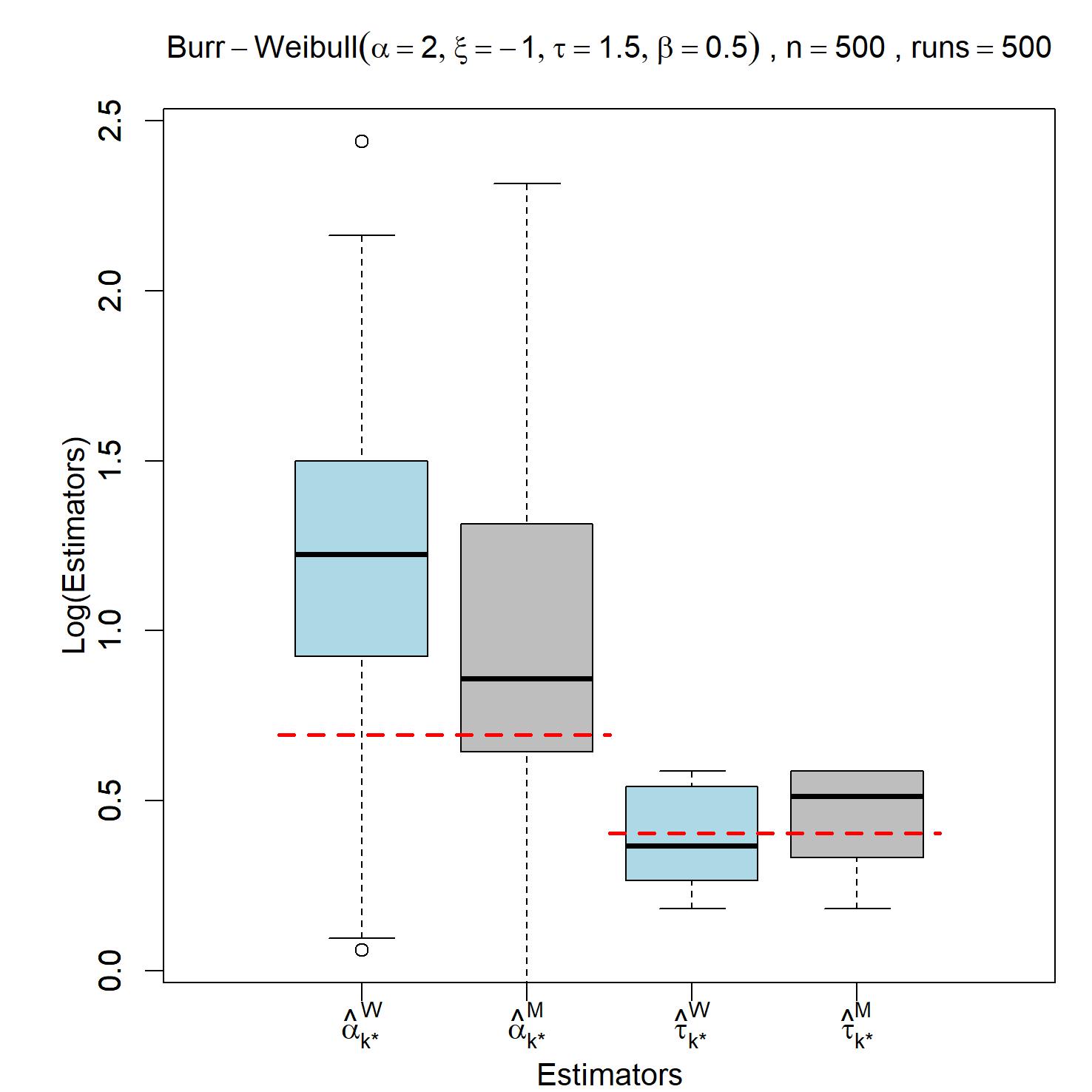}}
	\caption{Burr-Weibull($2.0, -1.0, 1.5, 0.5$). \textit{Top:} Mean (left) and RMSE (right) of  $\widehat{\alpha}^W_k$, $\widehat{\alpha}^M_k$ and $H_{k,n}$  as a function of $k$ ; \textit{Middle:} Mean (left) and RMSE (right) of $\widehat{\tau}^W_k$ and $\widehat{\tau}^M_k$ as a function of $k$; \textit{Bottom:} Boxplots of $\widehat{\alpha}^W_{\hat k}$, $\widehat{\alpha}^M_{\hat k}$, $\widehat{\tau}^W_{\hat k}$ and  $\widehat{\tau}^M_{\hat k}$ (log-scale). Horizontal dashed lines indicate the real parameters.}
	\label{SSs_BurrWeibull_ns1}
\end{figure}

\begin{figure}[ht]
\centering
	\subfloat{\includegraphics[width=0.5\textwidth, height = 0.5\textheight]{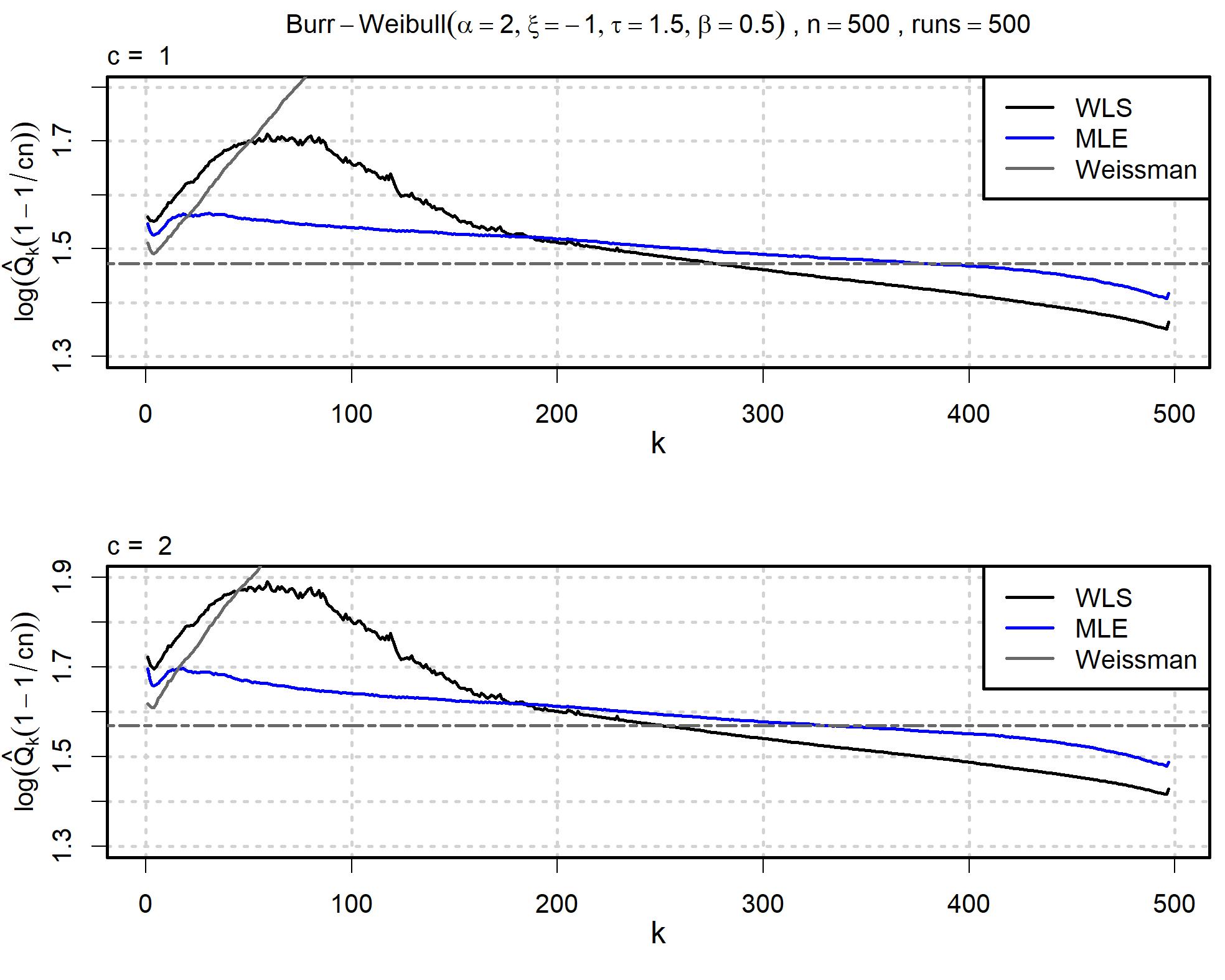}}
	\subfloat{\includegraphics[width=0.5\textwidth, height = 0.5\textheight]{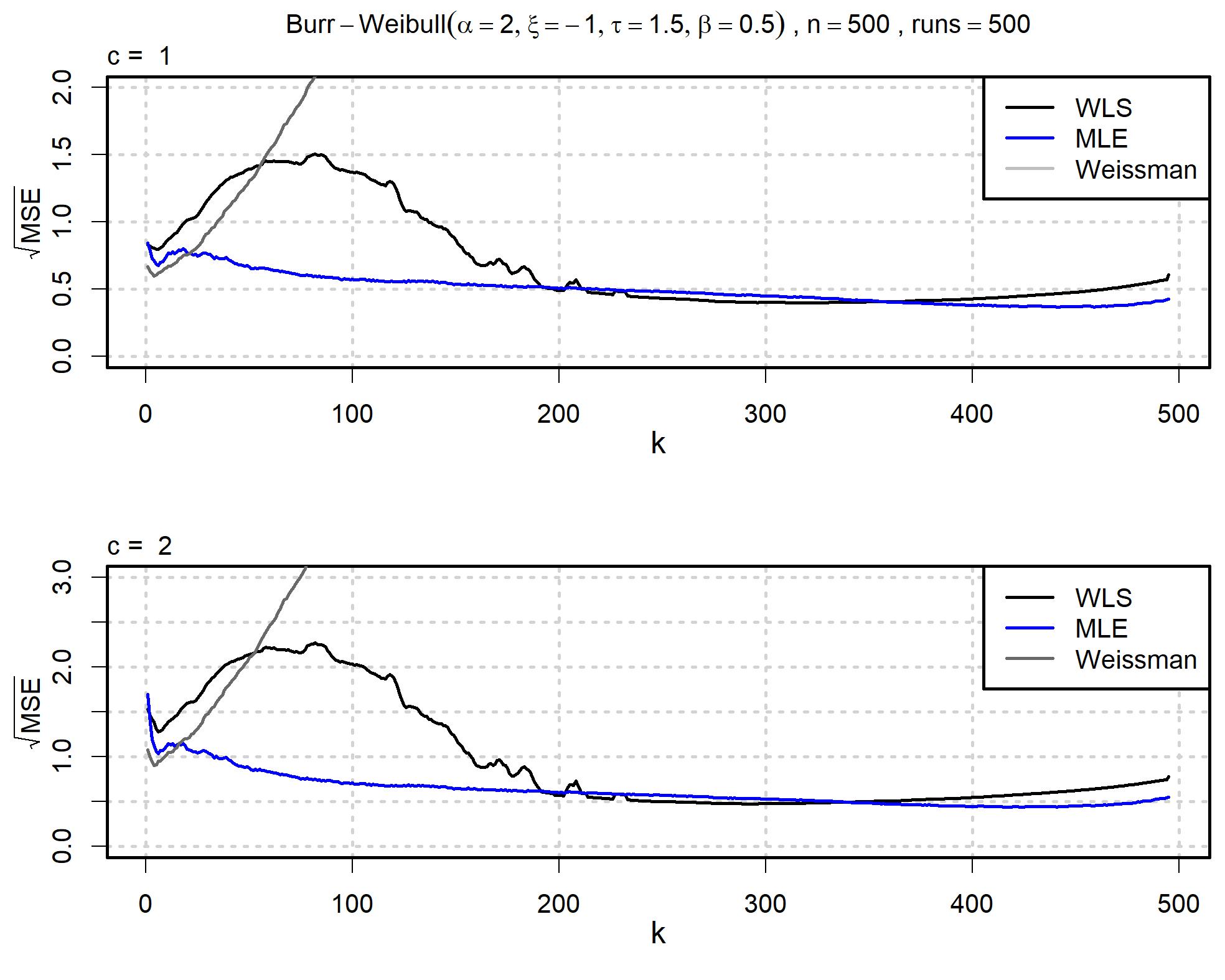}} \\
	\subfloat{\includegraphics[width=1\textwidth]{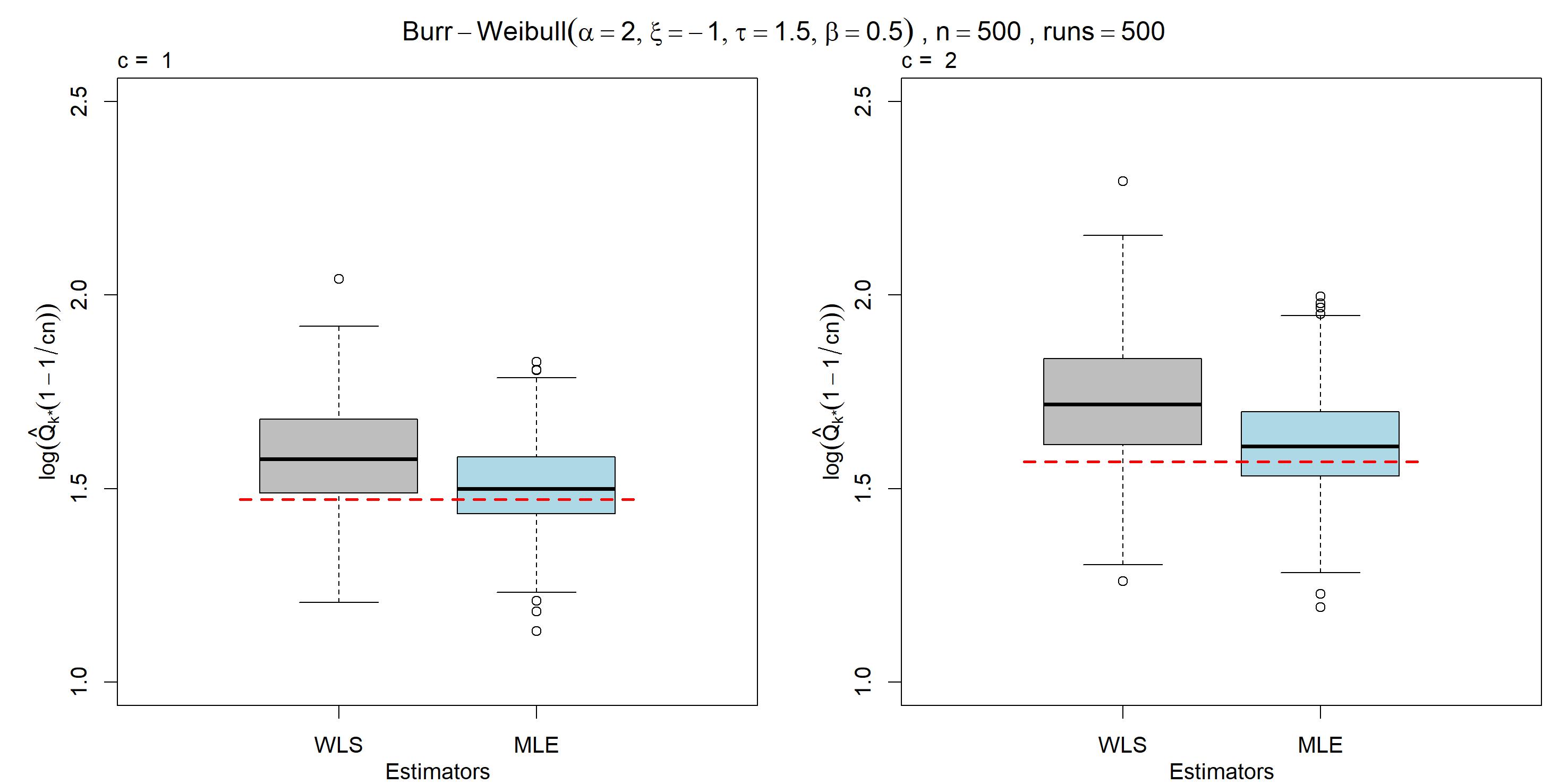}}
	\caption{Burr-Weibull($2.0, -1.0, 1.5, 0.5$): quantile estimates $\widehat{Q}^W_{p,k}$, $\widehat{Q}^M_{p,k}$  with $p={1 \over cn}$ with $c=1$ (top) and $c=2$ (middle).  Means  (left) and RMSE (right) as a function of $k$. Bottom line: boxplots of $\widehat{Q}^W_{p,\hat{k}}$, $\widehat{Q}^M_{p,\hat k}$ with $c=1$ (left) and $c=2$ (right). Horizontal dashed lines indicate the real parameters.}
	\label{SSs_BurrWeibull_Qk_ns1}
\end{figure}

\begin{figure}[ht]
\centering
	\subfloat{\includegraphics[width=0.5\textwidth, height = 0.5\textheight]{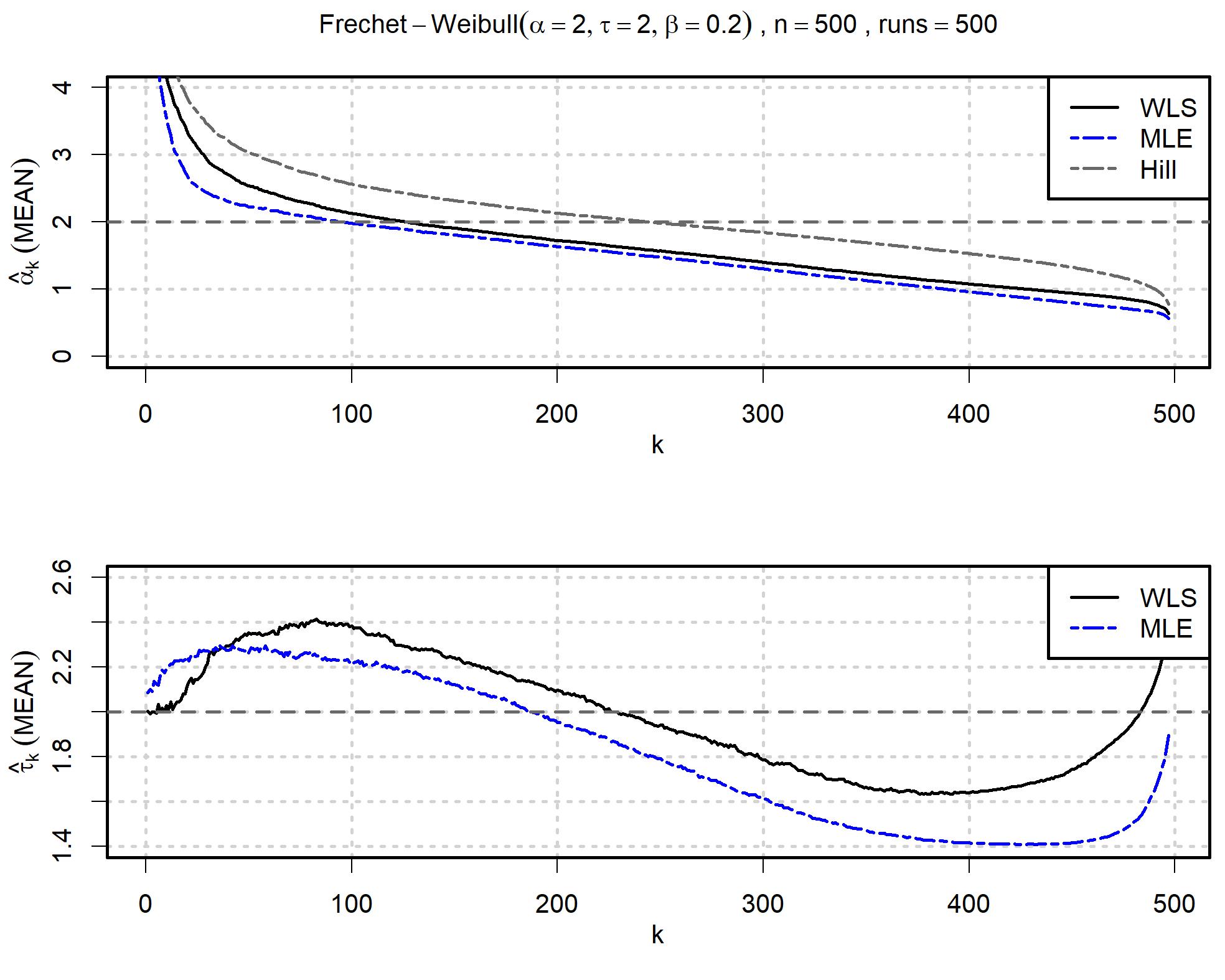}} 
	\subfloat{\includegraphics[width=0.5\textwidth, height = 0.5\textheight]{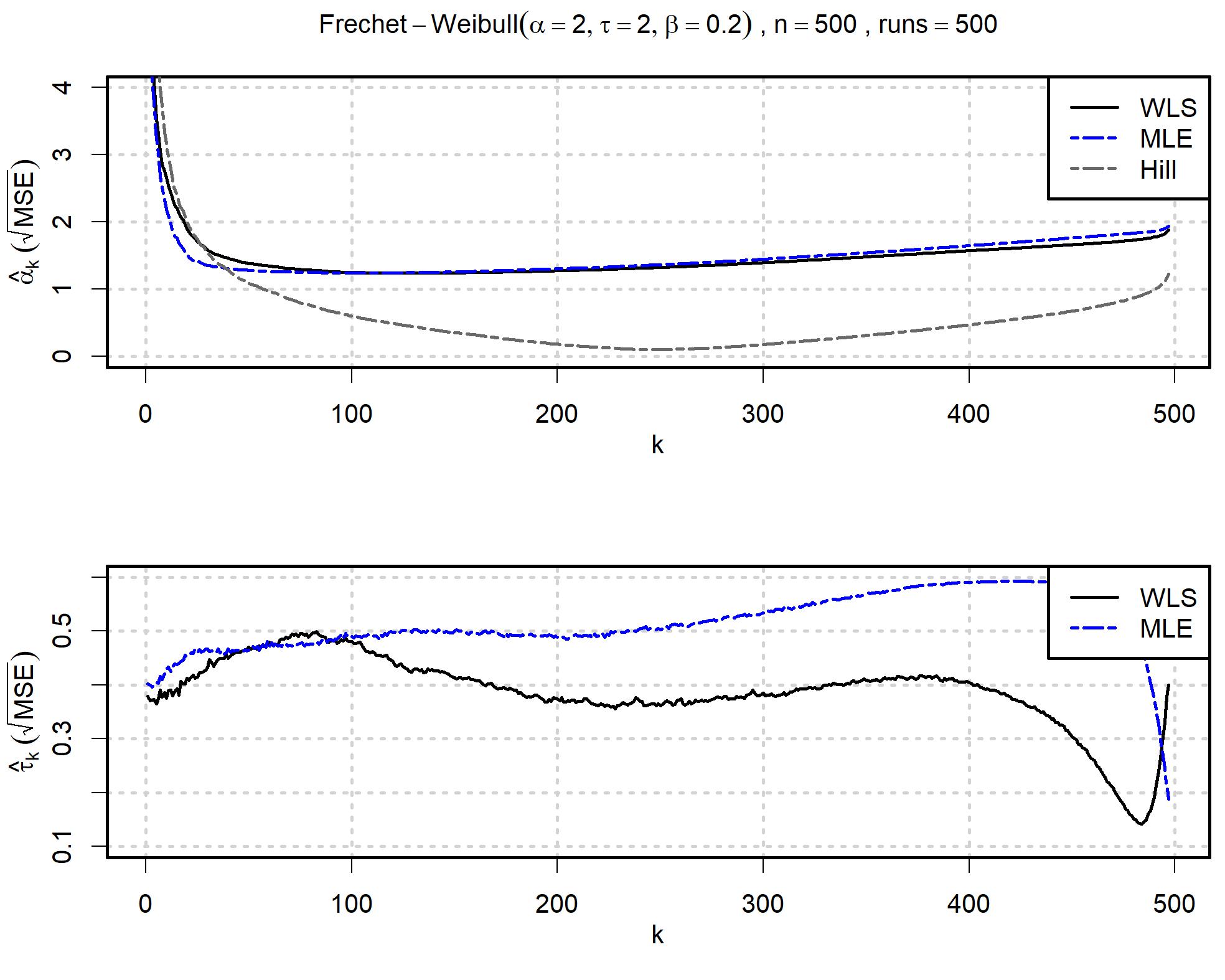}} \\
	\subfloat{\includegraphics[width=0.5\textwidth]{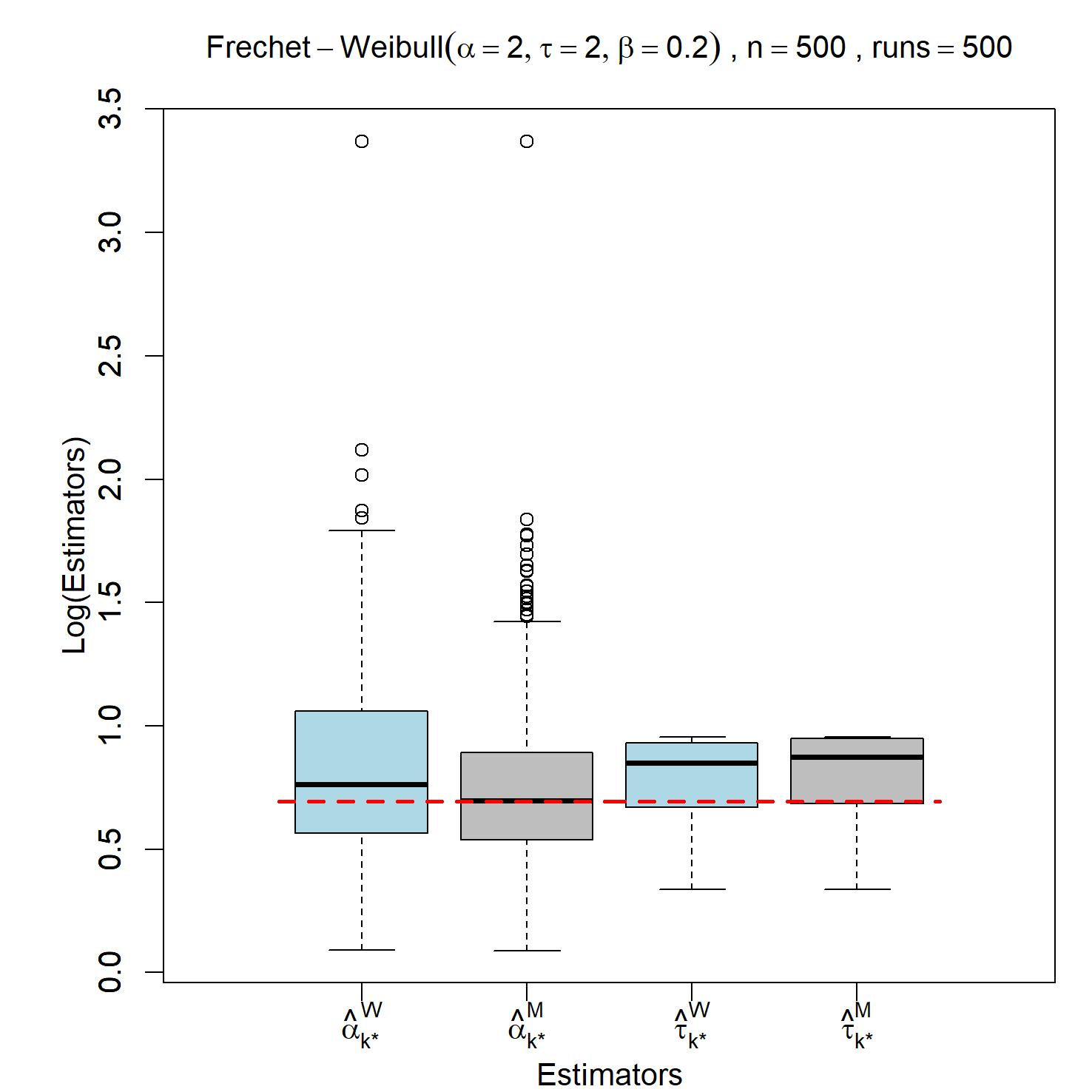}}
	\caption{Fr\'echet-Weibull($2.0, 2.0, 0.2$). \textit{Top:} Mean (left) and RMSE (right) of  $\widehat{\alpha}^W_k$, $\widehat{\alpha}^M_k$ and $H_{k,n}$  as a function of $k$ ; \textit{Middle:} Mean (left) and RMSE (right) of $\widehat{\tau}^W_k$ and $\widehat{\tau}^M_k$ as a function of $k$; \textit{Bottom:} Boxplots of $\widehat{\alpha}^W_{\hat k}$, $\widehat{\alpha}^M_{\hat k}$, $\widehat{\tau}^W_{\hat k}$ and  $\widehat{\tau}^M_{\hat k}$ (log-scale). Horizontal dashed lines indicate the real parameters.}
	\label{SSs_FrechetWeibull_ns3}
\end{figure}

\begin{figure}[ht]
\centering
	\subfloat{\includegraphics[width=0.5\textwidth, height = 0.5\textheight]{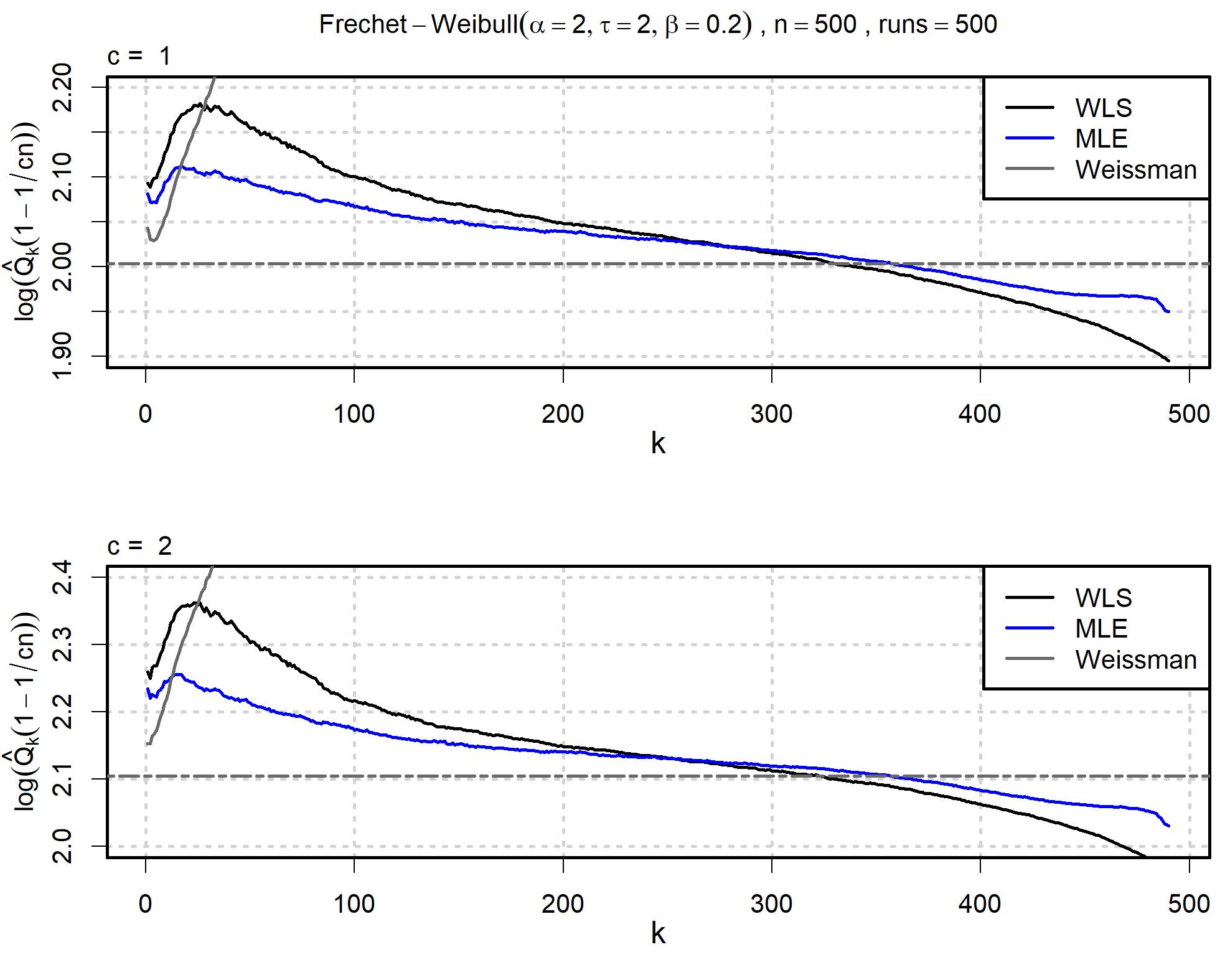}}
	\subfloat{\includegraphics[width=0.5\textwidth, height = 0.5\textheight]{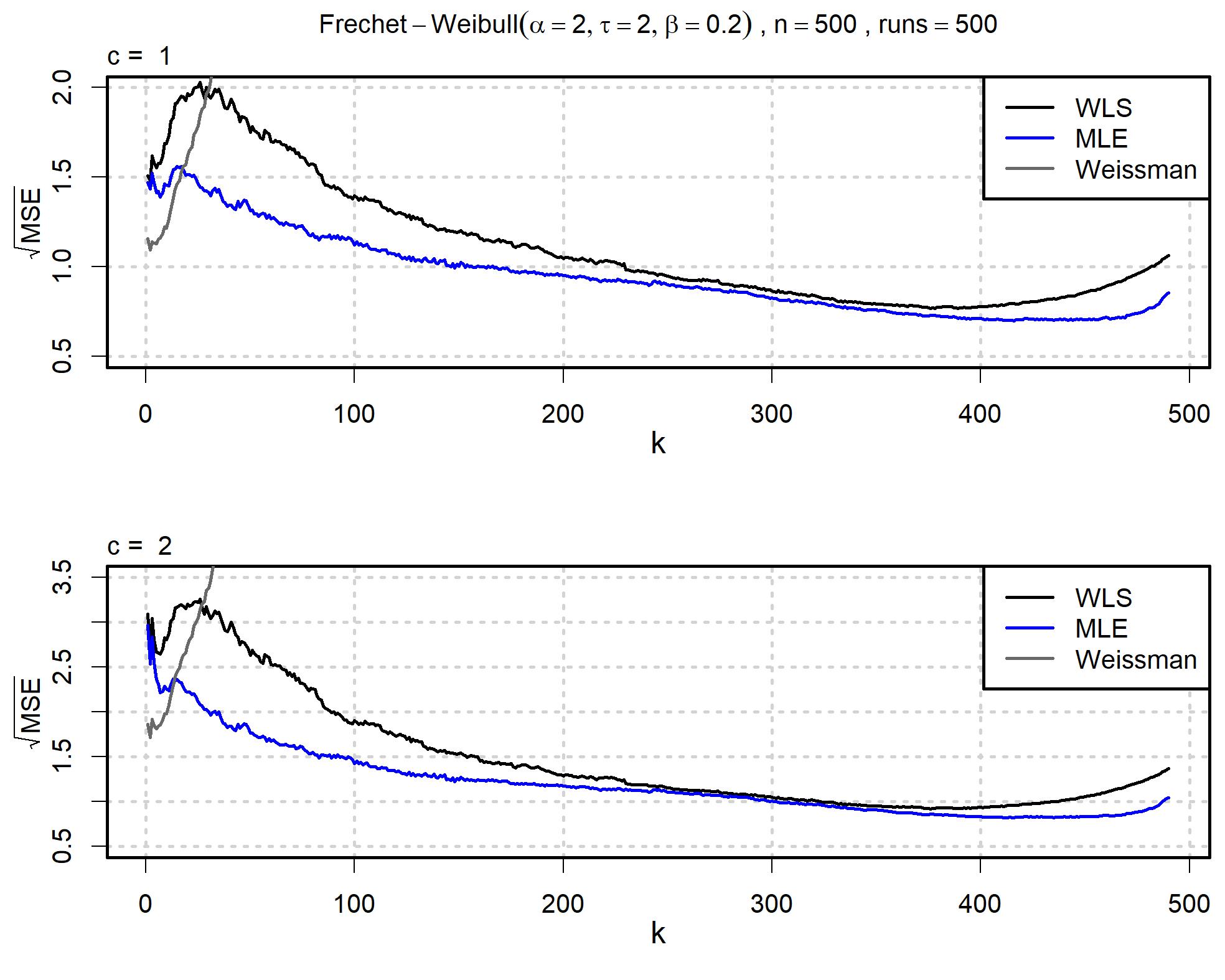}} \\
	\subfloat{\includegraphics[width=1\textwidth]{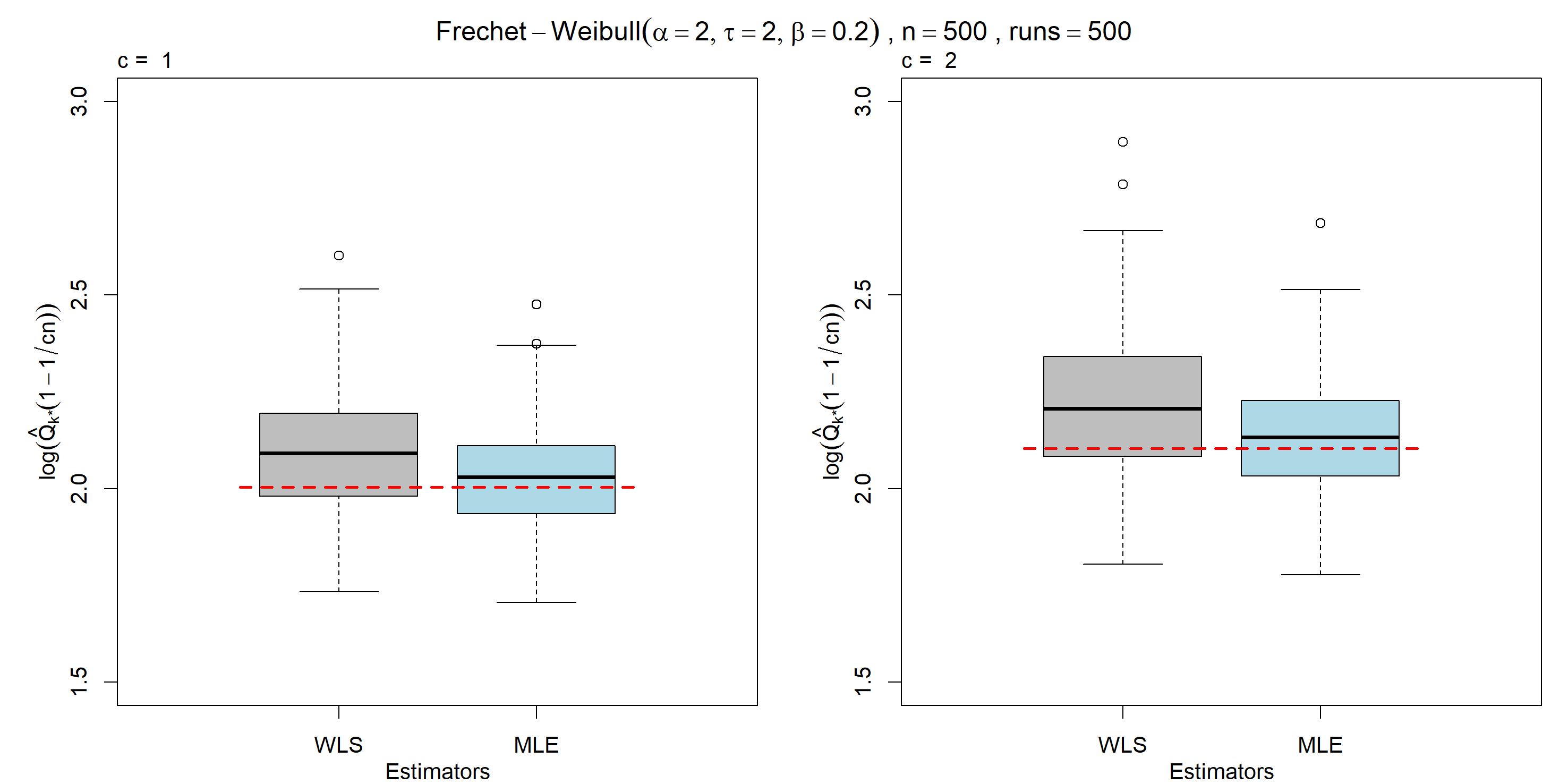}}
	\caption{Frechet-Weibull($2.0, 2.0, 0.2$): quantile estimates $\widehat{Q}^W_{p,k}$, $\widehat{Q}^M_{p,k}$ with $p={1 \over cn}$ with $c=1$ (top) and $c=2$ (middle).  Means  (left) and RMSE (right) as a function of $k$. Bottom line: boxplots of $\widehat{Q}^W_{p,\hat k}$, $\widehat{Q}^M_{p,\hat k}$ with $c=1$ (left) and $c=2$ (right). Horizontal dashed lines indicate the real parameters.}
	\label{SSs_FrechetWeibull_Qk_ns3}
\end{figure}

\begin{figure}[ht]
\centering
	\subfloat{\includegraphics[width=0.5\textwidth, height = 0.5\textheight]{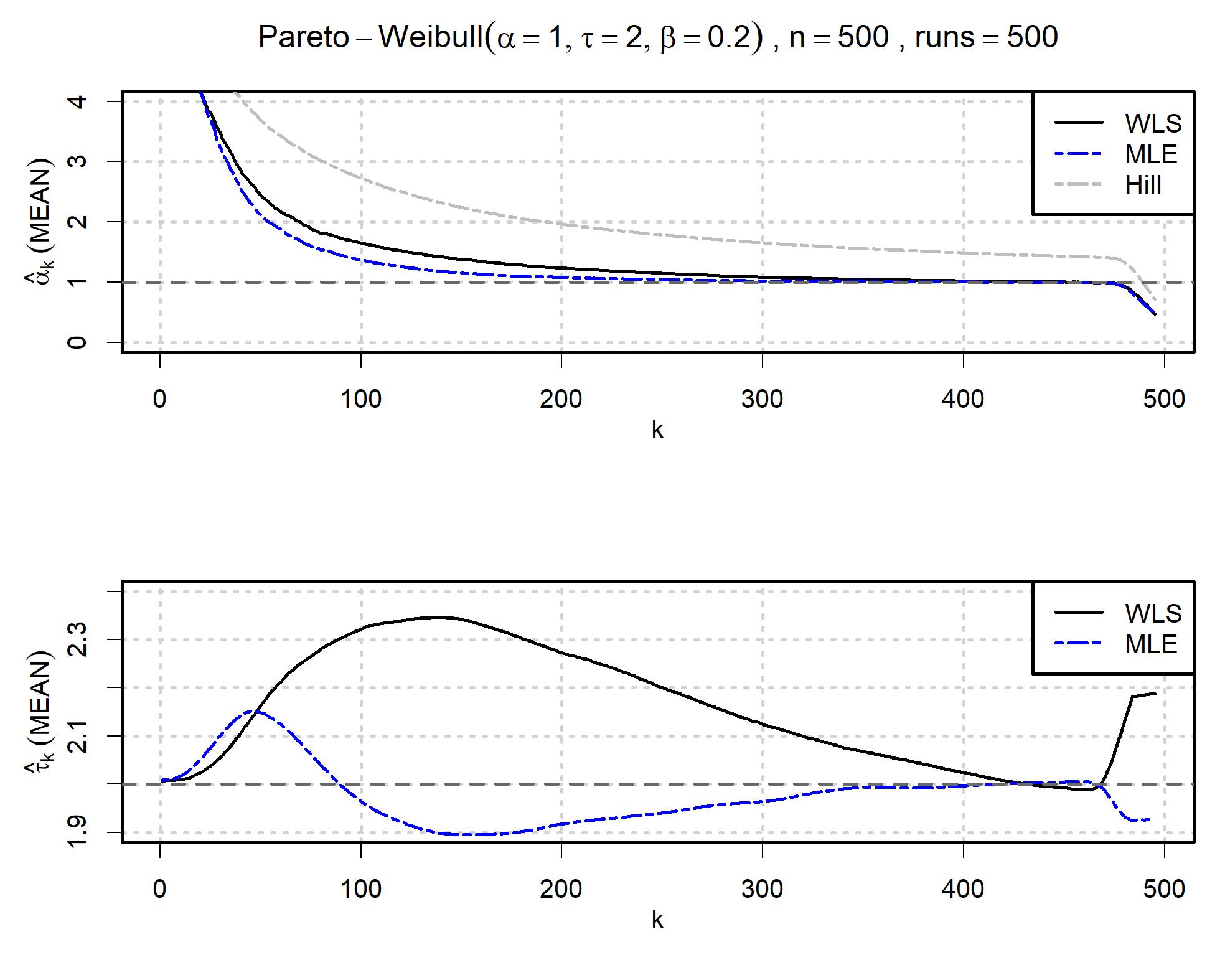}} 
	\subfloat{\includegraphics[width=0.5\textwidth, height = 0.5\textheight]{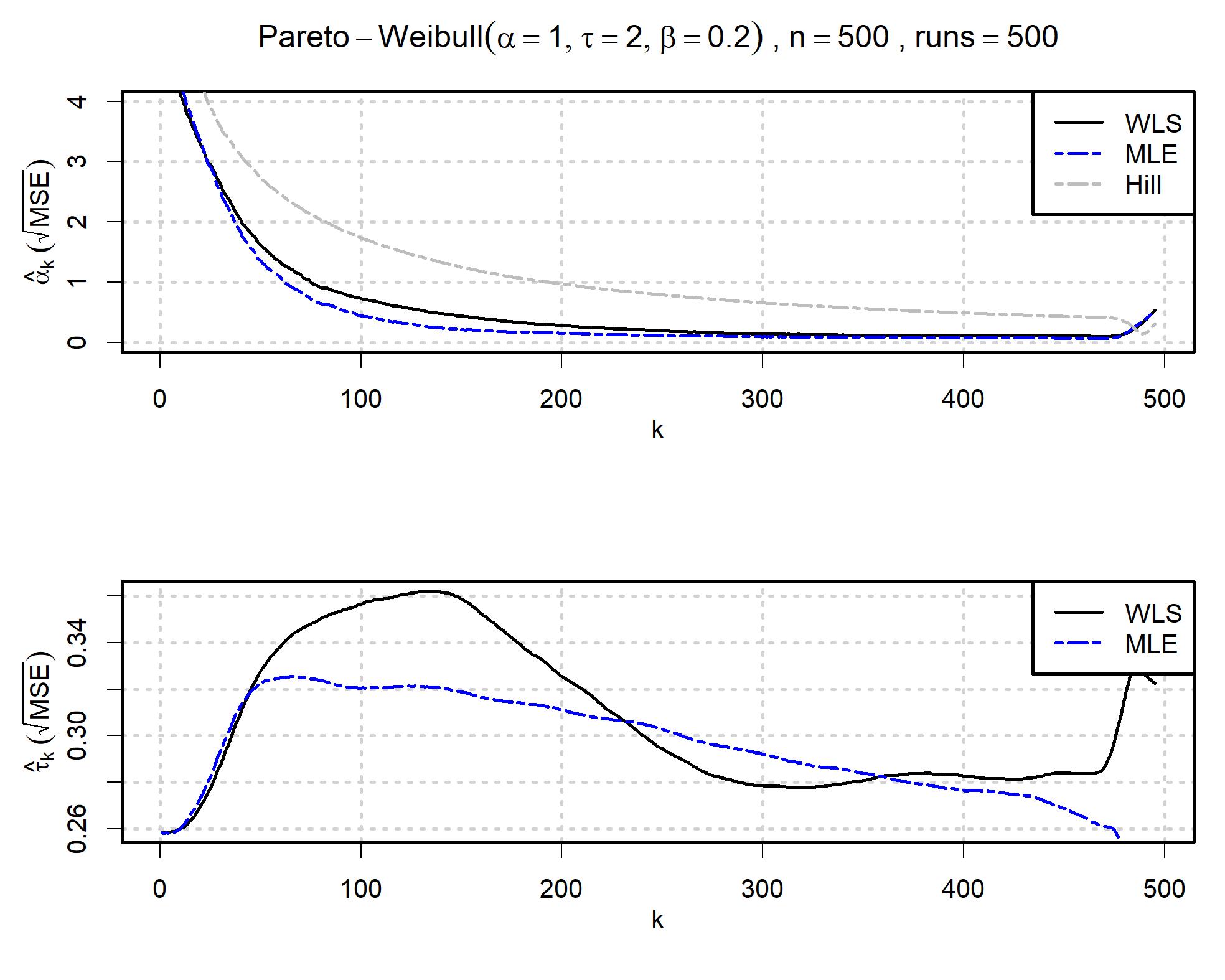}} \\
	\subfloat{\includegraphics[width=0.5\textwidth]{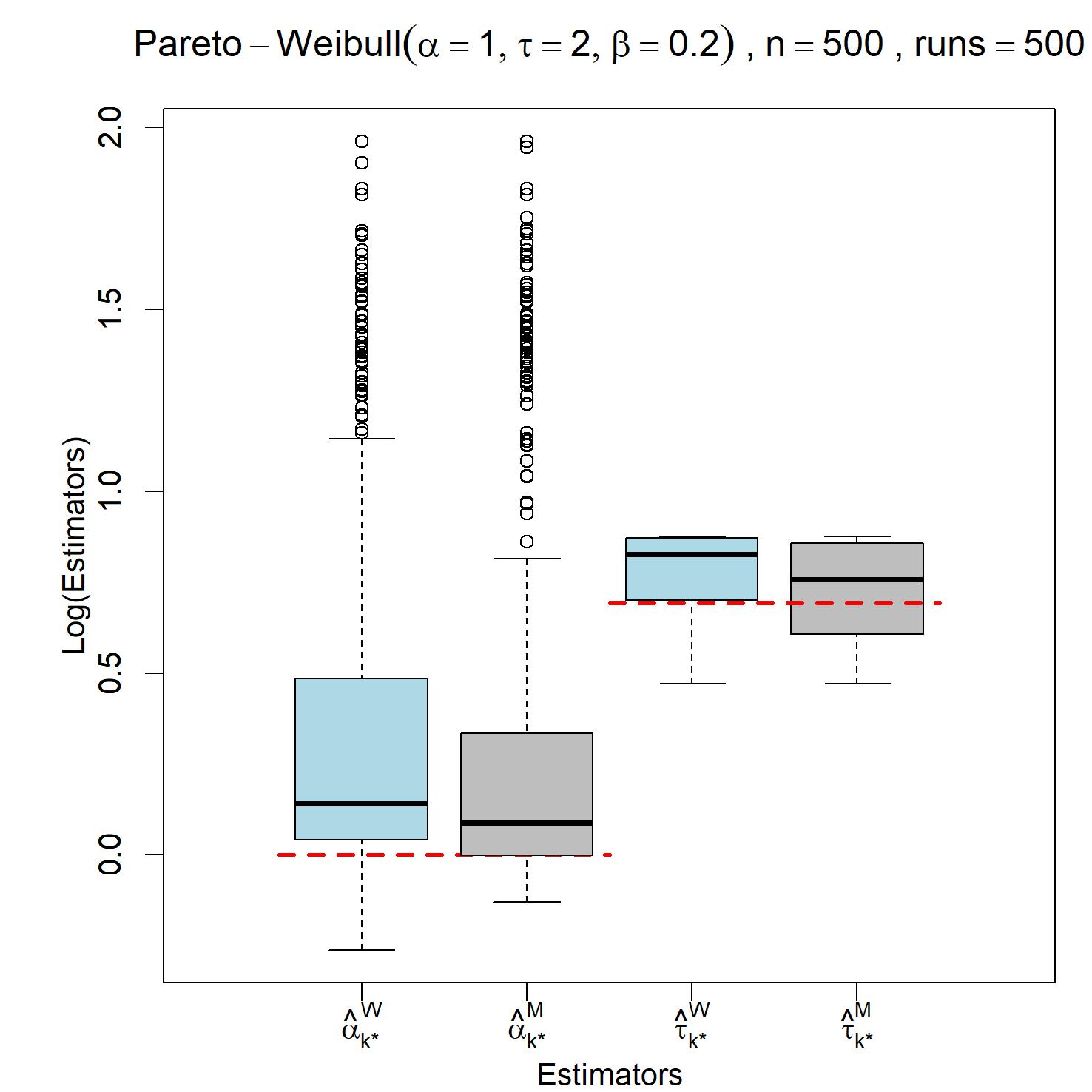}}
	\caption{Pareto-Weibull($1.0, 2.0, 0.2$).\textit{Top:} Mean (left) and RMSE (right) of  $\widehat{\alpha}^W_k$, $\widehat{\alpha}^M_k$ and $H_{k,n}$ as a function of $k$ ; \textit{Middle:} Mean (left) and RMSE (right) of $\widehat{\tau}^W_k$ and $\widehat{\tau}^M_k$ as a function of $k$; \textit{Bottom:} Boxplots of $\widehat{\alpha}^W_{\hat k}$, $\widehat{\alpha}^M_{\hat k}$, $\widehat{\tau}^W_{\hat k}$ and  $\widehat{\tau}^M_{\hat k}$ (log-scale). Horizontal dashed lines indicate the real parameters.}
	\label{SSs_ParetoWeibull_ns3}
\end{figure}

\begin{figure}[ht]
\centering
	\subfloat{\includegraphics[width=0.5\textwidth, height = 0.5\textheight]{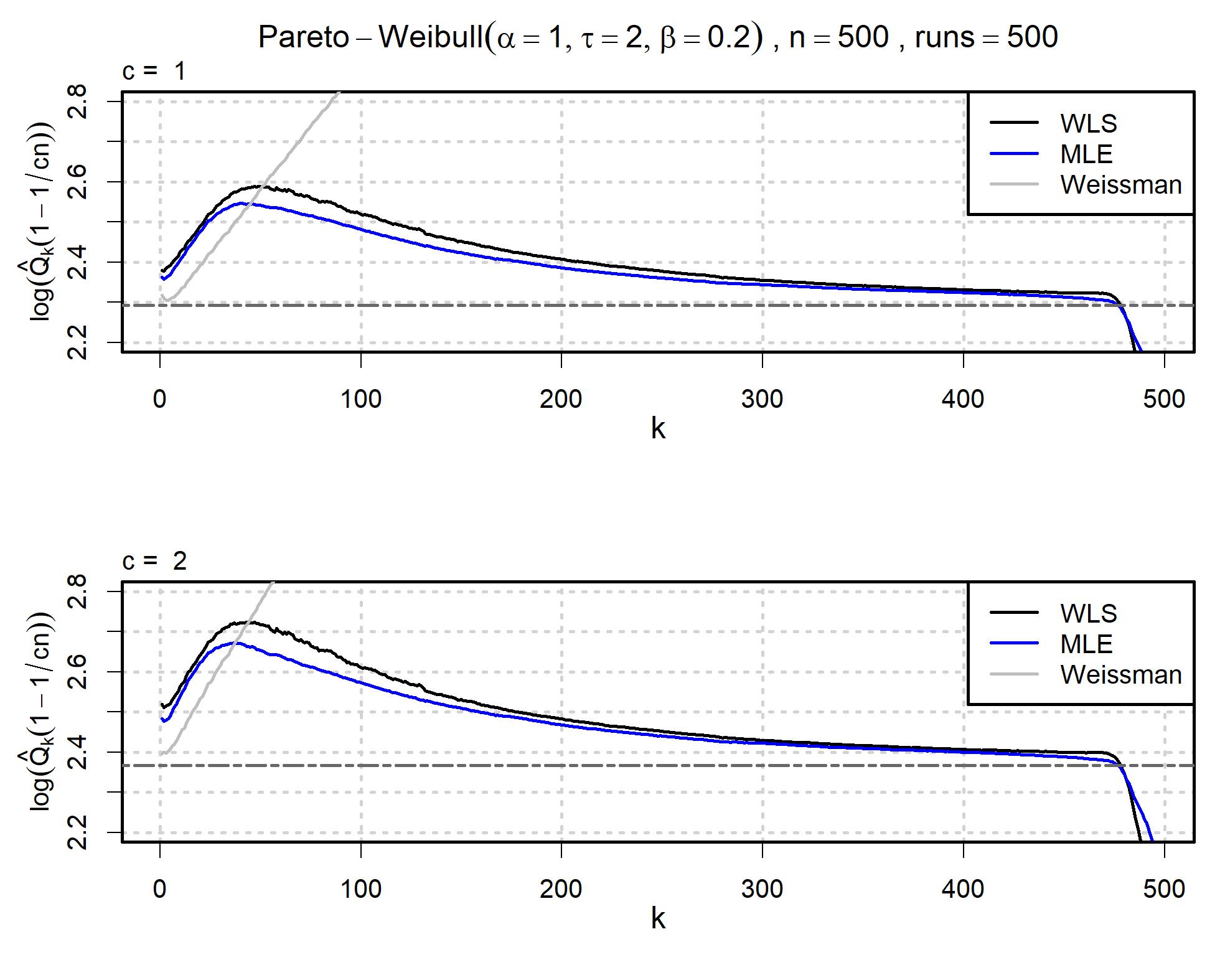}}
	\subfloat{\includegraphics[width=0.5\textwidth, height = 0.5\textheight]{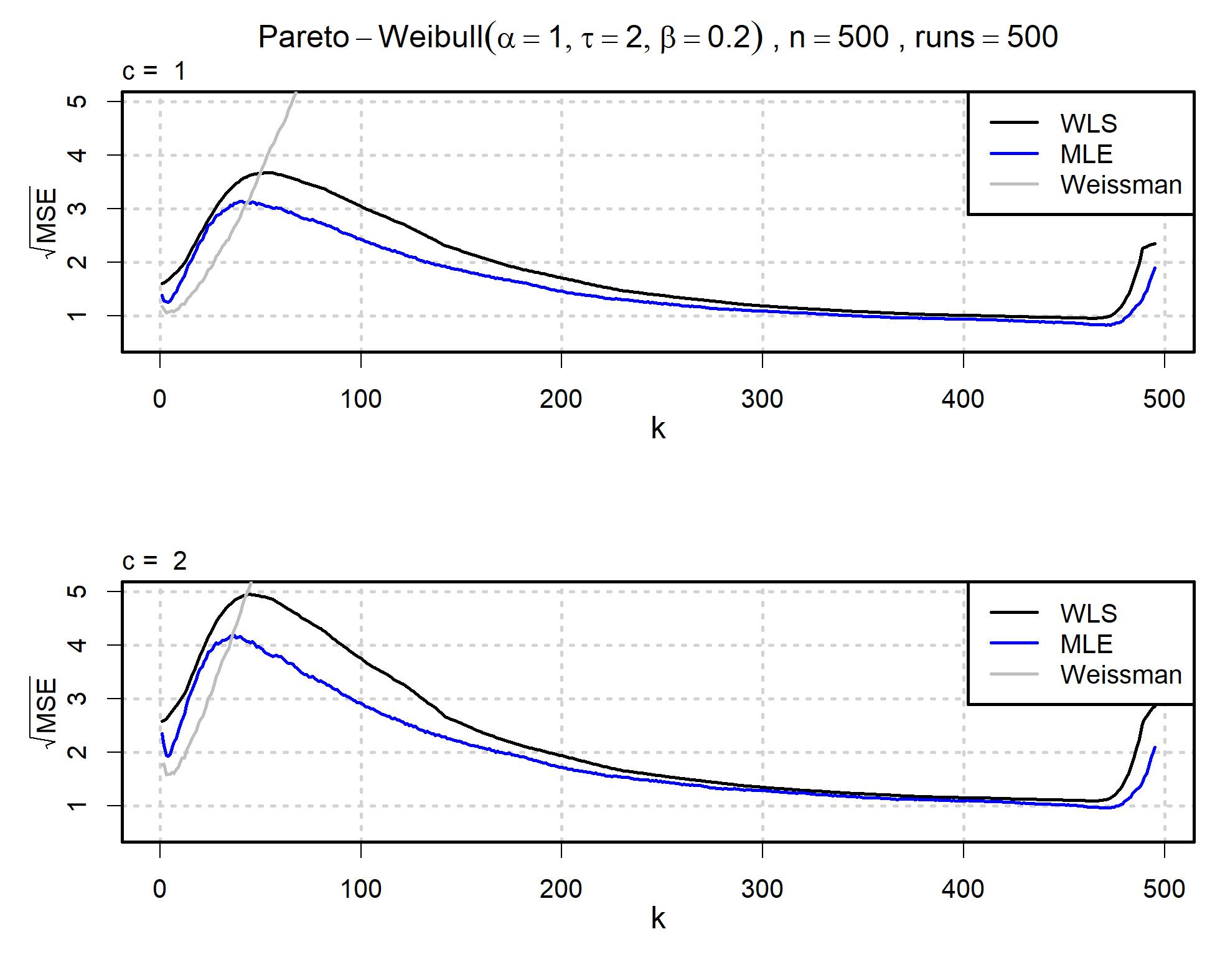}} \\ 
	\subfloat{\includegraphics[width=1\textwidth]{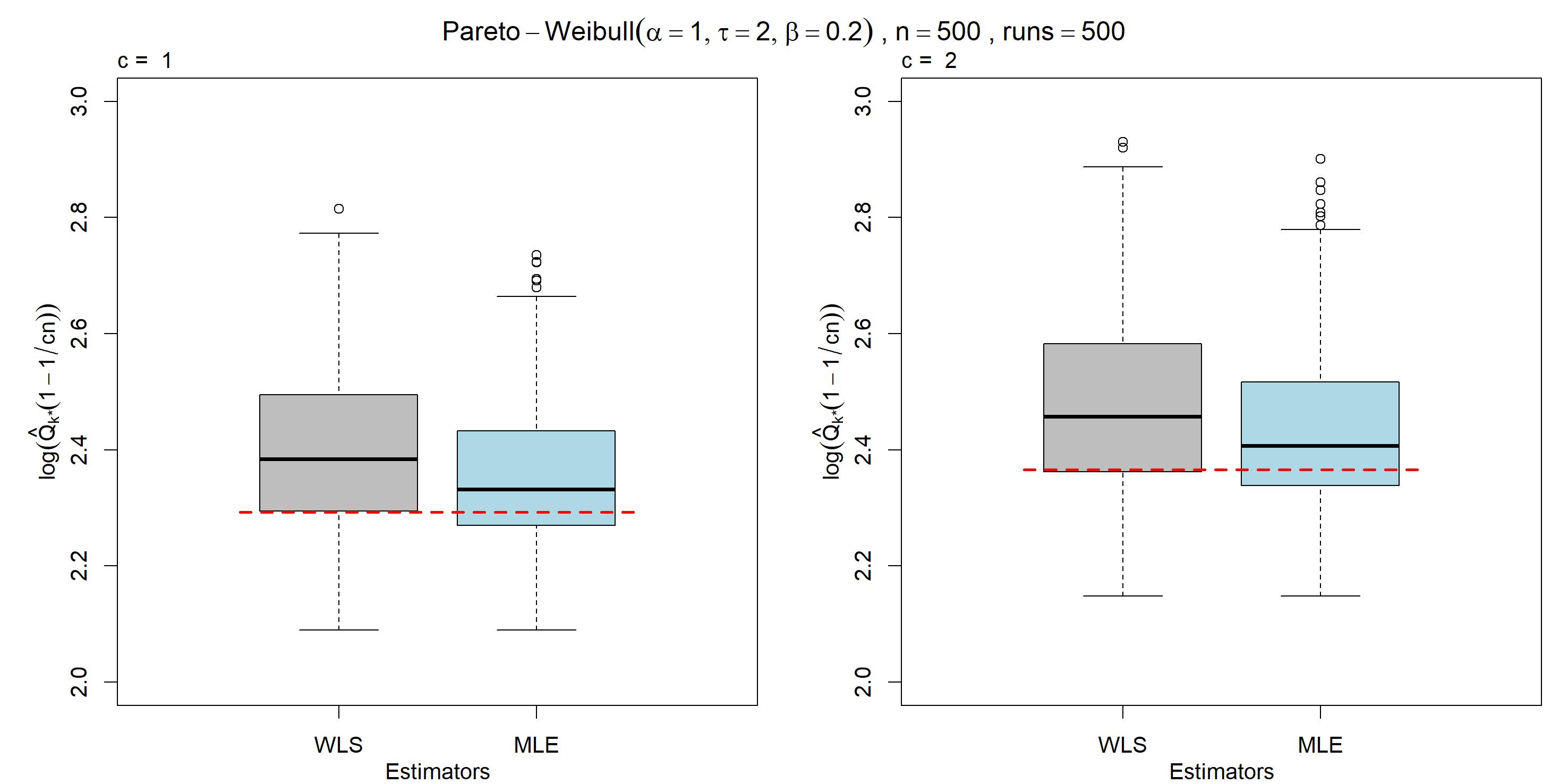}}
	\caption{Pareto-Weibull($1.0, 2.0, 0.2$): quantile estimates $\widehat{Q}^W_{p,k}$, $\widehat{Q}^M_{p,k}$ with $p={1 \over cn}$ with $c=1$ (top) and $c=2$ (middle).  Means  (left) and RMSE (right) as a function of $k$. Bottom line: boxplots of $\widehat{Q}^W_{p,\hat k}$, $\widehat{Q}^M_{p,\hat k}$ with $c=1$ (left) and $c=2$ (right). Horizontal dashed lines indicate the real parameters.}
	\label{SSs_ParetoWeibull_Qk_ns3}
\end{figure}

\begin{figure}[ht]
\centering
\subfloat{\includegraphics[width=0.5\textwidth, height = 0.5\textheight]{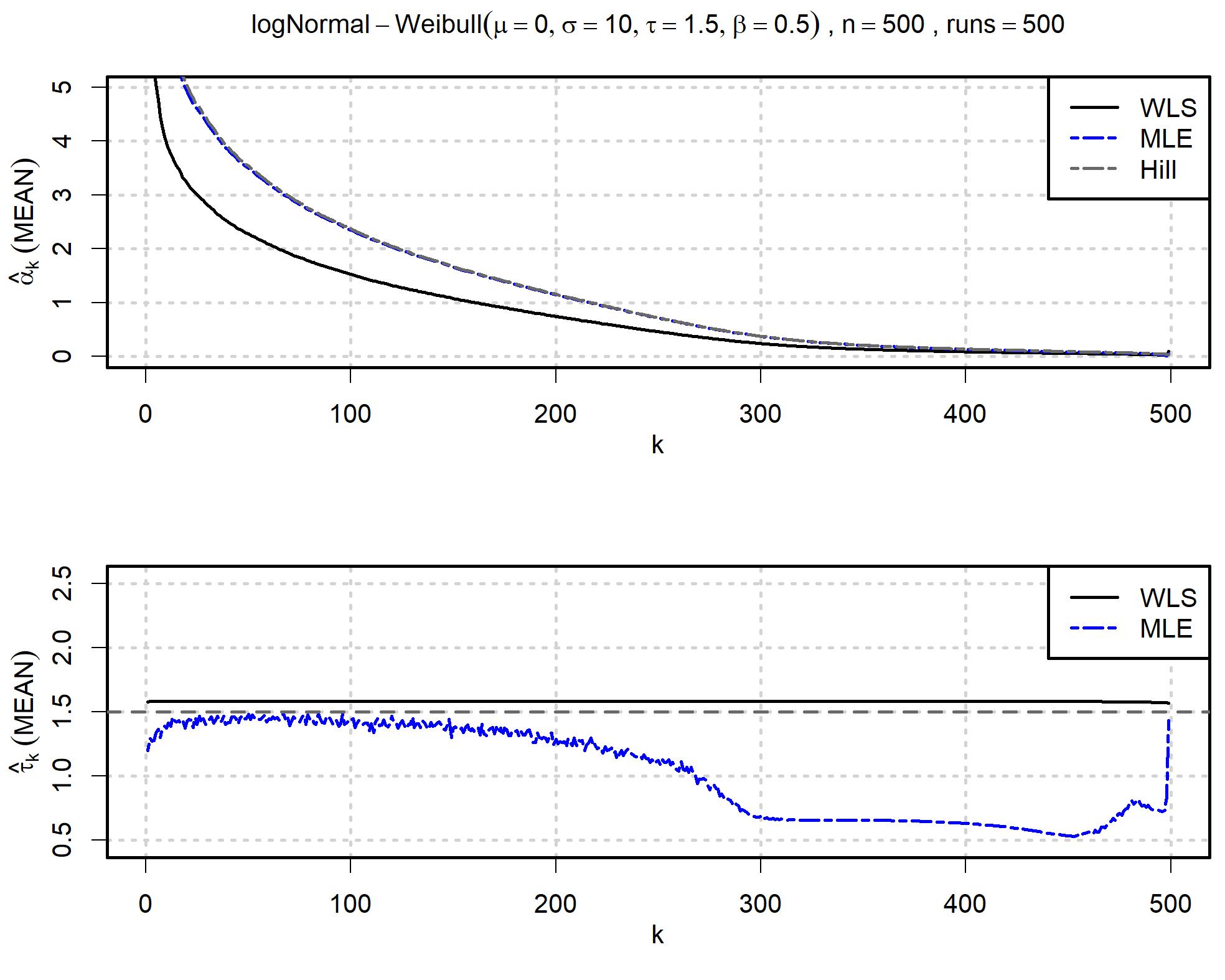}}
\subfloat{\includegraphics[width=0.5\textwidth, height = 0.5\textheight]{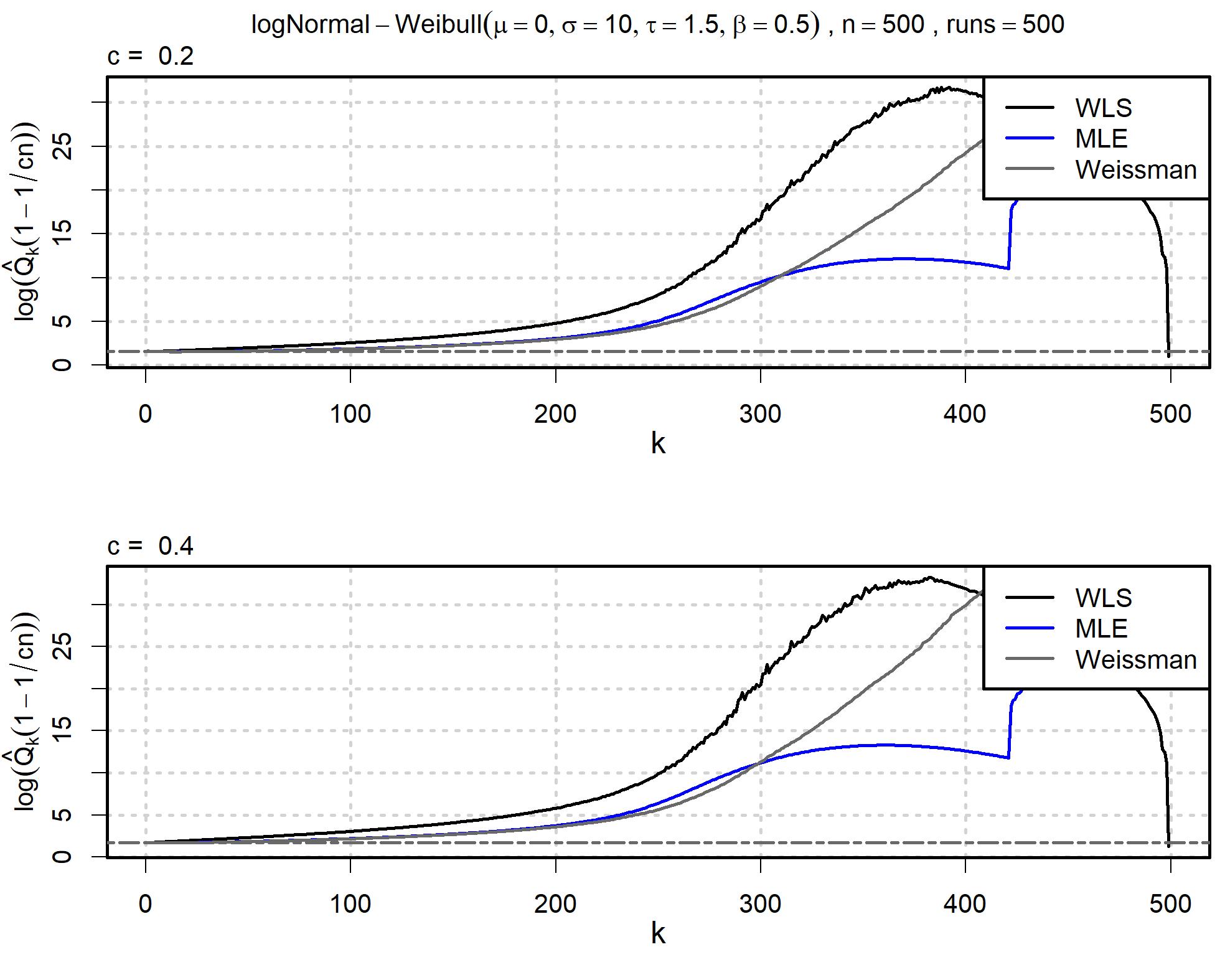}} \\
\subfloat{\includegraphics[width=0.5\textwidth]{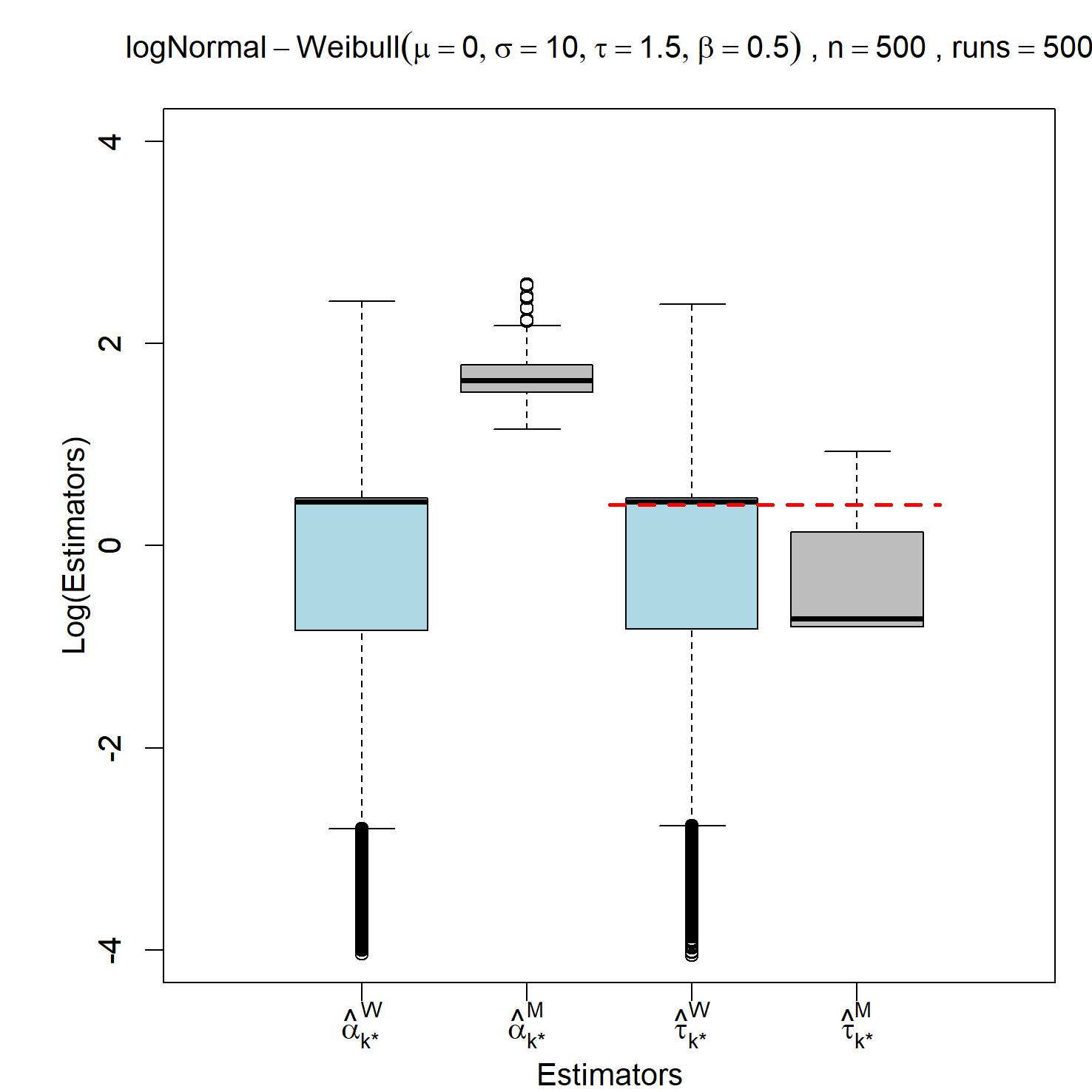}}
\caption{log-normal-Weibull($0.0, 100, 1.5, 0.5$). \textit{Top left:} $\widehat{\alpha}^W_k$, $\widehat{\alpha}^M_k$ and $H_{k,n}$ mean estimates as a function of $k$. \textit{Middle left:} $\widehat{\tau}^W_k$, $\widehat{\tau}^M_k$ mean estimates as a function of $k$. \textit{Right:} quantile estimates $\widehat{Q}^W_{p,k}$, $\widehat{Q}^M_{p,k}$ with $p={1 \over cn}$ with $c=0.2$ (top) and $c=0.4$ (middle). \textit{Bottom:} boxplots of $\widehat{\alpha}^W_{\hat k}$, $\widehat{\alpha}^M_{\hat k}$, $\widehat{\tau}^W_{\hat k}$ and $\widehat{\tau}^M_{\hat k}$ (log-scale). Horizontal dashed lines indicate the real parameters.}
\label{SSs_logNormalWeibull_ns1}
\end{figure}

\begin{figure}[ht]
\centering 
\subfloat{\includegraphics[width=1\textwidth]{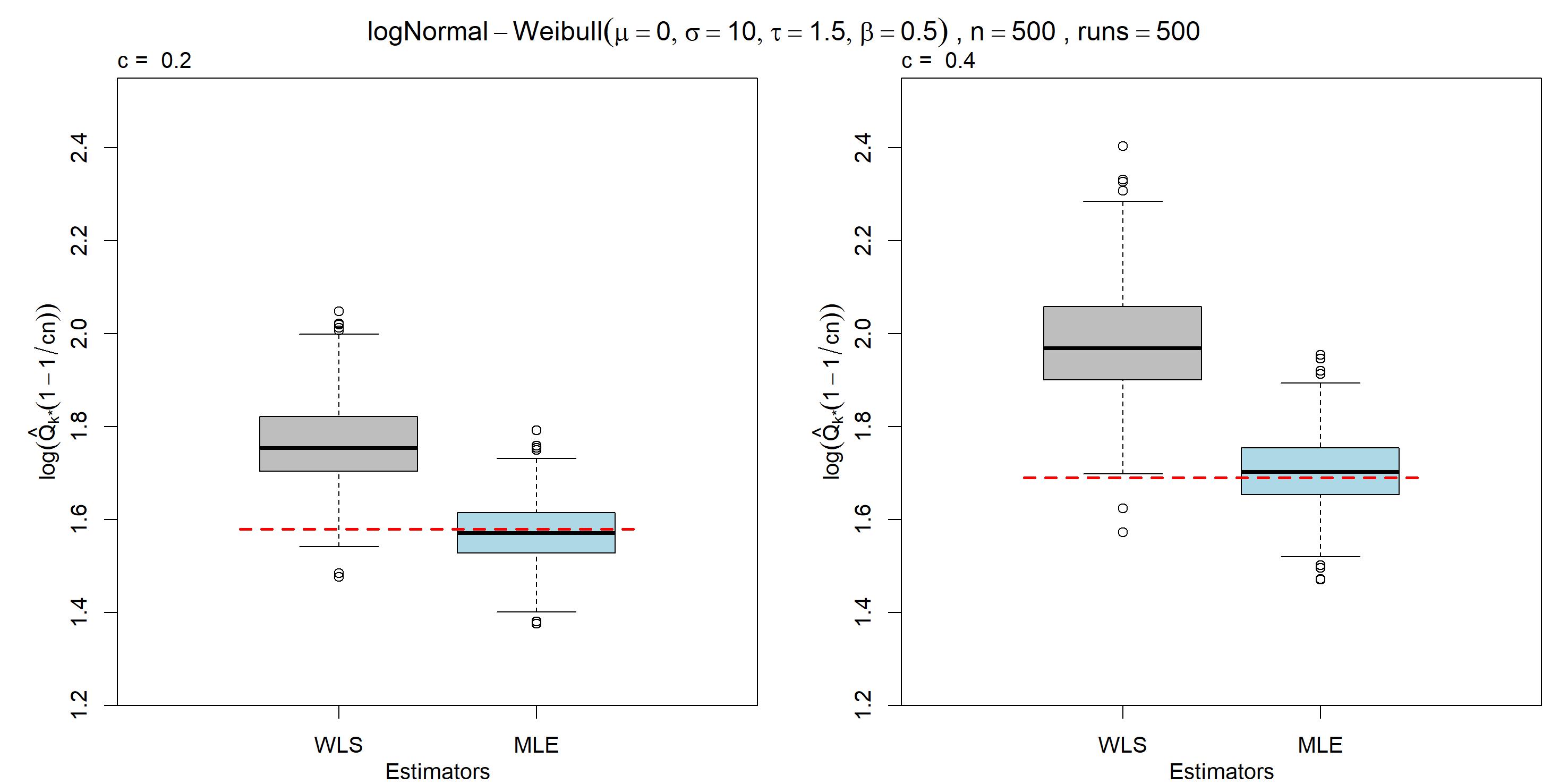}}
\caption{log-normal-Weibull($0.0, 100, 1.5, 0.5$): Boxplots of $\widehat{Q}^W_{p,\hat k}$, $\widehat{Q}^M_{p,\hat k}$ with $c=0.2$ (left) and $c=0.4$ (right). Horizontal dashed lines indicate the real values.}
\label{SSs_logNormalWeibull_Qk_ns1}
\end{figure}

\noindent Clearly the results for the MLE results $\hat\alpha ^M$, $\hat\tau ^M$, and $\hat{Q}_{p} ^M$ improve upon the weighted least squares based results. The results with the adaptive choice $\hat k$ of $k$ are promising, and again best for the MLE results. 
In case $\tau >1$ (see Figures \ref{SSs_BurrWeibull_ns1}, \ref{SSs_BurrWeibull_Qk_ns1}, \ref{SSs_FrechetWeibull_ns3}, \ref{SSs_FrechetWeibull_Qk_ns3}, \ref{SSs_ParetoWeibull_ns3}, \ref{SSs_ParetoWeibull_Qk_ns3}, \ref{SSs_logNormalWeibull_ns1} and \ref{SSs_logNormalWeibull_Qk_ns1}) when the tempering is quite strong, the results for the proposed methods are clearly improving upon the classical estimators $H_{k,n}$ and  $\hat{Q}^H_{p,k}$. Note that in these cases the $VaR$ estimates based on the MLE parameters taken at the adaptive value $\hat{k}$ show a rather small bias, even in case of the log-normal model which is situated outside our Pareto-type model assumption.

\begin{figure}[ht]
\centering
	\subfloat{\includegraphics[width=0.5\textwidth, height = 0.5\textheight]{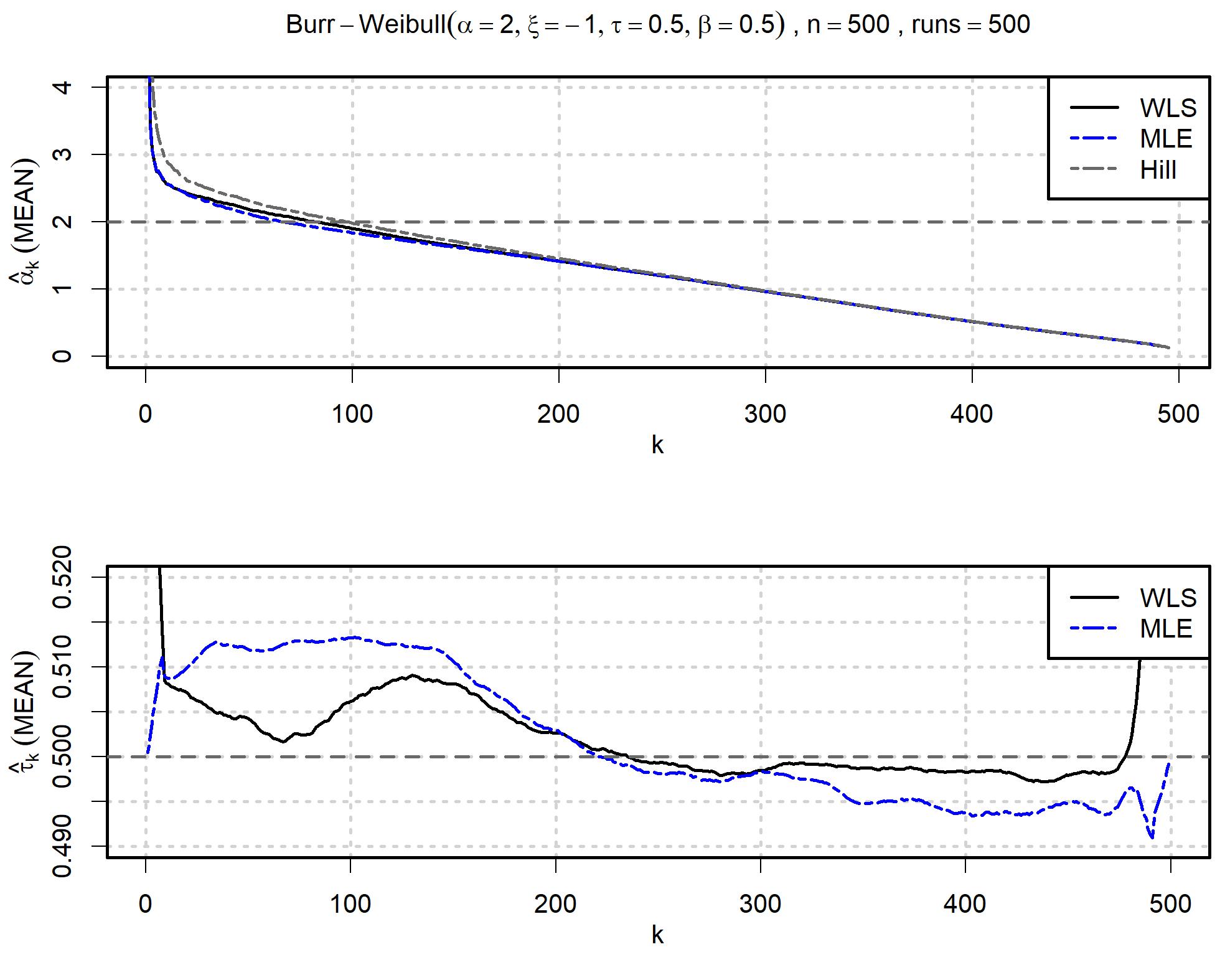}} 
	\subfloat{\includegraphics[width=0.5\textwidth, height = 0.5\textheight]{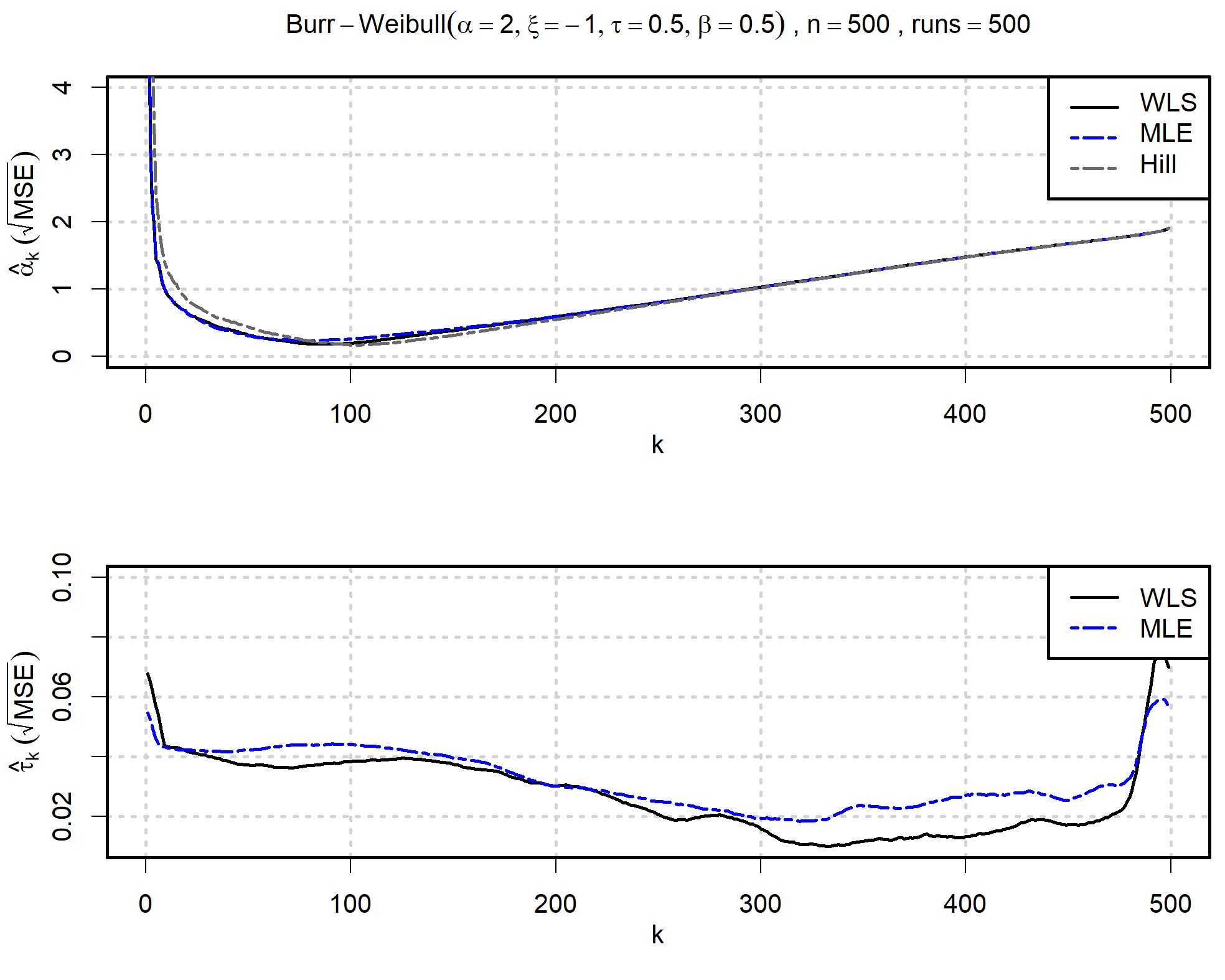}} \\
	\subfloat{\includegraphics[width=0.5\textwidth]{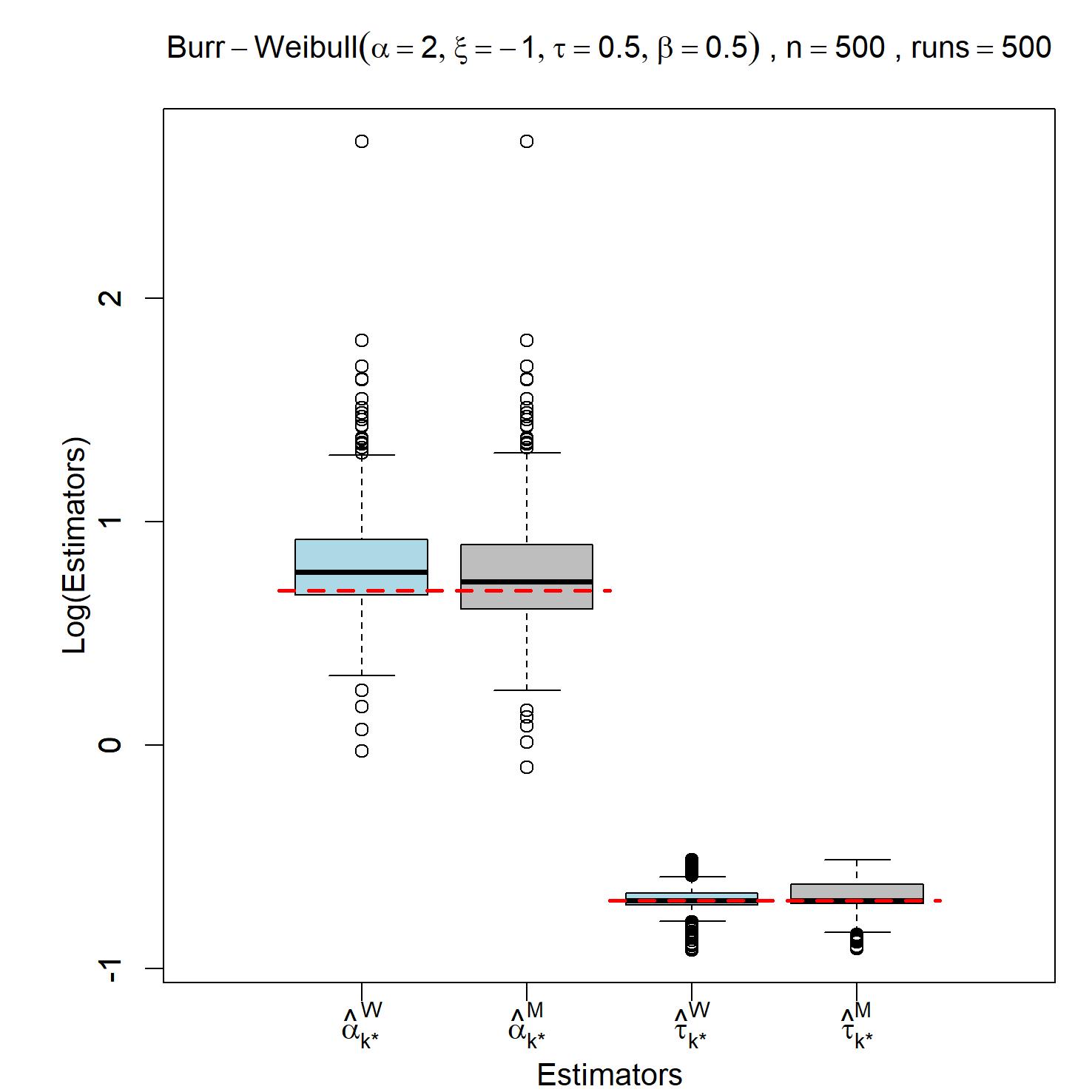}}
	\caption{Burr-Weibull($2.0, -1.0, 0.5, 0.5$). \textit{Top:} Mean (left) and RMSE (right) of  $\widehat{\alpha}^W_k$, $\widehat{\alpha}^M_k$ and $H_{k,n}$  as a function of $k$ ; \textit{Middle:} Mean (left) and RMSE (right) of $\widehat{\tau}^W_k$ and $\widehat{\tau}^M_k$ as a function of $k$; \textit{Bottom:} Boxplots of $\widehat{\alpha}^W_{\hat k}$, $\widehat{\alpha}^M_{\hat k}$, $\widehat{\tau}^W_{\hat k}$ and  $\widehat{\tau}^M_{\hat k}$ (log-scale). Horizontal dashed lines indicate the real parameters.}
	\label{SSs_BurrWeibull_ns2}
\end{figure}

\begin{figure}[ht]
\centering
	\subfloat{\includegraphics[width=0.5\textwidth, height = 0.5\textheight]{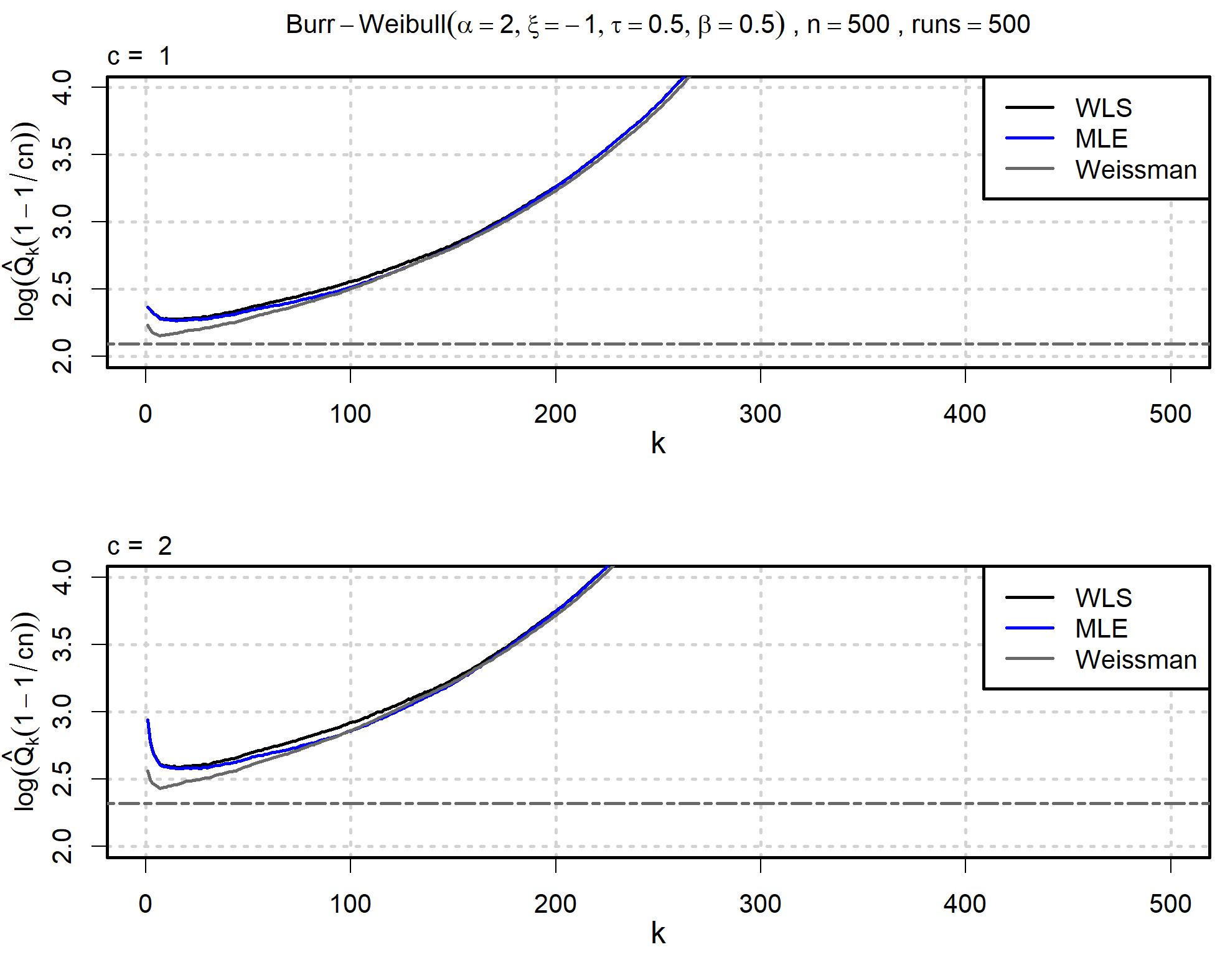}}
	\subfloat{\includegraphics[width=0.5\textwidth, height = 0.5\textheight]{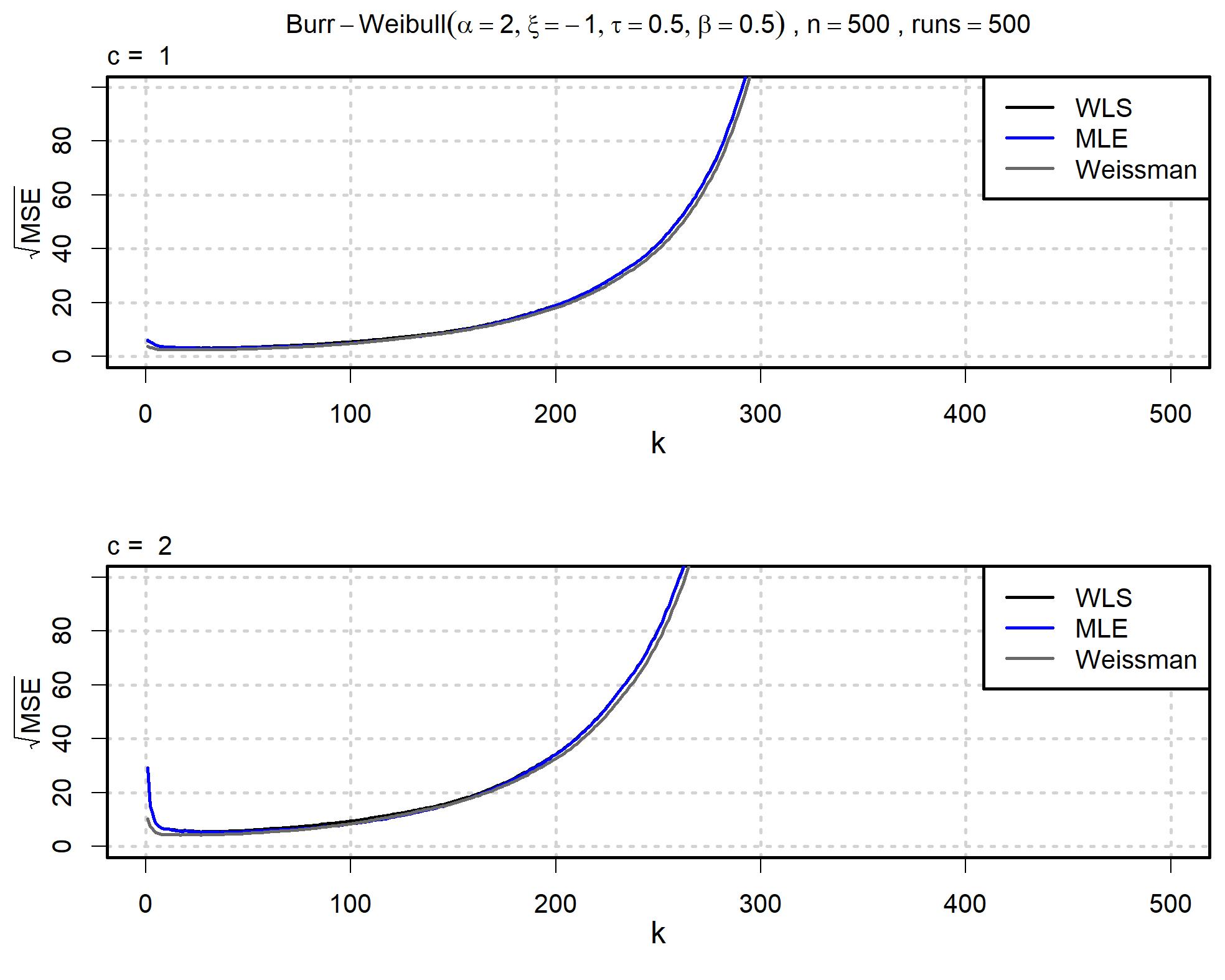}} \\
	\subfloat{\includegraphics[width=1\textwidth]{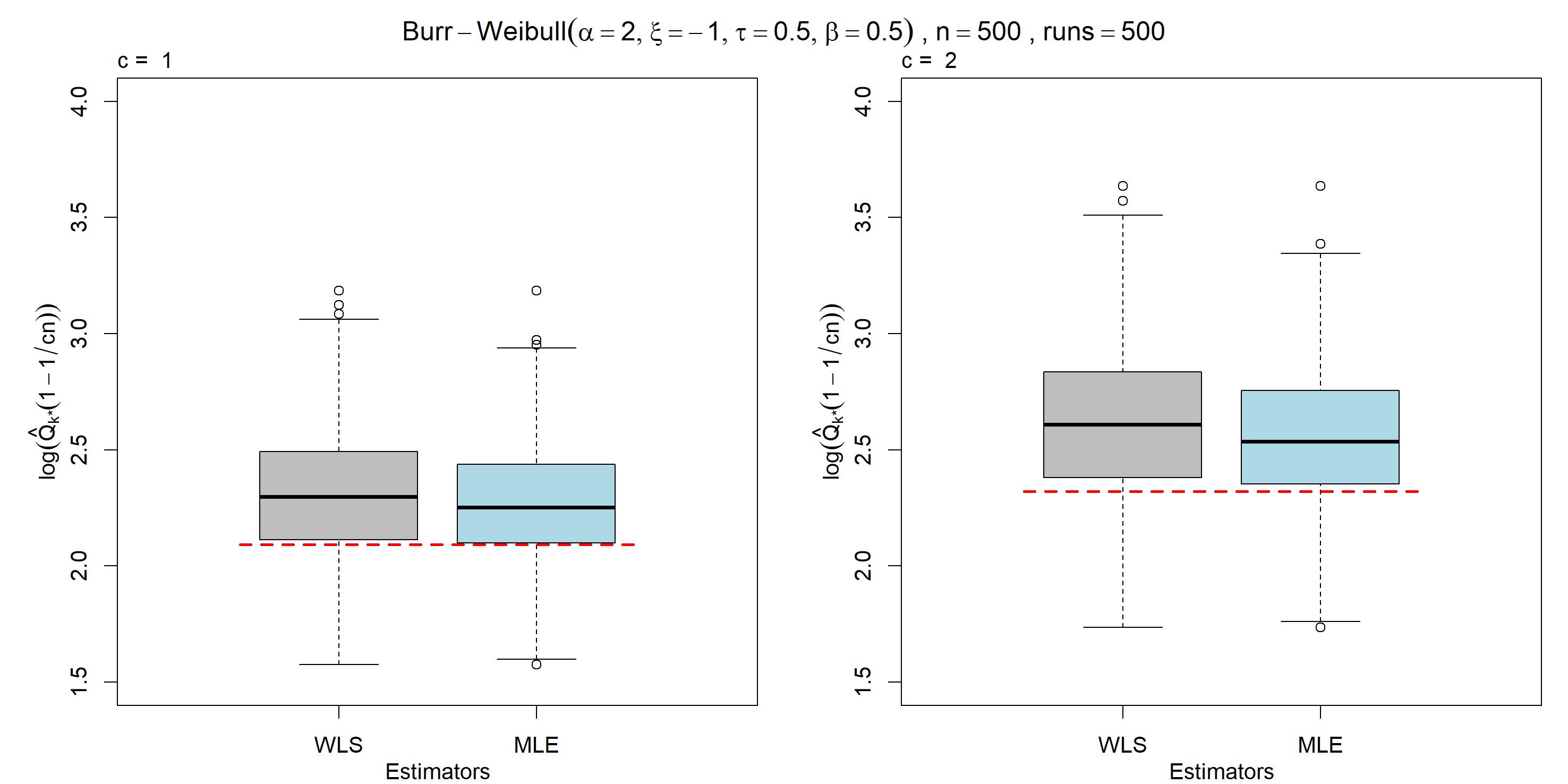}}
	\caption{Burr-Weibull($2.0, -1.0, 0.5, 0.5$):
quantile estimates $\widehat{Q}^W_{p,k}$, $\widehat{Q}^M_{p,k}$ with $p={1 \over cn}$ with $c=1$ (top) and $c=2$ (middle).  Means  (left) and RMSE (right) as a function of $k$. Bottom line: boxplots of $\widehat{Q}^W_{p,\hat k}$, $\widehat{Q}^M_{p,\hat k}$ with $c=1$ (left) and $c=2$ (right). Horizontal dashed lines indicate the real parameters.	}
	\label{SSs_BurrWeibull_Qk_ns2}
\end{figure}

\begin{figure}[ht]
\centering
	\subfloat{\includegraphics[width=0.5\textwidth, height = 0.5\textheight]{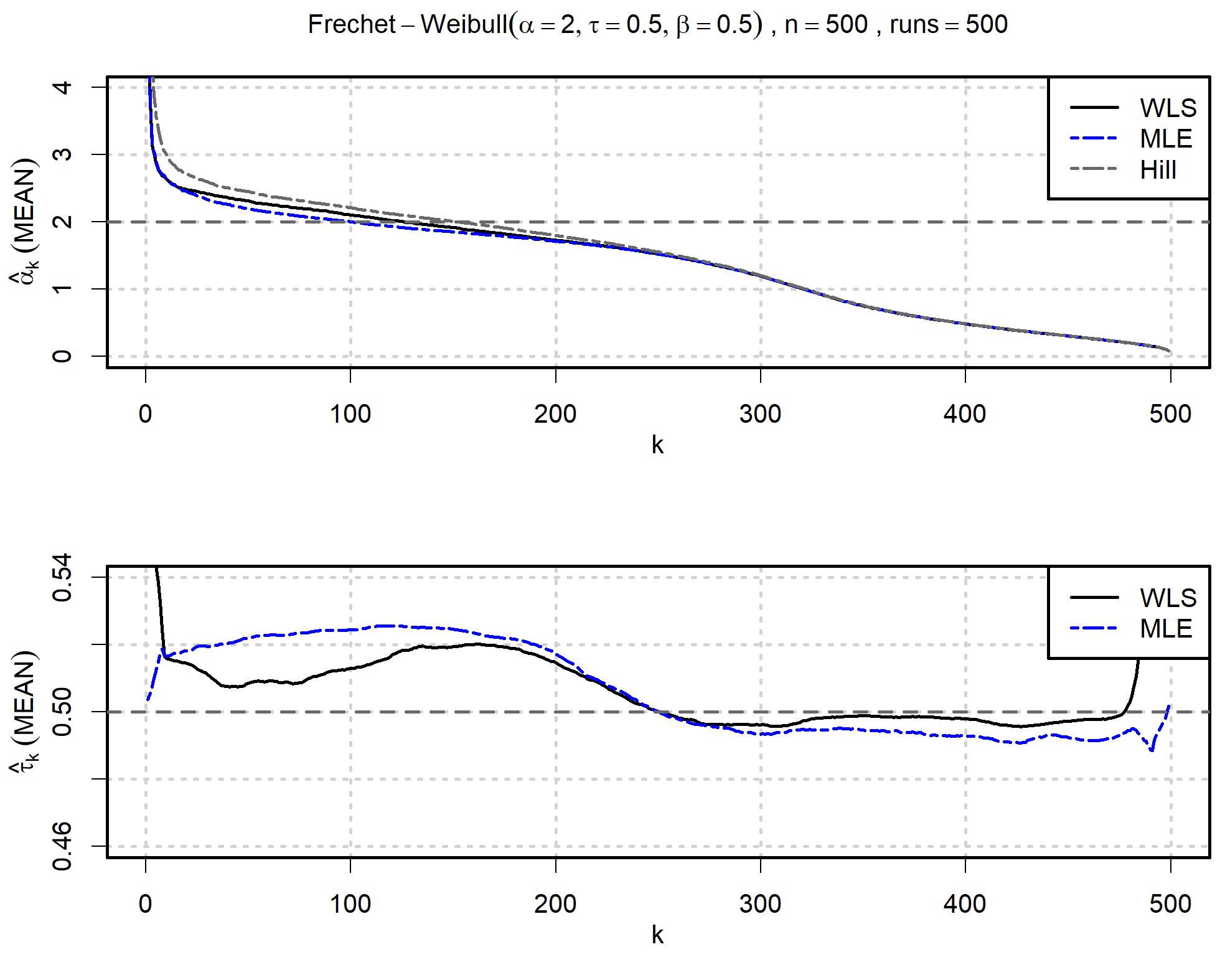}} 
	\subfloat{\includegraphics[width=0.5\textwidth, height = 0.5\textheight]{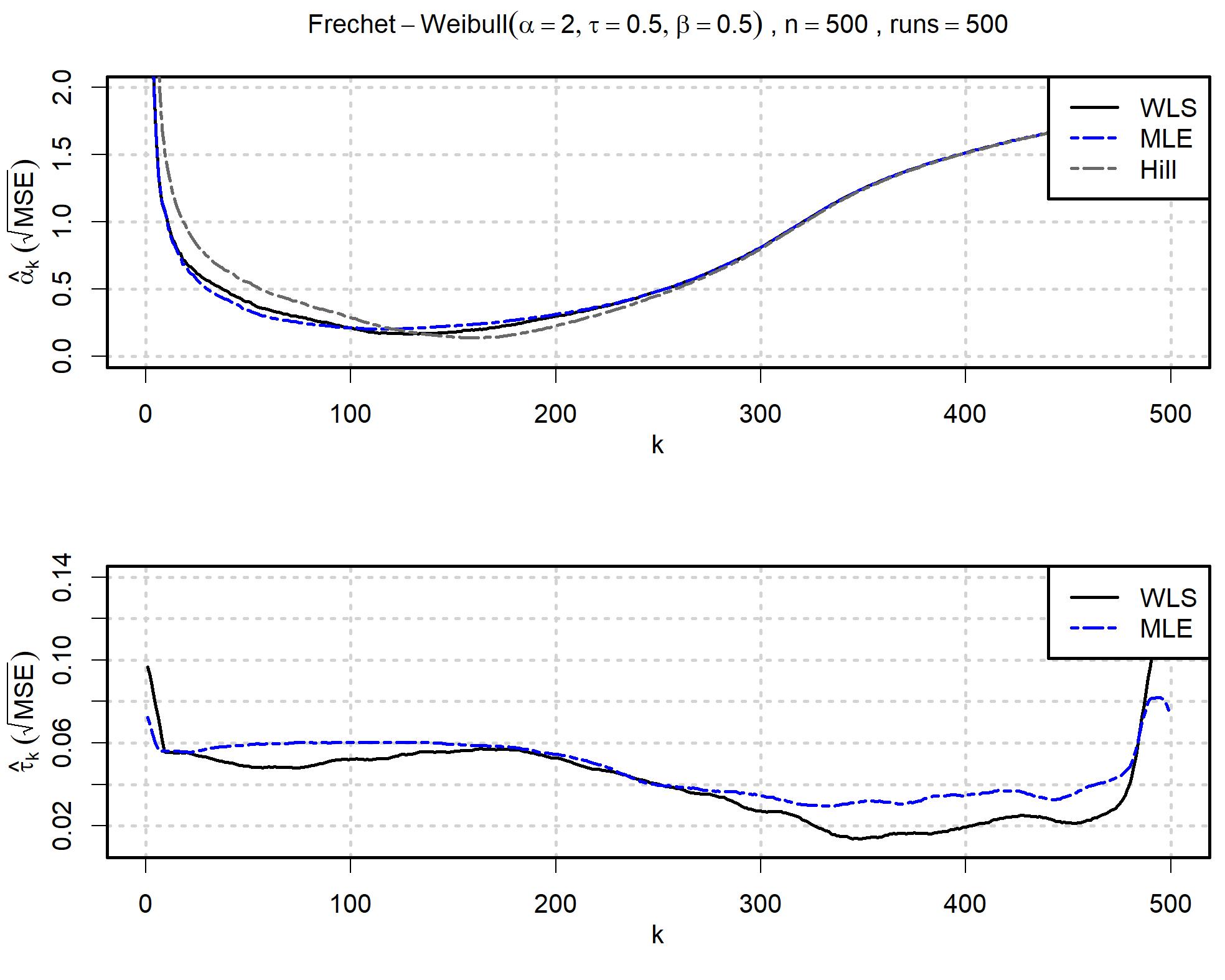}} \\
	\subfloat{\includegraphics[width=0.5\textwidth]{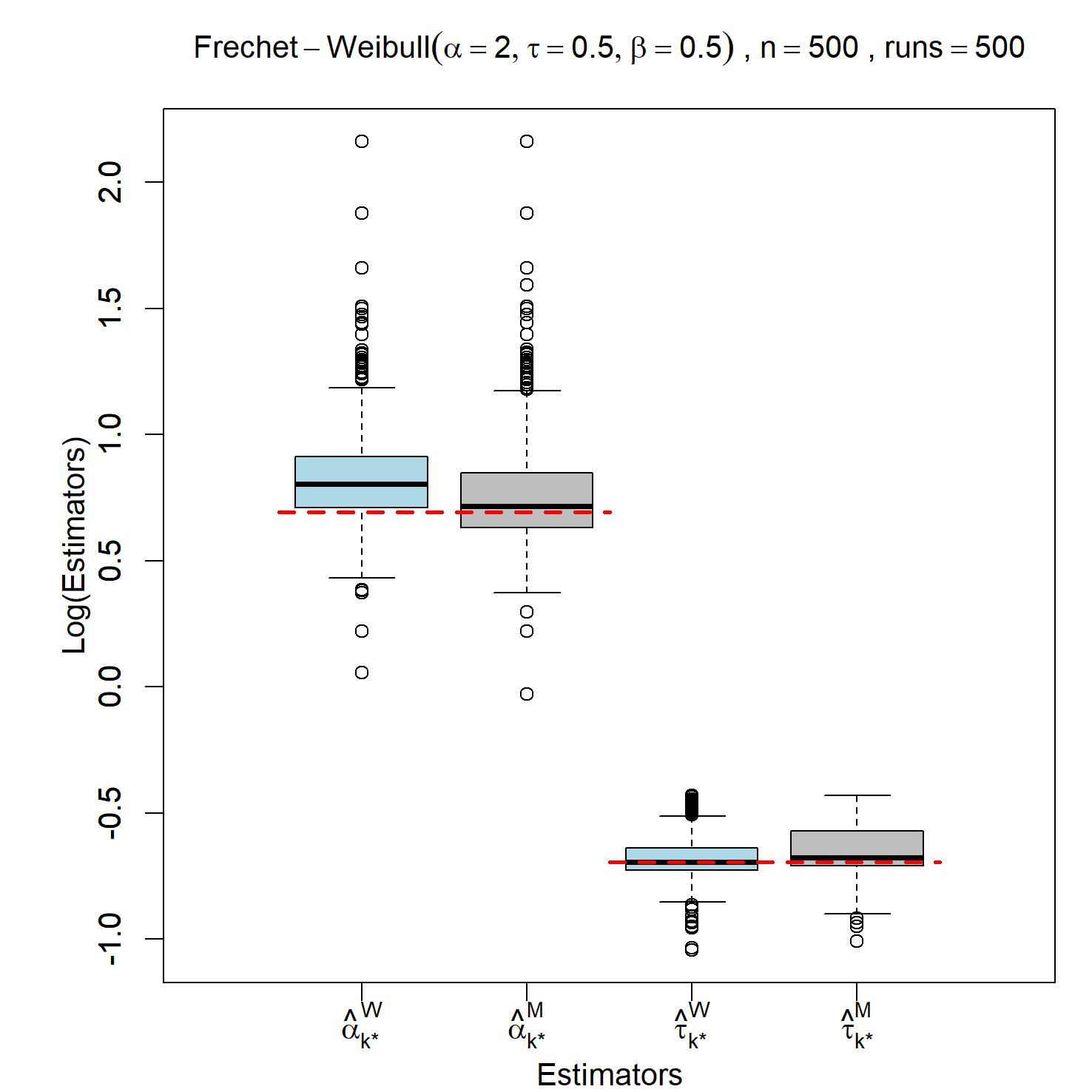}}
	\caption{Fr\'echet-Weibull($2.0, 0.50, 0.5$). \textit{Top:} Mean (left) and RMSE (right) of  $\widehat{\alpha}^W_k$, $\widehat{\alpha}^M_k$ and $H_{k,n}$ as a function of $k$ ; \textit{Middle:} Mean (left) and RMSE (right) of $\widehat{\tau}^W_k$ and $\widehat{\tau}^M_k$ as a function of $k$; \textit{Bottom:} Boxplots of $\widehat{\alpha}^W_{\hat k}$, $\widehat{\alpha}^M_{\hat k}$, $\widehat{\tau}^W_{\hat k}$ and  $\widehat{\tau}^M_{\hat k}$ (log-scale). Horizontal dashed lines indicate the real parameters.}
	\label{SSs_FrechetWeibull_ns1}
\end{figure}

\begin{figure}[ht]
\centering
	\subfloat{\includegraphics[width=0.5\textwidth, height = 0.5\textheight]{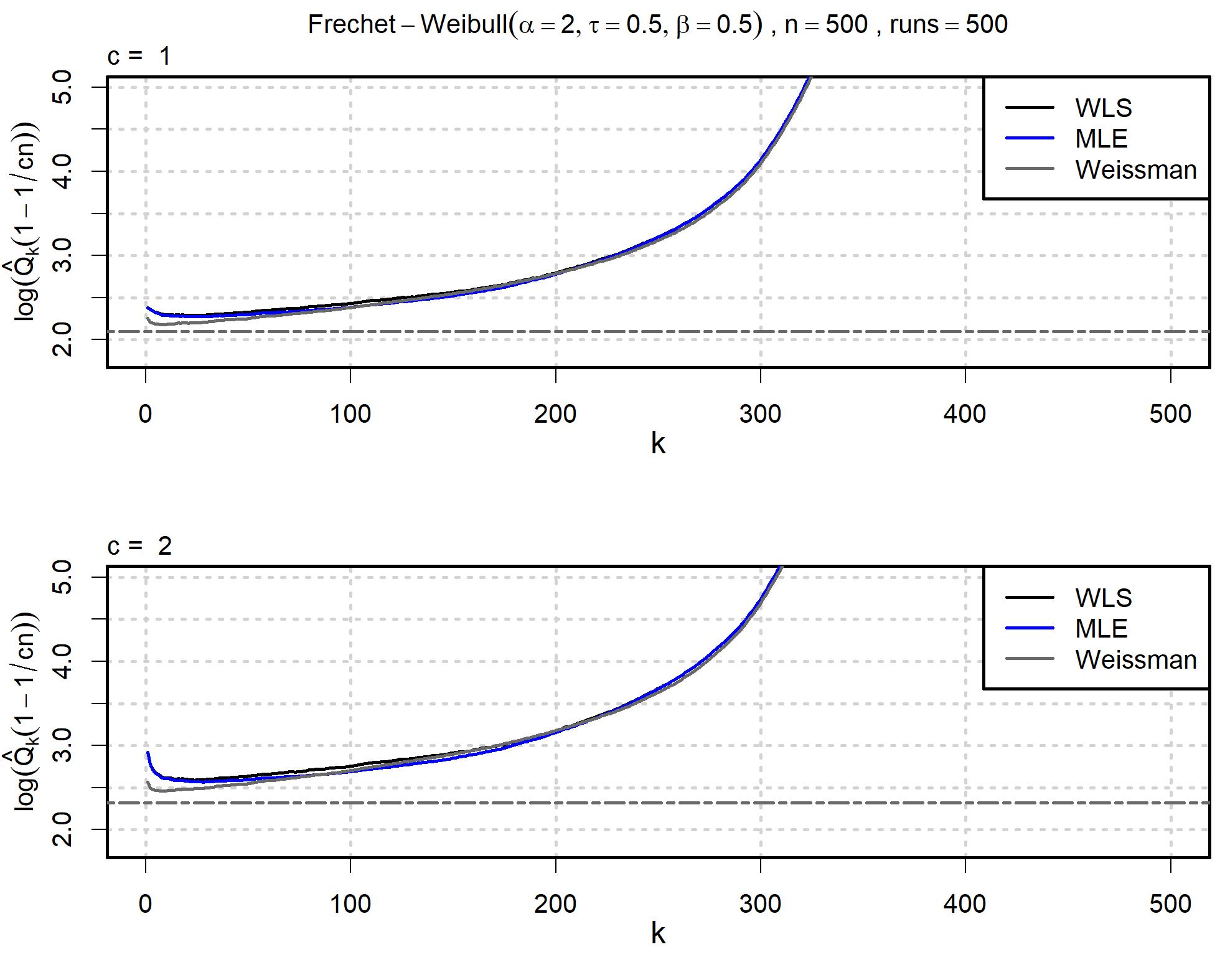}}
	\subfloat{\includegraphics[width=0.5\textwidth, height = 0.5\textheight]{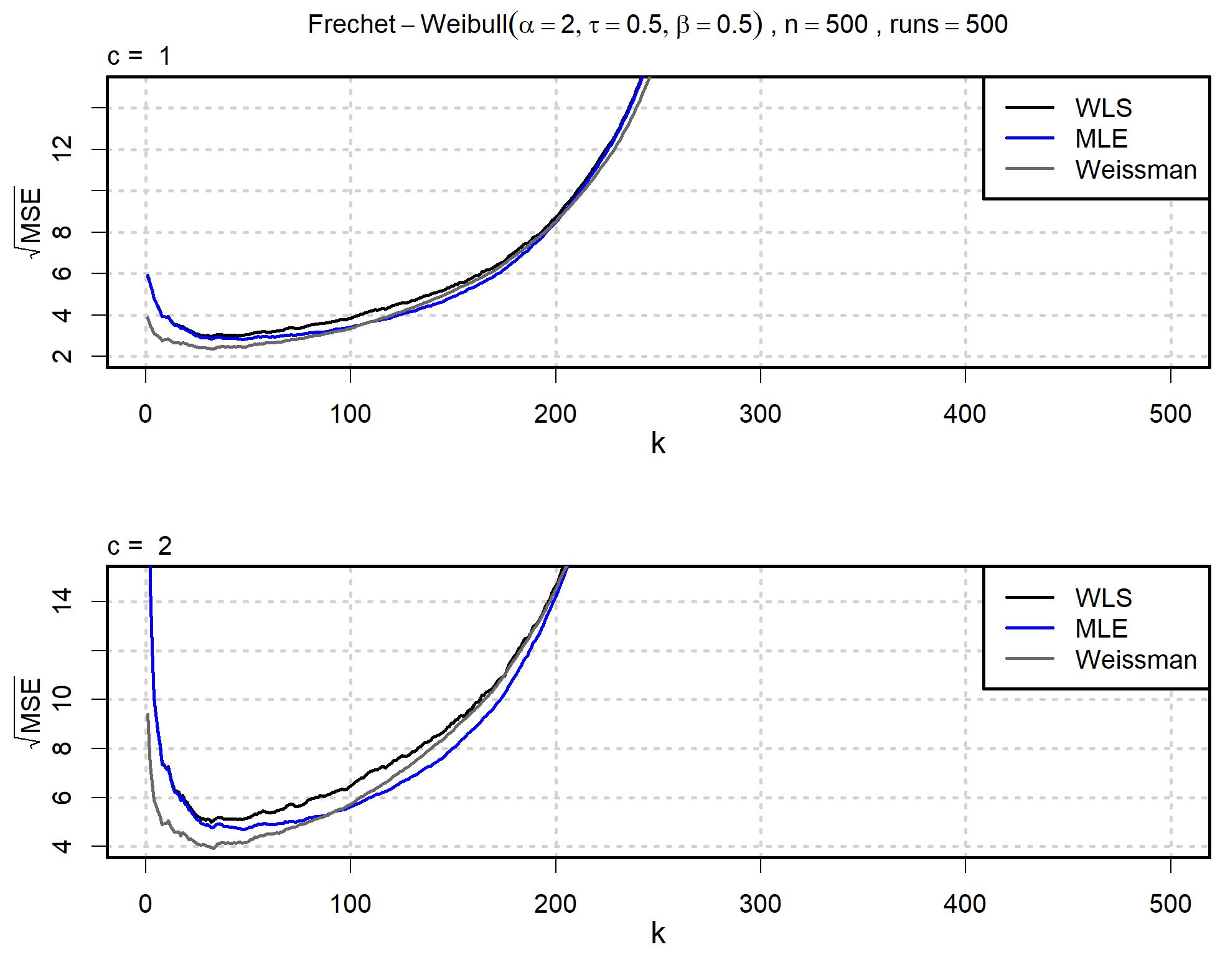}} \\
	\subfloat{\includegraphics[width=1\textwidth]{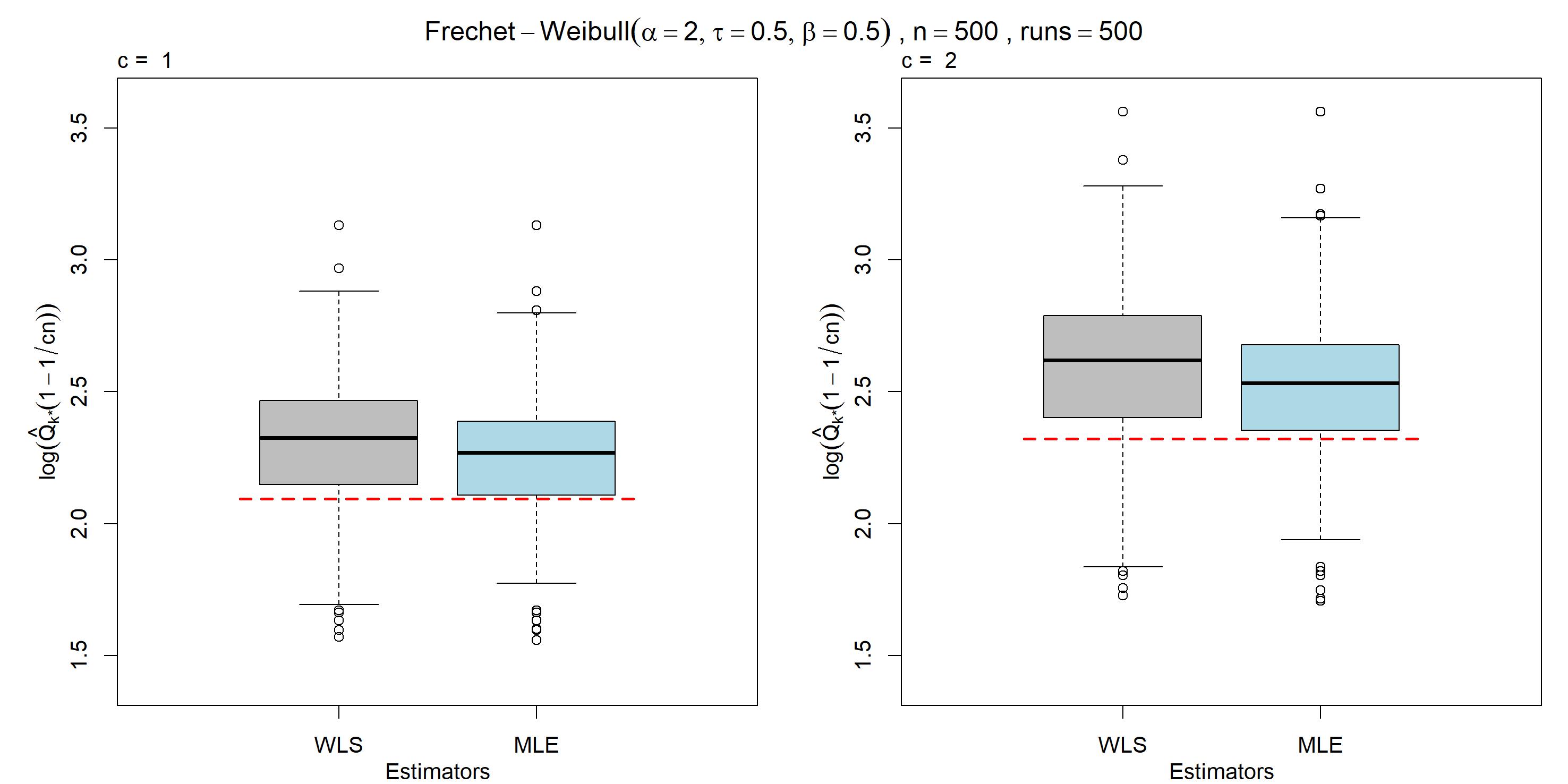}}
	\caption{Fr\'echet-Weibull($2.0, 0.5, 0.5$):
quantile estimates $\widehat{Q}^W_{p,k}$, $\widehat{Q}^M_{p,k}$ with $p={1 \over cn}$ with $c=1$ (top) and $c=2$ (middle).  Means  (left) and RMSE (right) as a function of $k$. Bottom line: boxplots of $\widehat{Q}^W_{p,\hat k}$, $\widehat{Q}^M_{p,\hat k}$ with $c=1$ (left) and $c=2$ (right). Horizontal dashed lines indicate the real parameters.	}
	\label{SSs_FrechetWeibull_Qk_ns1}
\end{figure}

In case $\tau <1$ (see Figures \ref{SSs_BurrWeibull_ns2}, \ref{SSs_BurrWeibull_Qk_ns2}, \ref{SSs_FrechetWeibull_ns1} and \ref{SSs_FrechetWeibull_Qk_ns1}),  hence under weaker tempering,  the bias and RMSE results are comparable with the classical estimators. The $VaR$ estimates at $\hat k$ tend to overestimate the correct value. As will become clear from the case studies in the next section, the Pareto and tempered Pareto  fits can lead to quite different extreme tail fits {\it per sample}.

\noindent We conclude  that the use of  classical estimators ignoring the tempering effect leads to serious overestimation of the risk measures, while the proposed method provides reasonable $VaR$ estimates especially for larger values of $\tau >1$. In case of smaller tempering with a heavier Weibull tail, improvements can be made concerning the adaptive choice of $k$. Another possibility is to search for bias reduced estimators as available in the non-tempering literature (see for instance Chapters 3 and 4 in \cite{beirlant2006statistics}).

\section{Insurance cases}\label{sec4}

We now apply the presented methods to the Norwegian and the Secura Re Belgian data sets  introduced in Section 1. In addition, we contrast the tail index estimates $\widehat{\alpha}^W_{k,n}$ and $\widehat{\alpha}^M_{k,n}$ with the values obtained for the truncated Pareto-type model proposed in Beirlant et al.\ \cite{beirlant2016tail}, where $\widehat{\alpha}^T_{k,n}$ is obtained as the solution to
\begin{align*}
H_{k,n} = \frac{1}{\alpha^T_{k,n}} + \frac{R_{k,n}^{\alpha^T_{k,n}} \log(R_{k,n})}{1-R_{k,n}^{\alpha^T_{k,n}},}
\end{align*}
with $R_{k,n} = X_{n-k,n} / X_{n,n}.$ The latter estimator was first proposed in Aban et al.\ \cite{aban2006parameter} as the conditional MLE based on the $k+1 ~(0\leq k <n)$ largest order statistics representing only the portion of the tail where the truncated Pareto approximation holds, see also \cite[Sec 4.2.3]{albrecher2017reinsurance}.\\
We then also measure the goodness-of-fit using QQ-plot \eqref{eqline} and the analogous expression for the truncated model.

\begin{figure}[ht]
\centering
\includegraphics[width=0.48\textwidth]{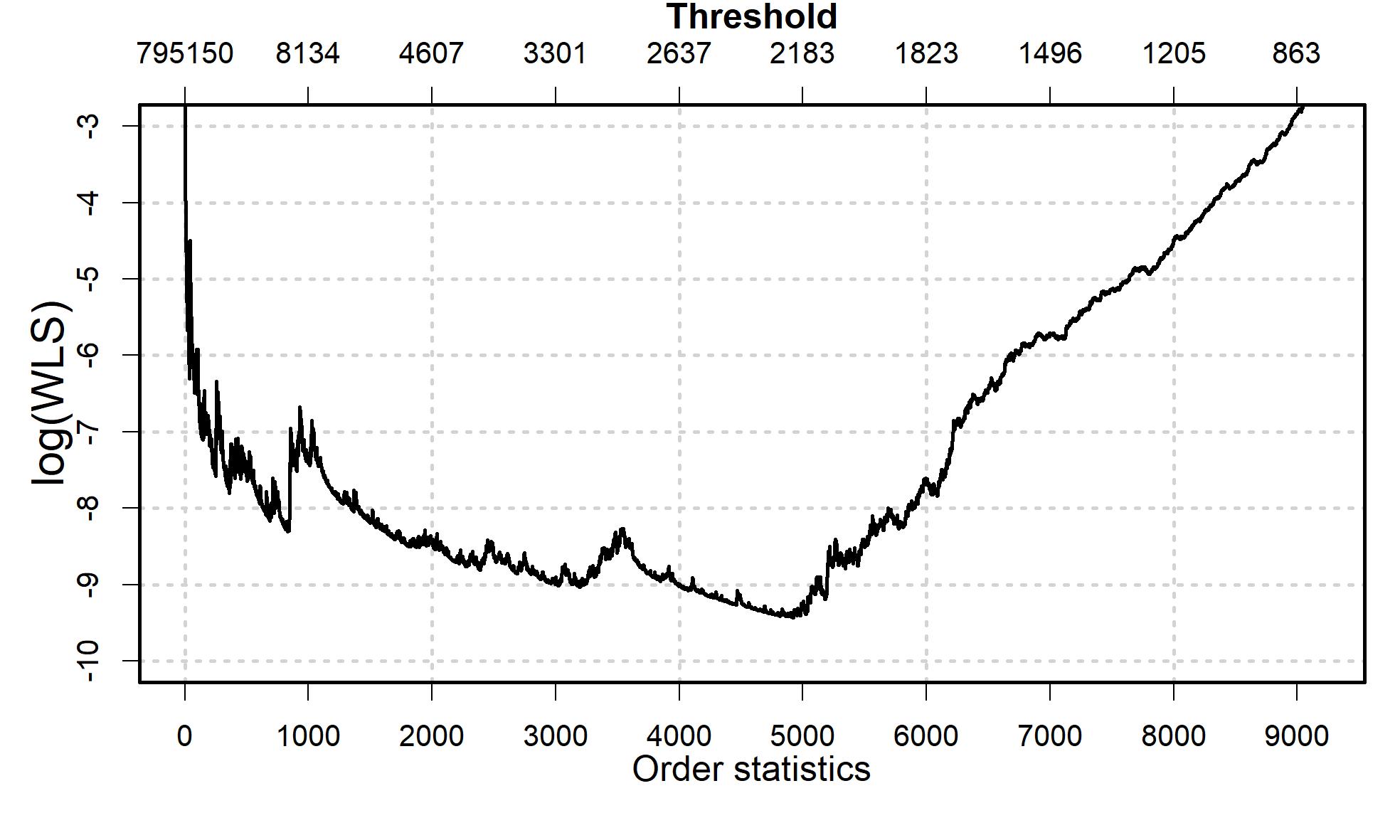}
		\subfloat{\includegraphics[width=0.5\textwidth]{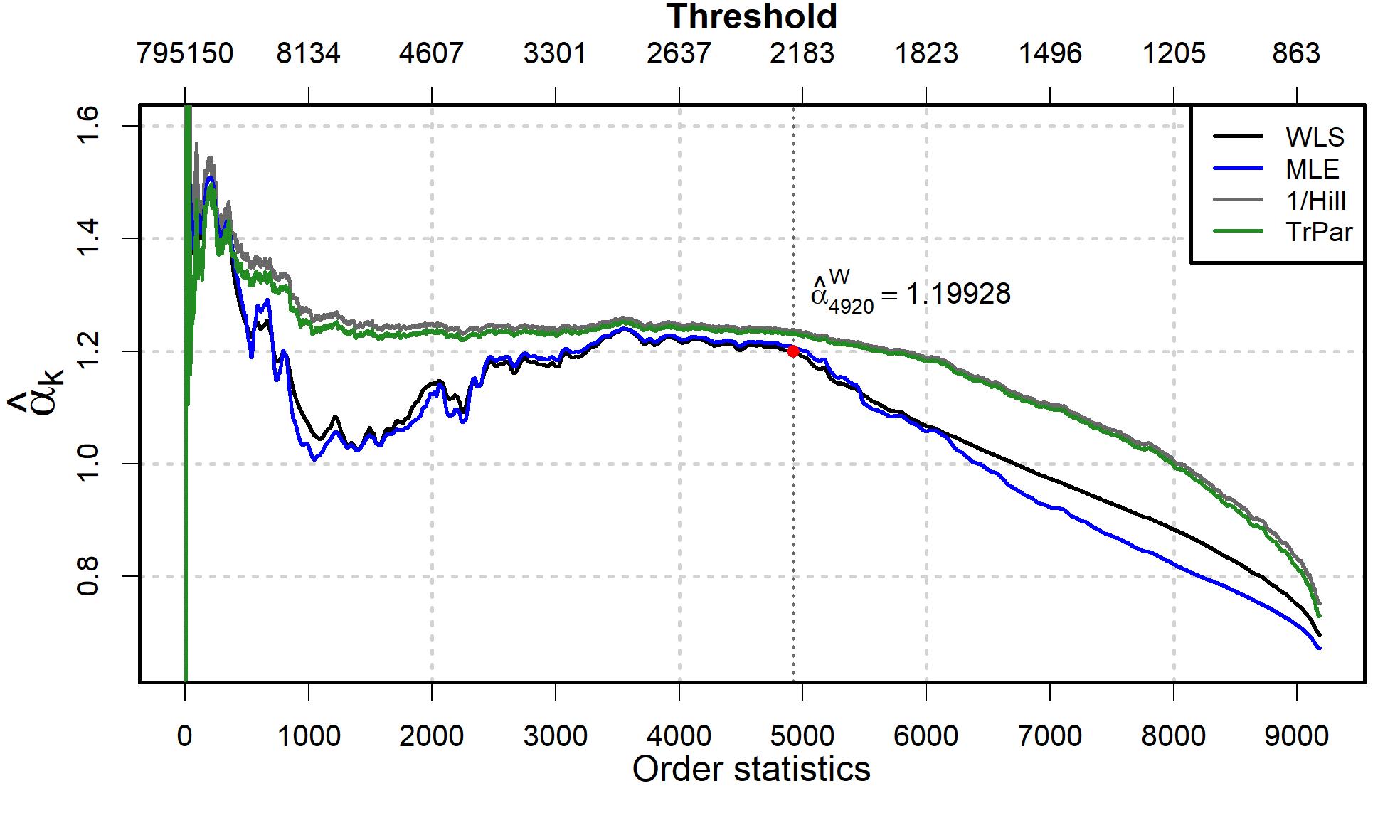}} \\
	\subfloat{\includegraphics[width=0.5\textwidth]{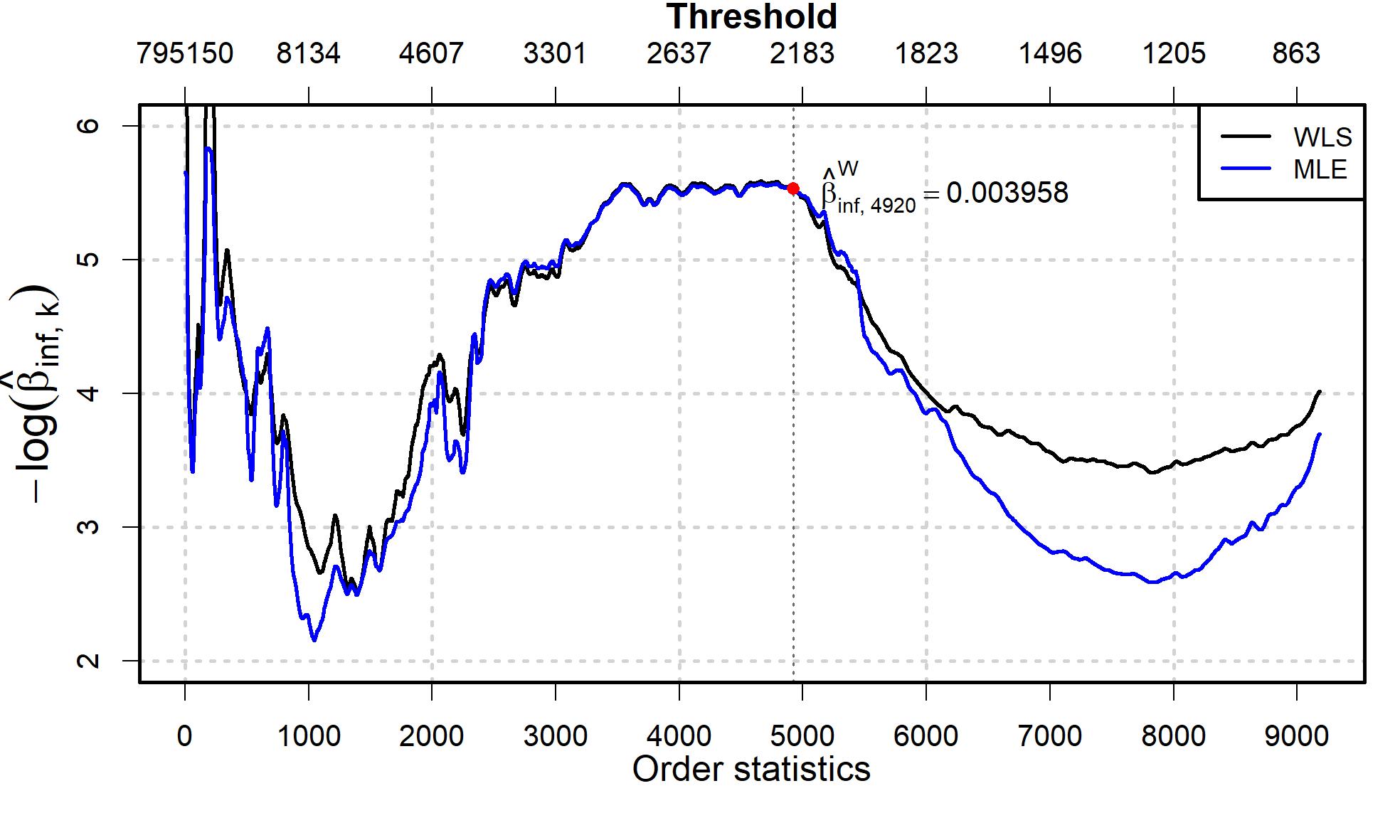}}  
	\subfloat{\includegraphics[width=0.5\textwidth]{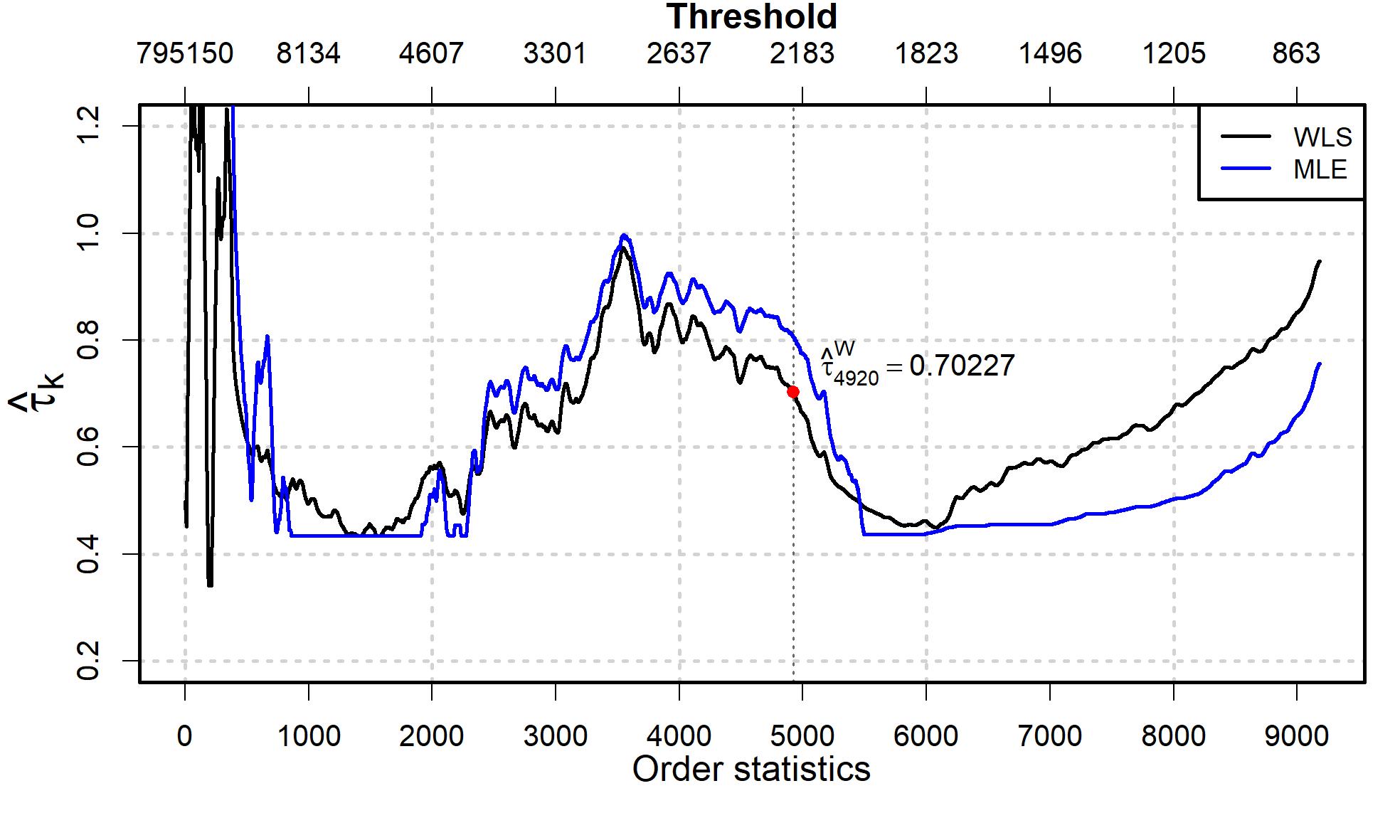}} \\
	\centering
	\subfloat{\includegraphics[width=0.70\textwidth]{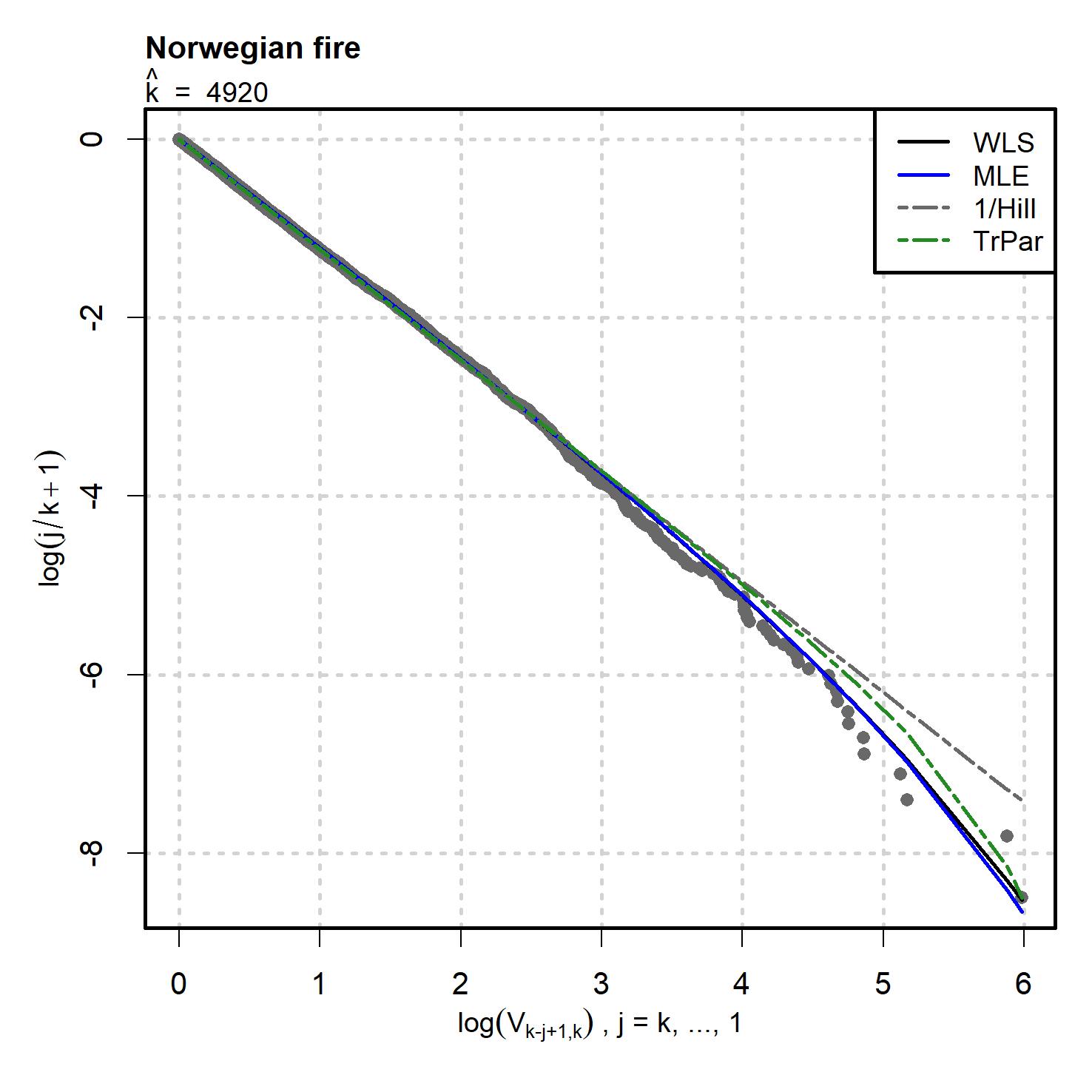}}
	\caption{\textit{Norwegian fire insurance data:} Top left: $SS_k$ from \eqref{SSk}; Top right: $\hat\alpha _k^W$, $\hat\alpha _k^M$, $H_{k,n}$ and $\hat{\alpha}_k^T$; 
	Middle left: $-\log\hat{\beta}_{\infty,k}^W$, $-\log\hat{\beta}_{\infty,k}^M$; Middle right: $\hat\tau _k^W$, $\hat\tau _k^M$; Bottom: $\log$-$\log$ plot with fit obtained from \eqref{eqline} with $k=\hat{k}=4915$ using MLE and WLS estimates, next to Pareto and truncated Pareto fit.}
	\label{Norwegian_Estimates}
\end{figure}

\begin{figure}[ht]
\centering
		\includegraphics[width=0.8\textwidth]{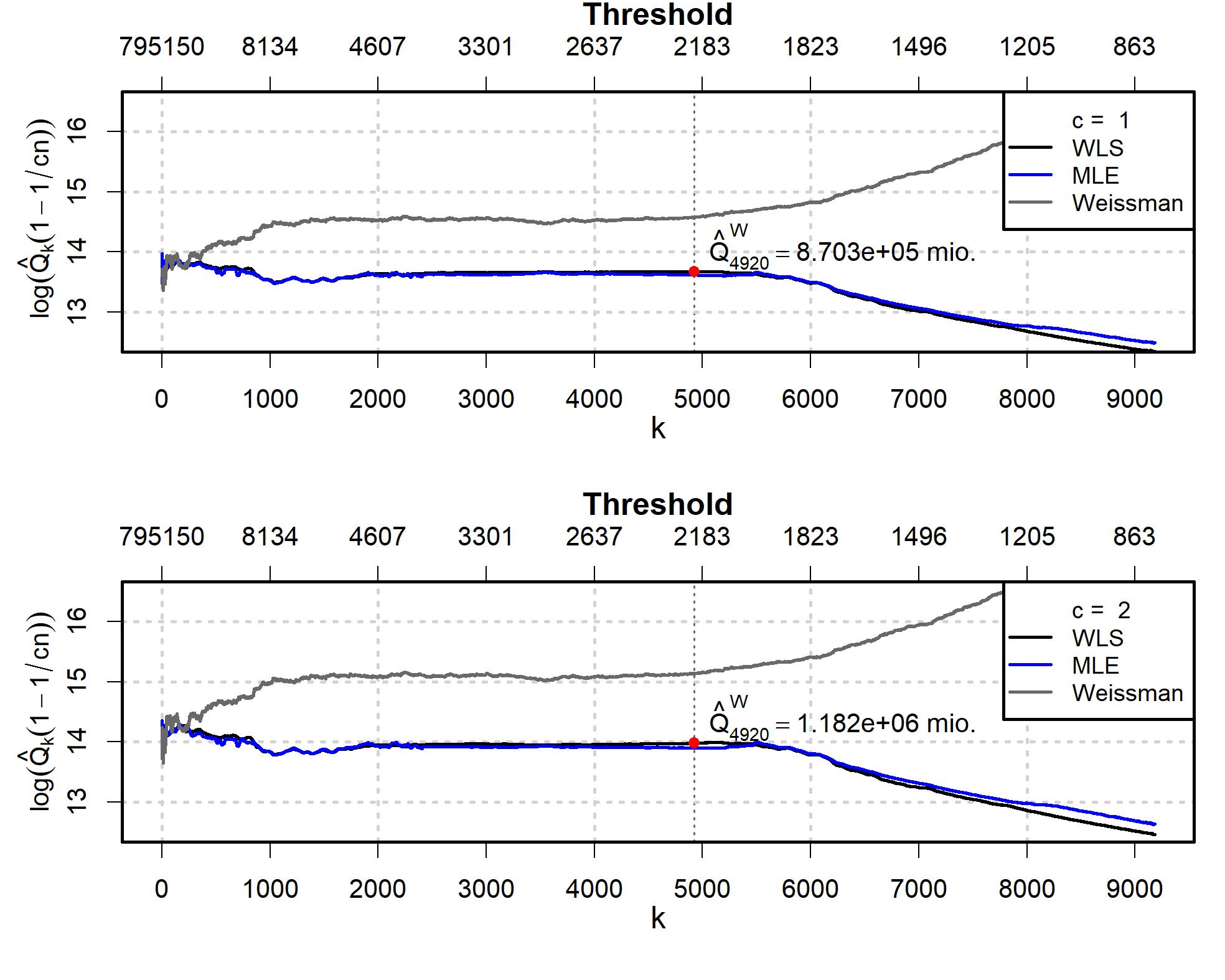}
	\caption{\textit{Norwegian fire insurance data}: $\hat{Q}^W_{p,k}$, $\hat{Q}^M_{p,k}$ and $\hat{Q}^H_{p,k}$ quantile estimates with $p=1/n$ (top) and $p=1/(2n)$ (bottom).}
	\label{Norwegian_Qestimates}
\end{figure}

For the Norwegian fire insurance data set, we find $\hat k= 4920$ from the plot of $SS_k$ from \eqref{SSk} in Figure \ref{Norwegian_Estimates}, where also the different parameter estimates as a function of $k$ can be found. The $\log$-$\log$ plot based on \eqref{eqline} at $k=4920$ shows a good tail fit for the tempered Pareto model, in contrast with the simple Pareto fit which will overfit tail probabilities and quantiles. This can be seen from  Figure \ref{Norwegian_Qestimates}
where for larger $k$, the classical Weissman estimates $\hat{Q}^H_{1/(cn),k}$ ($c=1,2$) lead to much larger estimates than those based on the proposed tempering modelling. Only when $k$ is really small, i.e.\ when restricting to the data situated in the bottom curved area of the  $\log$-$\log$ plot,  the classical linear Pareto fit is able to provide a reasonable representation of the most extreme data. Finally, note from the $\log$-$\log$ plot in Figure \ref{Norwegian_Estimates} that the truncated Pareto fit follows the linear Pareto fit except for the two final extreme points after which a sharp deviation is observed up to an estimated finite truncation point $T$ estimated at $\hat T_{\hat{k}}=1,211,106$, when using the estimation method proposed in \cite[Sec. 3, Eq. 19]{beirlant2016tail}.

\begin{figure}[ht]
\centering
\subfloat{\includegraphics[width=0.6\textwidth]{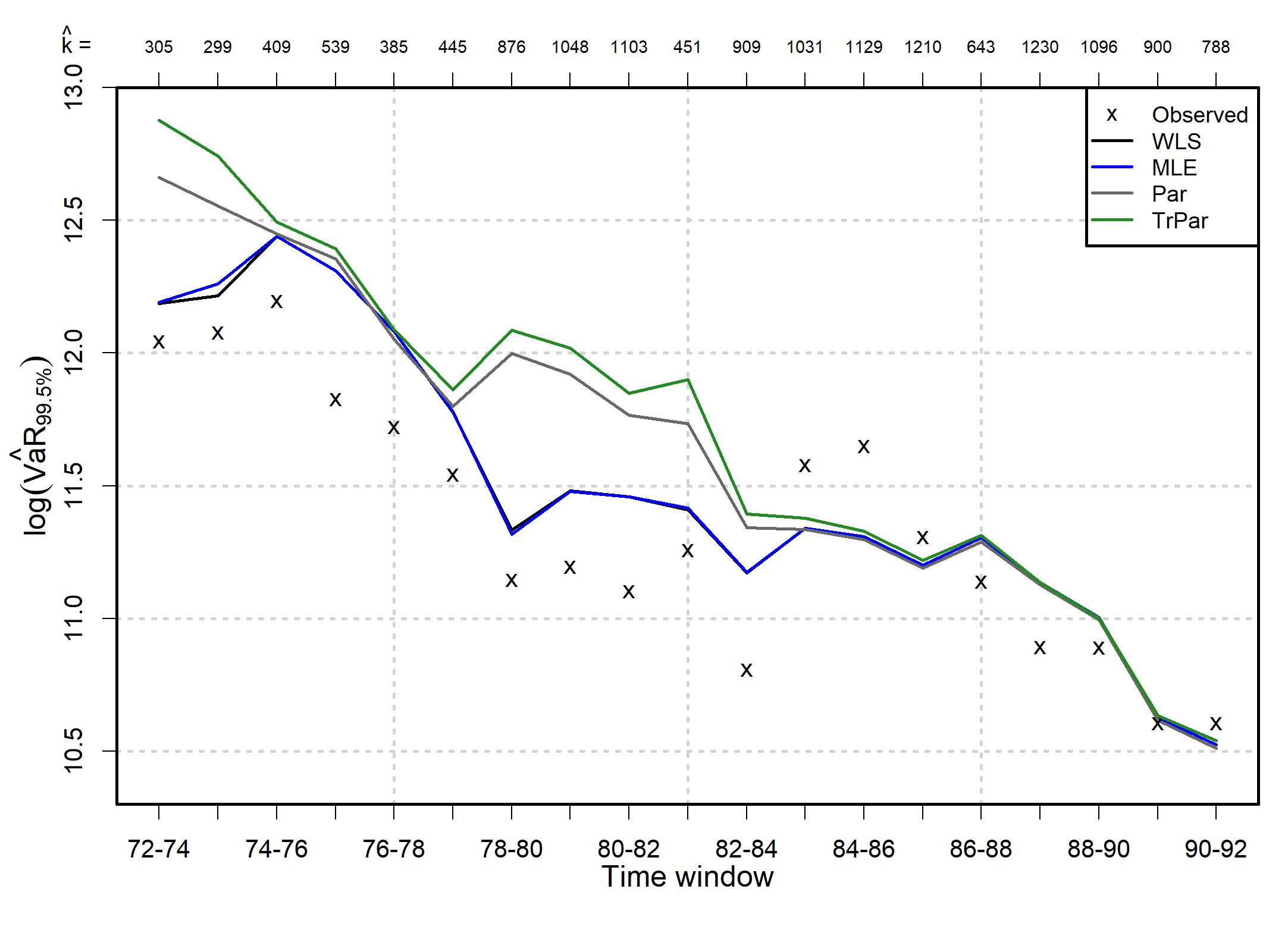}} \\
\subfloat{\includegraphics[width=0.6\textwidth]{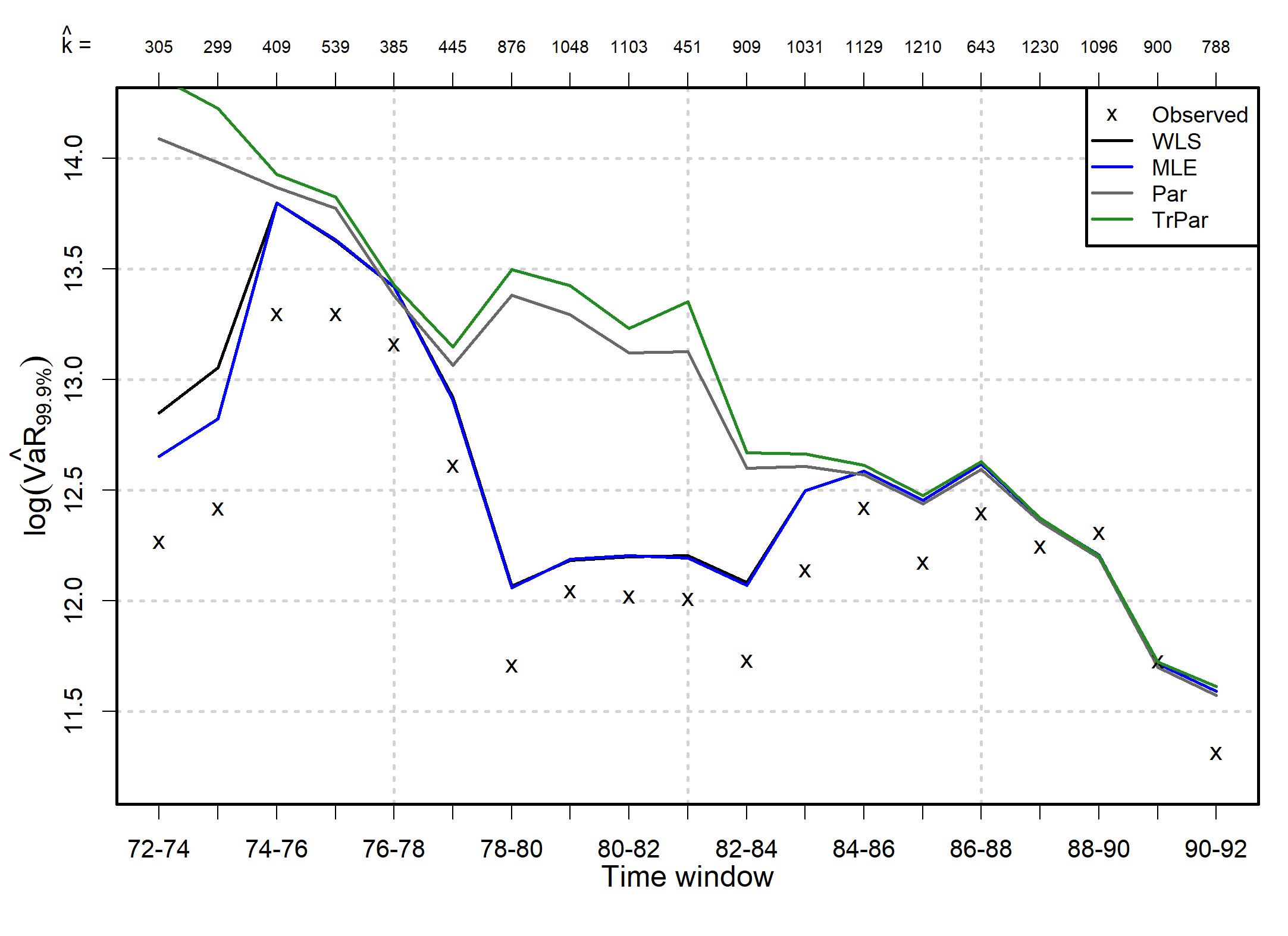}}
\caption{\textit{Norwegian fire insurance data}: log-$VaR(99.5\%)$  \textit{(top)} and log-$VaR(99.5\%)$ \textit{(bottom)} at $\hat{k}$ for tempered model (black and blue lines), Pareto (grey),  truncated Pareto (green) and observed  values (x). For each time window, $\hat{k}$ is displayed at the top margin.}
\label{fig_TD_VaR}
\end{figure}

\begin{figure}[ht]
\centering
\includegraphics[width=0.9\textwidth]{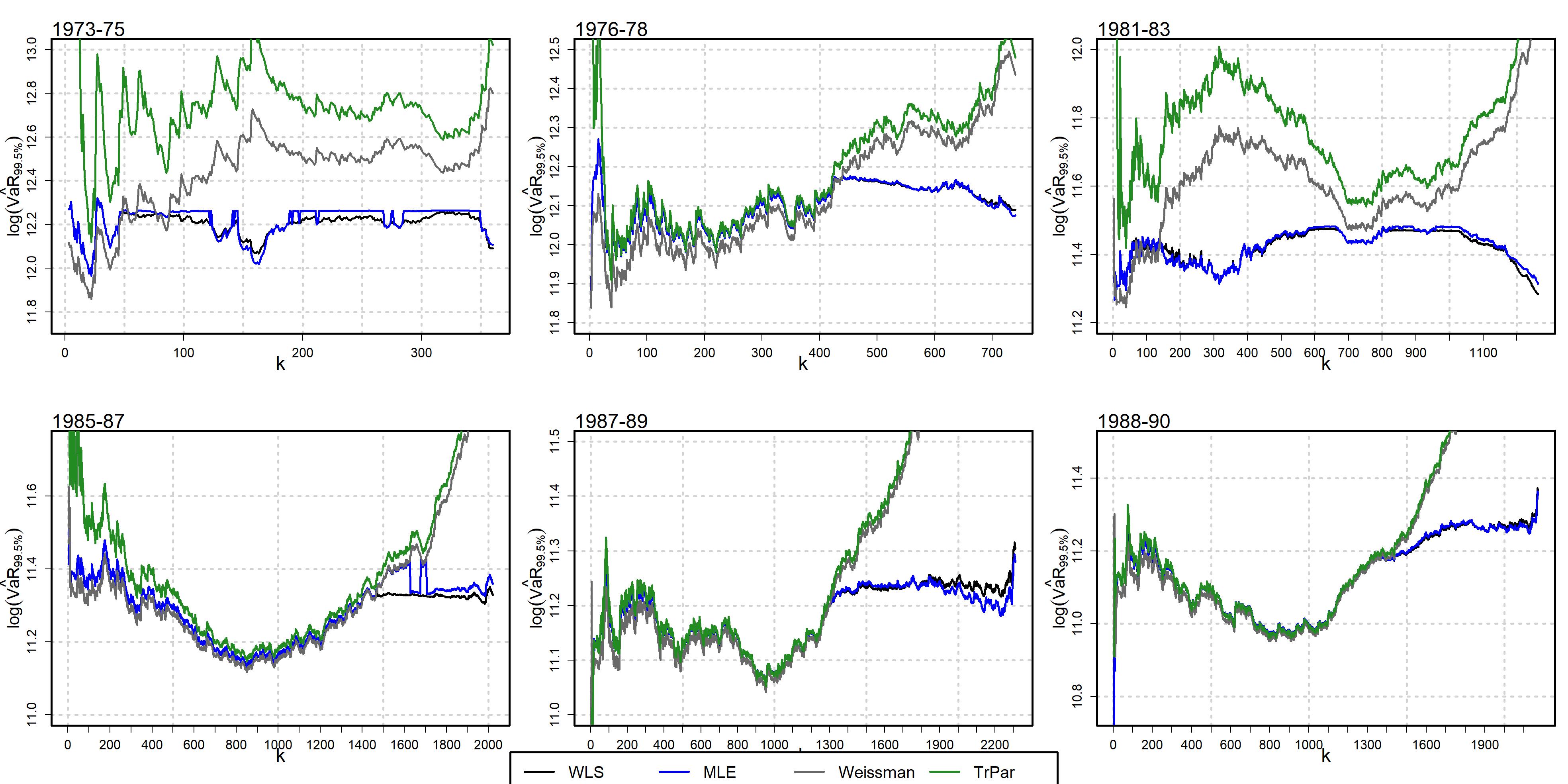}
\caption{\textit{Norwegian fire insurance data}: log $VaR(99.5\%)$ for tempered model (black and blue lines), Pareto (grey) and truncated Pareto (green) for selected time windows.}
\label{fig_TD_VaR_995}
\end{figure}

\begin{figure}[ht]
\centering
\includegraphics[width=0.9\textwidth]{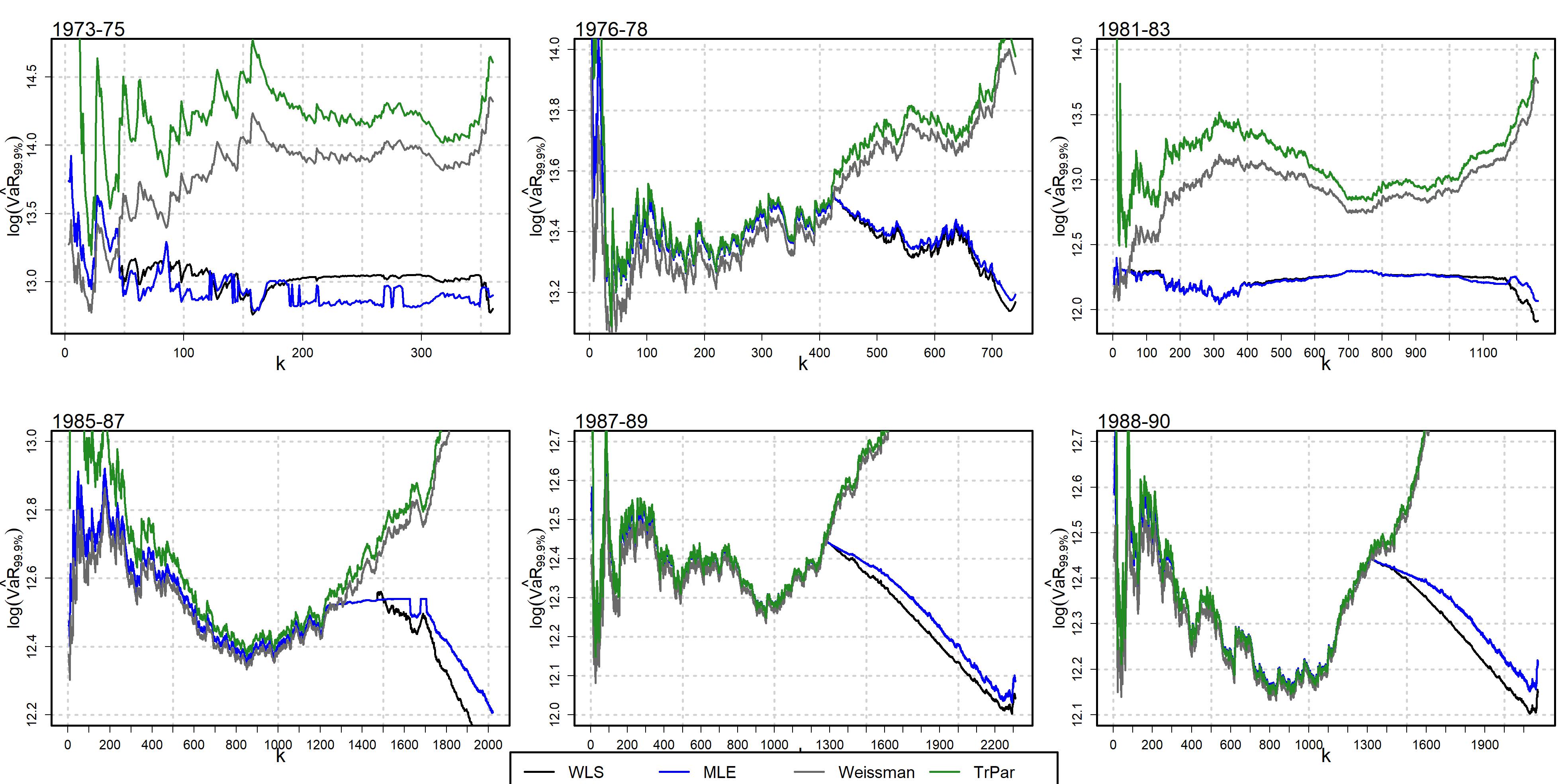}
\caption{\textit{Norwegian fire insurance data}: log $VaR(99.9\%)$ for tempered model (black and blue lines), Pareto (grey), truncated Pareto (green) for selected time windows.}
\label{fig_TD_VaR_999}
\end{figure}

In order to illustrate the possibility of extending the proposed method in a time-dependent regression context, we fitted the approach to  three-years sliding time windows. The size of the windows was selected to have at least 300 observations at each point in time. Figure \ref{fig_TD_VaR} shows the estimated $VaR$ at  99.5\% (top) and 99.9\% (bottom)  using the tempered Pareto approach with $\hat{k}$ selected using the proposed adaptive procedure, next to simple Pareto  and truncated Pareto modelling. We also compare  with the observed quantiles obtained using the standard \texttt{R} function, which estimates the quantiles as weighted averages of consecutive order statistics. The $VaR$ values based on the tempered Pareto model are situated between the observed and the Pareto and truncated Pareto fits, from which one can conclude that the tempered tail behaviour observed for the complete data set in the bottom frame in Figure \ref{Norwegian_Estimates} is also present conditional on a time window, leading to overestimation when using classical methods that ignore the proposed tempering. It is also worth noticing that  the $VaR$ at 99.5\%  values exhibit an overall decreasing trend with some stable behaviour between 1979 and 1987. Figures \ref{fig_TD_VaR_995} and \ref{fig_TD_VaR_999} show the respective $VaR$ estimates for all values of $k$ for some selected time windows.\\

\begin{figure}[ht]
\centering
\subfloat{\includegraphics[width=0.5\textwidth]{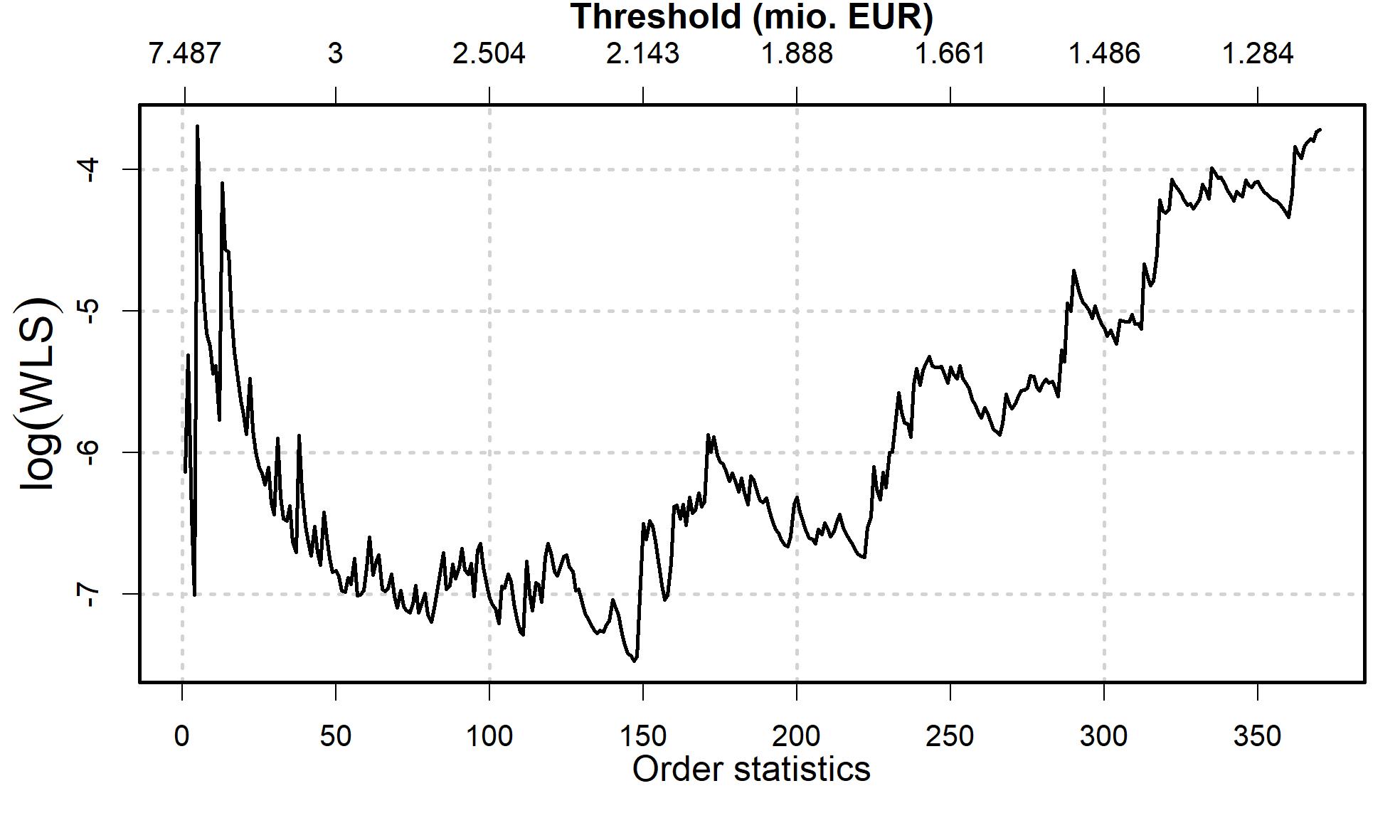}} 
	\subfloat{\includegraphics[width=0.5\textwidth]{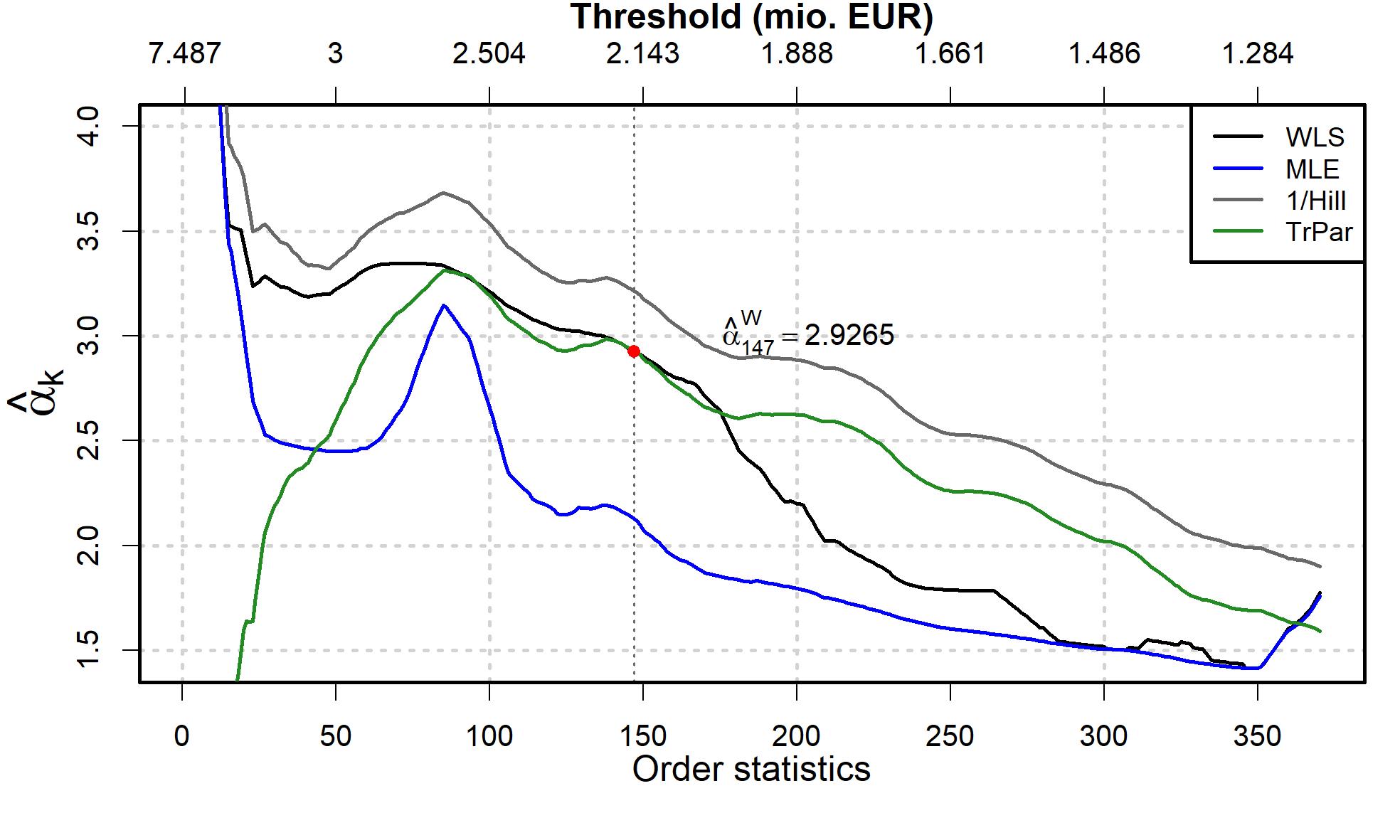}} \\
	\subfloat{\includegraphics[width=0.5\textwidth]{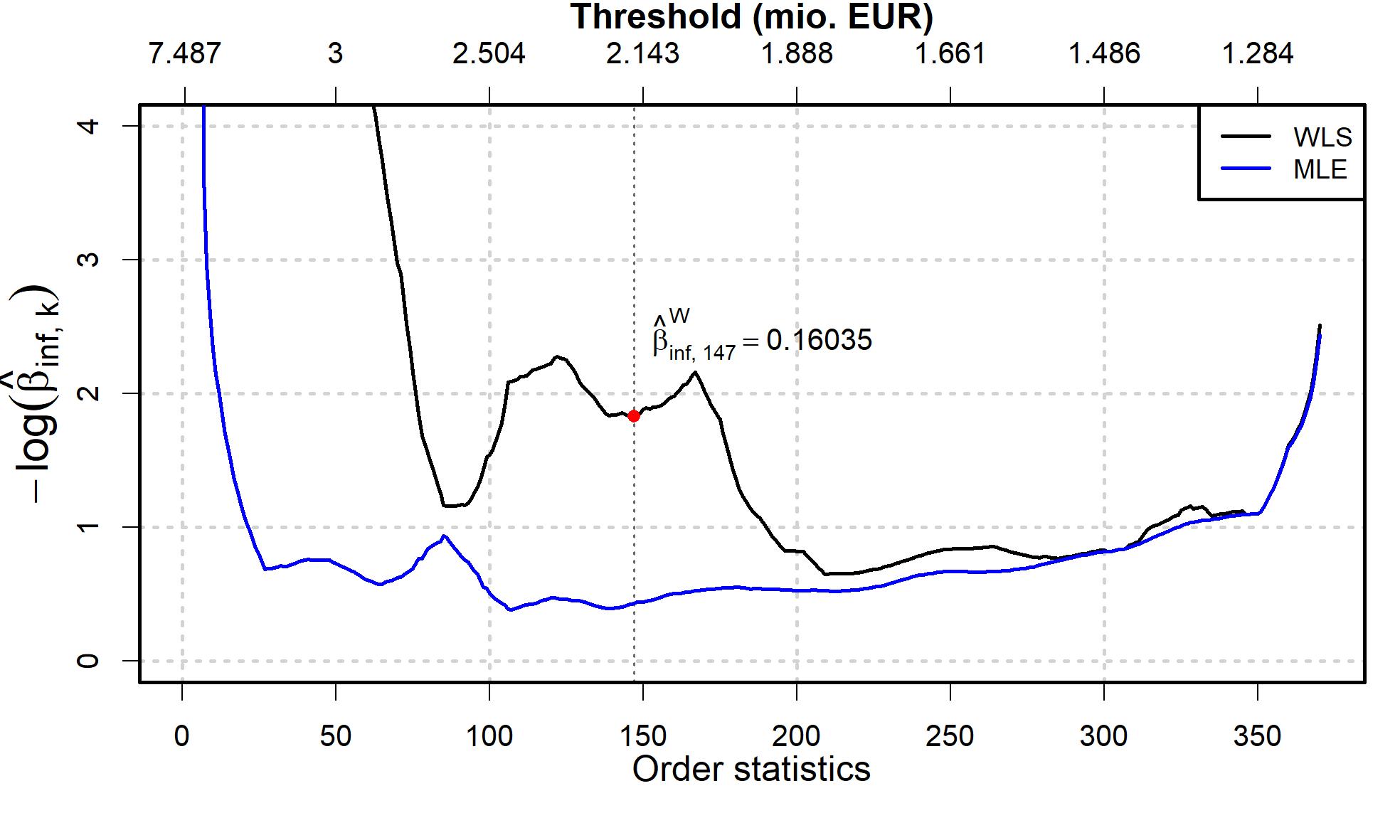}} 
	\subfloat{\includegraphics[width=0.5\textwidth]{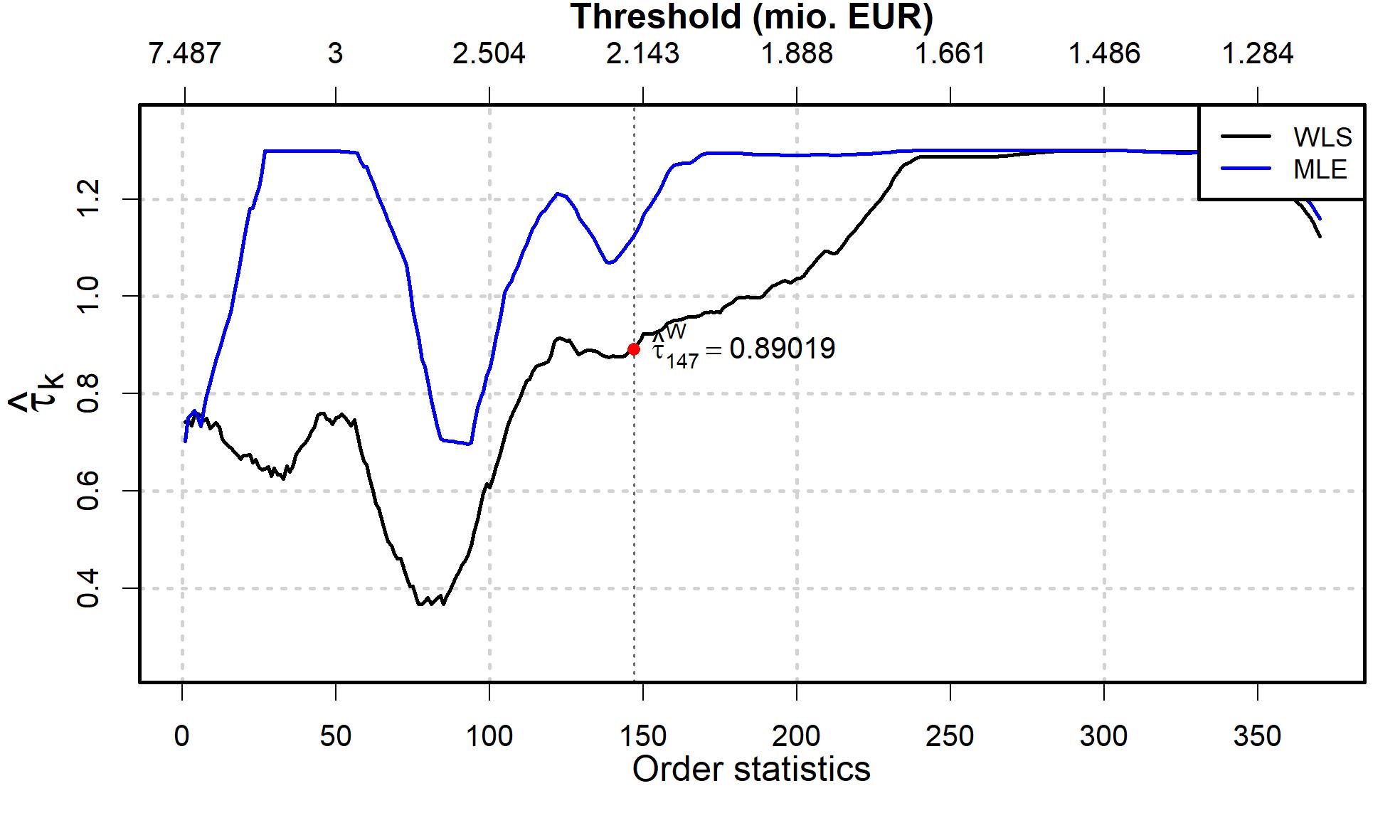}} \\
	\centering
	\subfloat{\includegraphics[width=0.8\textwidth]{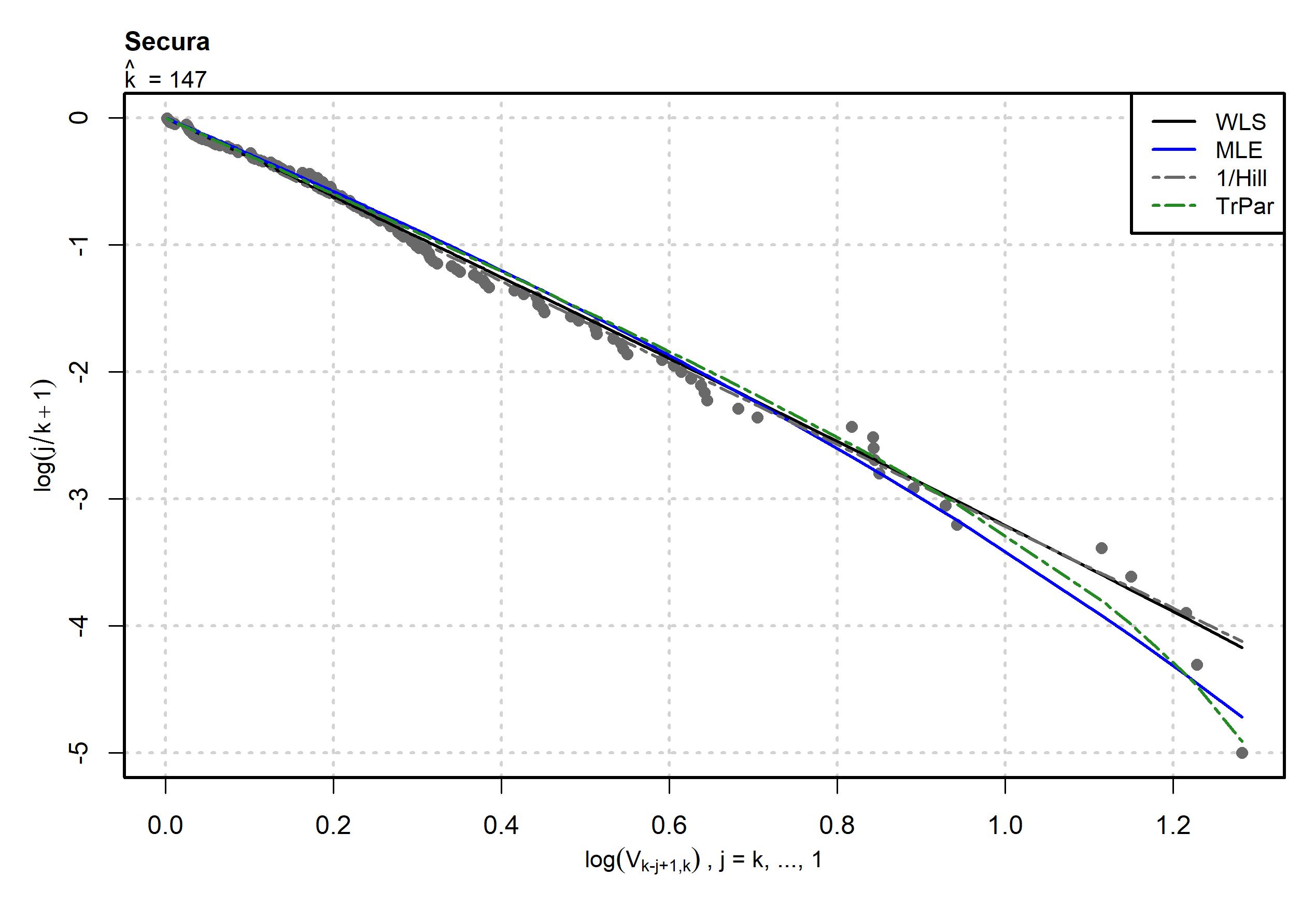}}
	\caption{\textit{Secura data set:} Top left: $SS_k$ from \eqref{SSk}; Top right: $\hat\alpha _k^W$, $\hat\alpha _k^M$, $H_{k,n}$ and $\hat{\alpha}_k^T$; 
	Middle left: $-\log\hat{\beta}_{\infty,k}^W$, $-\log\hat{\beta}_{\infty,k}^M$; Middle right: $\hat\tau _k^W$, $\hat\tau _k^M$; Bottom: $\log$-$\log$ plot with fit obtained from \eqref{eqline} with $k=\hat{k}=147$ using MLE and WLS estimates, next to Pareto and truncated Pareto fit.}
	\label{Besecura_estimates}
\end{figure}

In Figure \ref{Besecura_estimates}, the respective results are given for the Secura Re Belgium data set. Here the best tempered Pareto fit is found at $\hat{k}= 147$, with  the corresponding $\log$-$\log$ plot given in the bottom figure. Here the tempered Pareto WLS fit closely follows the linear Pareto fit, while the MLE fit shows too much bending near the largest data. Both the Pareto and WLS tempered Pareto fit do miss the deviation
at the top two data, which however is taken into account in the truncated Pareto analysis with $\hat{T}_{\hat{k}}=8,967,620 = e^{16.009}.$ While this deviation can be considered as statistically non-significant, it makes sense to consider the truncated Pareto fit here since Belgian car insurance contracts do show explicit upper limits. Another motivation for a truncated model is that the extreme quantile estimates  $\hat{Q}_{1/(cn),\hat k}^M $  hardly change from $c=1$ to $c=2$, namely around the value $e^{16}$. 

\begin{figure}[ht]
\centering
		\subfloat{\includegraphics[width=0.8\textwidth]{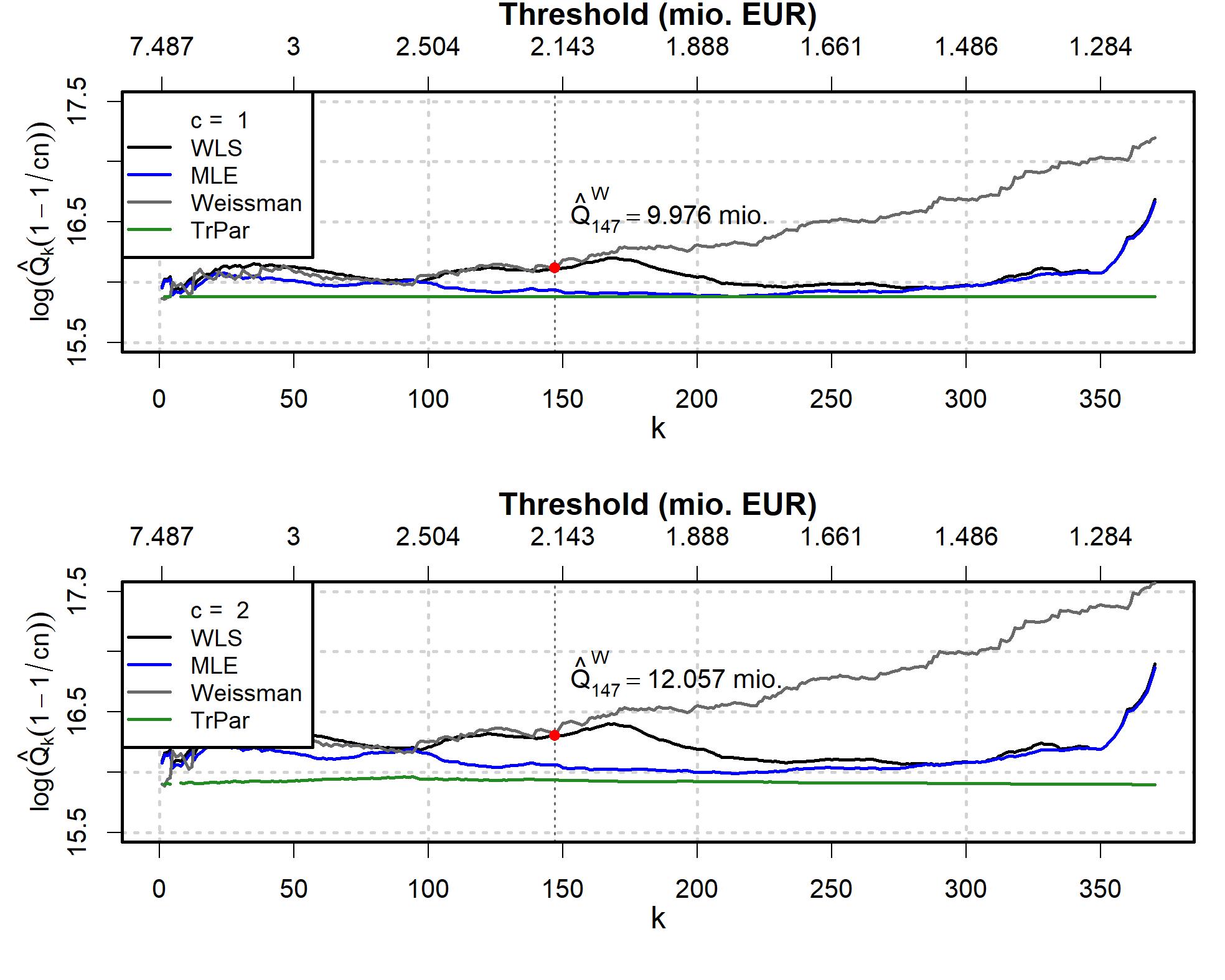}}
	\caption{\textit{Secura data set}: $\hat{Q}^W_{p,k}$, $\hat{Q}^M_{p,k}$ and $\hat{Q}^H_{p,k}$ quantile estimates with $p=1/n$ (top) and $p=1/(2n)$ (bottom).}
	\label{Besecura_quantiles}
\end{figure}


\section{Conclusion}\label{secconc}
 In this paper we addressed the fitting of Pareto-type distributions with a tempering component of Weibull type at large values. We extend earlier results for exponential tempering  on strict Pareto tails,  provide a Peaks over Threshold (POT) approach, develop  estimation procedures and provide asymptotic properties of the proposed estimators. Finally, we present a simulation study and also apply the developed methods to actual insurance data, discussing challenges in the implementation and how to overcome them. 
The estimation of $VaR$ values at extreme quantile levels shows improvements compared to more classical extreme value estimation methods that ignore the considered tempering effect. These improvements are more pronounced with growing  tempering effect. 
\\
Further research concerning the generalization to a regression context and the use of tempered Pareto-Weibull models in composed or splicing models, will be taken up in the future.  

\section{Acknowledgements}
\noindent The authors are grateful to the anonymous referees for their helpful comments and suggestions that lead to improvements of the paper. H.A. acknowledges financial support from the Swiss National Science Foundation Project 200021\_191984.


\bibliographystyle{plain}
\bibliography{References}
\addresseshere


\section{Appendix: Proof of Theorem 2.1}

Using Taylor expansions of the likelihood equations in  $\hat{\boldsymbol\theta}_t$ around the correct value $\boldsymbol\theta$ leads to the following system of three equations, with $\tilde{\boldsymbol{\theta}} =(\talpha,\tlambda,\ttau)$ situated in between $\hat{\boldsymbol\theta}_t$ and $\boldsymbol{\theta}$:
\begin{eqnarray}
&& \hspace{-2cm}
\sqrt{n \overF (t)}(\halpha_t - \alpha)\frac{1}{n\overF (t)}\sum_{j = 1}^{n} \frac{1}{\Big(\talpha + \tlambda\ttau (\Vjk)^{\ttau}\Big)^2}1_{(X_j > t)}\nonumber \\  
&&  \hspace{-1cm}+ \sqrt{n \overF(t)}(\hat\lambda _t - \lambda)\frac{1}{n\overF(t)}\sum_{j = 1}^{n} \frac{\tilde\tau (\Vjk)^{\tilde\tau}}{\Big(\tilde\alpha + \tilde\lambda\tilde\tau (\Vjk)^{\tilde\tau} \Big)^2}1_{(X_j > t)} 
\nonumber \\ && \hspace{-1cm}
+\sqrt{n \overF(t)}(\hat\tau _t - \tau)\frac{1}{n\overF(t)}\sum_{j = 1}^{n}\tilde\lambda 
\frac{(\Vjk)^{\tilde\tau}(1+\tilde\tau\log \Vjk )}{\Big(\tilde\alpha + \tilde\lambda\tilde\tau (\Vjk)^{\tilde\tau} \Big)^2} \Ij \nonumber\\
&
&  \hspace{-0.5cm}= \sqrt{n \overF(t)}\bigg(\frac{1}{n\overF(t)}\sum_{j = 1}^{n} \Big\{\frac{1}{\alpha + \lambda\tau (\Vjk)^\tau} - \log\Vjk\Big\}1_{(X_j > t)}\bigg)
\label{ML1}\\
&& \hspace{-2cm}
\sqrt{n \overF (t)}(\halpha_t - \alpha)\frac{1}{n\overF (t)}\sum_{j = 1}^{n} \frac{\tilde\tau (\Vjk)^{\tilde\tau}}{\Big(\talpha + \tlambda\ttau (\Vjk)^{\ttau}\Big)^2}1_{(X_j > t)}\nonumber \\  
&&  \hspace{-1cm}+ \sqrt{n \overF(t)}(\hat\lambda _t - \lambda)\frac{1}{n\overF(t)}\sum_{j = 1}^{n} \frac{\tilde\tau^2 (\Vjk)^{2\tilde\tau}}{\Big(\tilde\alpha + \tilde\lambda\tilde\tau (\Vjk)^{\tilde\tau} \Big)^2}1_{(X_j > t)} 
\nonumber \\ && \hspace{-1cm}
+\sqrt{n \overF(t)}(\hat\tau _t - \tau)\frac{1}{n\overF(t)}\sum_{j = 1}^{n}\Big( 
\frac{\talpha(\Vjk)^{\tilde\tau}(1+\tilde\tau\log \Vjk )}{\Big(\tilde\alpha + \tilde\lambda\tilde\tau (\Vjk)^{\tilde\tau} \Big)^2} -(\Vjk)^{\ttau} \log \Vjk
\Big) \Ij \nonumber\\
&
&  \hspace{-0.5cm}= \sqrt{n \overF(t)}\bigg(\frac{1}{n\overF(t)}\sum_{j = 1}^{n} \Big\{\frac{\tau(\Vjk)^\tau}{\alpha + \lambda\tau (\Vjk)^\tau} - (\Vjk)^\tau+1\Big\}1_{(X_j > t)}\bigg)
\label{ML2}\\
&& \hspace{-2cm}
\sqrt{n \overF (t)}(\halpha_t - \alpha)\frac{1}{n\overF (t)}\sum_{j = 1}^{n} \frac{\tlambda (\Vjk)^{\tilde\tau}
(1+\ttau\log \Vjk )}{\Big(\talpha + \tlambda\ttau (\Vjk)^{\ttau}\Big)^2}1_{(X_j > t)}\nonumber \\  
&&  \hspace{-1cm}+ \sqrt{n \overF(t)}(\hat\lambda _t - \lambda)\frac{1}{n\overF(t)}\sum_{j = 1}^{n} 
\Big( 
\frac{\talpha(\Vjk)^{\tilde\tau}(1+\tilde\tau\log \Vjk )}{\Big(\tilde\alpha + \tilde\lambda\tilde\tau (\Vjk)^{\tilde\tau} \Big)^2} -(\Vjk)^{\ttau} \log \Vjk
\Big) 1_{(X_j > t)} 
\nonumber \\ && \hspace{-1cm}
+\sqrt{n \overF(t)}(\hat\tau _t - \tau)\frac{\tlambda}{n\overF(t)}\sum_{j = 1}^{n}
\left( 
\frac{\tlambda(\Vjk)^{\tilde\tau}(1+2\tilde\tau\log \Vjk )+\talpha\ttau (\log \Vjk )^2}{\Big(\tilde\alpha + \tilde\lambda\tilde\tau (\Vjk)^{\tilde\tau} \Big)^2} \right.\nonumber\\
&&\hspace{5cm} \left.
-\frac{2\log \Vjk}{\tilde\alpha + \tilde\lambda\tilde\tau (\Vjk)^{\tilde\tau}}
+ (\log \Vjk)^2
\right) (\Vjk)^{\ttau}\Ij \nonumber\\
&
&  \hspace{-0.5cm}= \sqrt{n \overF(t)}\bigg(\frac{\lambda	}{n\overF(t)}\sum_{j = 1}^{n}
 \Big\{\frac{(\Vjk)^{\tau}(1+\tau \log\Vjk) }{\alpha + \lambda\tau (\Vjk)^\tau} - (\Vjk)^\tau  \log \Vjk \Big\}
 1_{(X_j > t)}\bigg)
\label{ML3}
\end{eqnarray}
The coefficients of $\sqrt{n \overF(t)}(\hat\alpha _t - \alpha)$, $\sqrt{n \overF(t)}(\hat\lambda _t - \lambda)$
and $\sqrt{n \overF(t)}(\hat\tau _t - \tau)$ on the left hand sides of \eqref{ML1}, \eqref{ML2} and \eqref{ML3} now converge in probability to the corresponding elements of ${\bf I}$. For instance for 
\[
{\bf I}_{1,1,n,t}(\alpha, \lambda, \tau) := 
\frac{1}{n\overF (t)}\sum_{j = 1}^{n} \frac{1}{\Big(\alpha + \lambda\tau (\Vjk)^{\tau}\Big)^2}1_{(X_j > t)}
\]
we have
\begin{eqnarray*}
\mathbb{E}({\bf I}_{1,1,n,t}(\alpha, \lambda, \tau))
 &=& -\int_t^\infty
\frac{1}{\Big(\alpha + \lambda\tau ({x \over t})^{\ttau}\Big)^2}
d {\overline{F} (x)\over \overline{F} (t)}\\
&=& -\int_1^\infty
\frac{1}{\Big(\alpha + \lambda\tau u^{\tau}\Big)^2}
d \overline{F}_t (u)\\
&\to&  -\int_1^\infty
\frac{1}{\Big(\alpha + \lambda\tau u^{\tau}\Big)^2}
d \overline{F} _{\alpha,\lambda,\tau}(u)=I_{1,1},
\end{eqnarray*}
as $t \to \infty$ using the consistency of ML estimators and assumption $(\mathcal{M})$. The convergence of 
${\bf I}_{1,1,n,t}(\talpha, \tlambda, \ttau)$ to ${\bf I}_{1,1}$ then follows from 
\[
\mbox{Var}\left(
{\bf I}_{1,1,n,t}(\alpha, \lambda, \tau)
\right) = O\left((n\overline{F}(t))^{-1}\right) \mbox{ and }
{\bf I}_{1,1,n,t}(\talpha, \tlambda, \ttau)-
{\bf I}_{1,1,n,t}(\alpha, \lambda, \tau) =o_p(1)
\] as $n,t \to \infty$ using the consistency of the ML estimators. \\

\noindent
Next the asymptotic normal distribution of  the right hand sides of \eqref{ML1}-\eqref{ML3}
\begin{eqnarray}
&& \hspace{-2.5cm}\sqrt{n\overF(t)}\left(
\frac{1}{n\overF(t)}\sum_{j = 1}^{n} \Big\{\frac{1}{\alpha + \lambda\tau (\Vjk)^\tau} - \log\Vjk\Big\}1_{(X_j > t)} \right. \nonumber\\
&&, 
\frac{1}{n\overF(t)}\sum_{j = 1}^{n} \Big\{\frac{\tau(\Vjk)^\tau}{\alpha + \lambda\tau (\Vjk)^\tau} - (\Vjk)^\tau+1\Big\}1_{(X_j > t)}\nonumber\\
&&, \left.
\frac{\lambda	}{n\overF(t)}\sum_{j = 1}^{n}
 \Big\{\frac{(\Vjk)^{\tau}(1+\tau \log\Vjk) }{\alpha + \lambda\tau (\Vjk)^\tau} - (\Vjk)^\tau  \log \Vjk \Big\}
 1_{(X_j > t)}
\right)
\label{RHS}
\end{eqnarray}
is derived.
\\
Concerning the first component
\begin{align*}
\frac{1}{n\overF(t)}\E\Bigg(\sum_{j = 1}^{n} \Big\{\frac{1}{\alpha + \lambda\tau (\Vjk)^\tau} - \log\frac{X_j}{t} \Big\}1_{(X_j > t)}\Bigg) &= -\frac{1}{\overF(t)}\int_{t}^{\infty} \Big\{\frac{1}{\alpha + \betain \frac{x}{t}} - \log\Big(\frac{x}{t}\Big)\Big\} d\overF(x) \\
&= -\int_{1}^{\infty} \Big\{\frac{1}{\alpha + \lambda\tau u^{\tau}} - \log u\Big\} d\overF_t(u),
\end{align*}
with $\overF_t(u) = \Prob(X/t > u | X>t) = u^{-\alpha}(1+D t^{\rho} h_{\rho}(u))e^{-\lambda (u^\tau-1)}$ using the second order slow variation condition \eqref{2ndorder}, so that
\begin{equation*}
-\frac{d\overF_t(u)}{du} 
= u^{-\alpha - 1}e^{-\lambda(u^\tau-1)}(\alpha + \lambda\tau u^{\tau}) + Dt^{\rho}\; u^{-\alpha - 1}
e^{-\lambda(u^\tau-1)}
\{h_{\rho}(u)[\alpha + \lambda\tau u^\tau] - u^{\rho}\}. 
\end{equation*}
Using partial integration one easily checks that
\[
\int_{1}^{\infty} \Big\{\frac{1}{\alpha + \lambda\tau u^{\tau}} - \log u\Big\}u^{-\alpha - 1}e^{-\lambda(u^\tau-1)}(\alpha + \lambda\tau u^{\tau})du=0,
\]
so that the expected value of the first component is given by $Dt^\rho \, b_1$, leading to the asymptotic bias expression of $\hat{\alpha}_t$ as given in Theorem 2.1, and similar calculations lead to the bias of $\hat\lambda_t$ and $\hat\tau_t$.
\\

\noindent 
So it remains to derive the asymptotic variances and covariances of the vector in \eqref{RHS}. The variance of the first component is derived  from
\begin{eqnarray*}
&& \hspace{-1cm}
\frac{1}{n\overF(t)}\sum_{j = 1}^{n}\mathbb{E} \Big\{\frac{1}{\alpha + \lambda\tau (\Vjk)^\tau} - \log\Vjk\Big\}^2 1_{(X_j > t)} \\&=&
\frac{1}{\overF(t)}\mathbb{E}\left(
 \Big\{\frac{1}{\alpha + \lambda\tau ({X \over t})^\tau} - \log {X \over t}\Big\}^2 1_{(X > t)}
 \right)\\
 &=&-\int_1^\infty \left( 
\frac{1}{(\alpha + \lambda\tau u^\tau )^2} 
 - \frac{2\log u}{\alpha + \lambda\tau u^\tau }
 + (\log u)^2
 \right) d\overline{F}_t (u)
\\
&\to & 
 -\int_1^\infty \left( 
\frac{1}{(\alpha + \lambda\tau u^\tau )^2} 
 - \frac{2\log u}{\alpha + \lambda\tau u^\tau }
 + (\log u)^2
 \right) du^{-\alpha}e^{-\lambda (u^\tau-1)},
\end{eqnarray*} 
as $n,t \to \infty$.
Using partial integration one finds that
$\int_1^\infty (\frac{2\log u}{\alpha + \lambda\tau u^\tau }
 - (\log u)^2)du^{-\alpha}e^{-\lambda (u^\tau-1)}=0,$
 so that the asymptotic variance of the first component in \eqref{RHS} equals ${\bf I}_{1,1}$.
 In the same way one finds that the asymptotic variance covariance matrix of \eqref{RHS} equals ${\bf I}$. \\
 
 \noindent
 Hence
 \begin{equation}
 \left( {\bf I} + o_p(1)\right)\sqrt{n\overline{F}(t)}(\hat{\boldsymbol{\theta}}_t -
 \boldsymbol{\theta}) = \mathcal{N}_3 \left((D\nu){\bf b},{\bf I}\right)+ o_p(1),
 \end{equation}
from which the result follows.

\end{document}